\documentclass[table]{amsart}
\usepackage[utf8]{inputenc}
\usepackage[top=1in,left=1in,right=1in,bottom=1in]{geometry}
\usepackage{booktabs}
\usepackage{todonotes}
\usepackage{longtable}

\usepackage{amsmath}
\usepackage{amssymb}
\usepackage{amsthm}
\usepackage{amsfonts,mathrsfs}
\usepackage{amscd}
\usepackage{array}
\usepackage{bigstrut}
\usepackage{color}
\usepackage{listings}
\usetikzlibrary{decorations.pathreplacing}
\usetikzlibrary{matrix,arrows.meta}
\usepackage[center]{caption}
\usepackage{enumitem}
\usepackage{pgfplots}
\usepackage{subfigure}
\usepackage{hyperref}
\usepackage{calligra}
\usepackage[T1]{fontenc}
\usepackage{flafter}
\usepackage{relsize}
\usepackage{mathtools}

\definecolor{TableShade}{rgb}{0.88,1,1}

\usepackage{mathrsfs}

\usepackage{subfigure}
\usepackage{wrapfig}
\usepackage{tikz-cd}
\usepackage{graphicx}
\usepackage{caption}
\newcolumntype{L}{>{\centering\arraybackslash}p{3cm}}

\newtheorem{theorem}{Theorem}[section]
\newtheorem*{theorem*}{Theorem}
\newtheorem{lemma}[theorem]{Lemma}

\newtheorem{condition}[theorem]{Condition}
\newtheorem*{result}{Main Result}

\newtheorem{corollary}[theorem]{Corollary}
\newtheorem*{corollary*}{Corollary}
\newtheorem{proposition}[theorem]{Proposition}

\theoremstyle{definition}
\newtheorem{definition}[theorem]{Definition}

\theoremstyle{remark}
\newtheorem*{remark}{Remark}
\newtheorem{example}[theorem]{Example}

\newtheorem{convention}[theorem]{Convention}

\newcommand{\CC}{\mathbb{C}}

\newcommand{\ZZ}{\mathbb{Z}}
\newcommand{\QQ}{\mathbb{Q}}

\newcommand{\HH}{\mathcal{H}} % state space
\newcommand{\cQ}{\mathcal{Q}} % curly Q

\newcommand{\cK}{\mathcal{K}}

\newcommand{\cO}{\mathcal{O}}

\newcommand{\Ba}{\mathsf B}
\newcommand{\Ca}{\mathsf C}
\newcommand{\Ra}{\mathsf R}
\newcommand{\Ea}{\mathsf E}

\newcommand{\orbit}[1]{\mathcal{O}_{#1}}

\DeclareMathOperator{\Hom}{Hom}  % Hom
\DeclareMathOperator{\SL}{SL}   % SL: special linear group

\DeclareMathOperator{\sign}{sign}
\DeclareMathOperator{\GL}{GL}
\DeclareMathOperator{\Fix}{Fix}

\DeclareMathOperator{\diag}{diag}
\DeclareMathOperator{\age}{age}

\DeclareMathOperator{\id}{id}
\DeclareMathOperator{\Span}{Span}
\DeclareMathOperator{\conj}{conj}
\DeclareMathOperator{\bideg}{bideg}
\DeclareMathOperator{\Stab}{Stab}

\usepackage{epsfig}

\usepackage{xcolor} %this is so we can use colored fonts

\makeatletter
\def\imod#1{\allowbreak\mkern10mu({\operator@font mod}\,\,#1)}
\makeatother

\title{Mirror Map for Landau-Ginzburg models with nonabelian groups}
\author{Annabelle Clawson}
\author{Drew Johnson}
\author{Duncan Morais}
\author{Nathan Priddis}
\author{Caroline B. White}
\date{\today}

\begin{document}

%Let \(a, b, n\) be integer with \(n>0\).
\maketitle

\begin{abstract}
	BHK mirror symmetry as introduced by Berglund--H\"ubsch and Marc Krawitz between Landau--Ginzburg (LG) models has been the topic of much study in recent years. An LG model is determined by a potential function and a group of symmetries. BHK mirror symmetry is only valid when the group of symmetries is comprised of the so-called diagonal symmetries. 
	Recently, an extension to BHK mirror symmetry to include nonabelian symmetry groups has been conjectured.  
	In this article, we provide a mirror map at the level of state spaces between the LG A-model state space and the LG B-model state space for the mirror model predicted by the BHK mirror symmetry extension for nonabelian LG models. 
	We introduce two technical conditions, the Diagonal Scaling Condition, and the Equivariant $\Phi$ condition, under which a bi-degree preserving isomorphism of state spaces (the mirror map) is guaranteed to exist, and we prove that the condition is always satisfied if the permutation part of the group is cyclic of prime order. 
\end{abstract}

\tableofcontents

\section{Introduction}

Mirror symmetry is a phenomenon that was initially discovered by theoretical physicists, and has since been the vehicle for several major breakthroughs in mathematics. 

Understood at its most basic level, mirror symmetry is an equivalence between the variation of K\"ahler structures (the A-model) of a Calabi--Yau manifold $X$ and the variation of complex structures (the B-model) on the \emph{mirror} Calabi--Yau $\check{X}$. In particular, one first prediction of mirror symmetry is that the Hodge structure of $\check{X}$ can be obtained by rotating the Hodge diamond of $X$ by ninety degrees, i.e., $h^{2,1}(X)=h^{1,1}(\check{X})$ and $h^{2,1}(\check{X})=h^{1,1}(X)$. 

The physics can also be described using what is known as a Landau--Ginzburg (LG) model. The predictions of mirror symmetry for LG models are not as easily given a geometric meaning, but there is also an LG A-model and an LG B-model and a similar mirror symmetry statement which is predicted for these LG models. 

An LG model is determined by a pair $(W,G)$, where $W$ is a potential and $G$ is a group of symmetries. To each such pair, one can associate what is known as a \emph{state space} with a bidegree---much like the cohomology of a CY manifold and its associated Hodge diamond. At the most basic level, mirror symmetry for LG models can be formulated as an isomorphism between the state space for the LG A-model of a pair $(W,G)$ and the state space for the LG B-model of the \emph{mirror pair} $(W^\vee, G^\vee)$. The purpose of this article is to establish this isomorphism, especially in the case in which the group $G$ is not abelian. 

This isomorphism has already been established in \cite{Krawitz} in the case that the group $G$ is made up of by what are known as \emph{diagonal symmetries} (see Section~\ref{sec:ms-diagonal} for a definition). Aside from establishing the mirror map in the diagonal case, one of the most important contributions of Krawitz in \cite{Krawitz}, and independently Berglund and H\"ubsch in \cite{BH} and Berglund and Henningson in \cite{BergHenn}, was determining what the mirror model should be. In other words, given a pair $(W,G)$, can we say what the mirror pair $(W^\vee, G^\vee)$ should be? This rule has come to be known as BHK mirror symmetry (for diagonal LG models). We will review this construction in Section~\ref{sec:ms-diagonal}.

Until recently, only the diagonal case has been investigated, because it was unknown how to find the mirror model $(W^\vee,G^\vee)$ when $G$ was not a group of diagonal symmetries. A few years ago, however, a paper by Mukai \cite{M} gave many examples of the dimensions of the various state spaces. Soon thereafter, several authors conjectured a rule for finding the mirror $(W^\vee, G^\vee)$ for certain nonabelian groups $G$. In \cite{Ebel1}, Ebeling and Gusein-Zade attributed the conjecture to A. Takahashi and then gave some evidence that this conjecture is correct based on the dimensions of what are known as Milnor fibers. Independently, in \cite{PWW}, Priddis, Ward, and Williams made the same conjecture. and then gave an isomorphism between certain subspaces of the A-model and B-model state spaces when the polynomial $W$ is of Fermat type. Furthermore, in \cite{basalaev}, Basalaev and Ionov provide a mirror map between certain subspaces of the respective state spaces for the case when $W$ is a Fermat polynomial and $H$ is either maximal or a certain subgroup known as $\SL_W$. In that paper, the authors also provide a sate space isomorphism when $W$ is the Fermat quintic, but for \emph{any} subgroup $S\subseteq S_N$. And finally, in \cite{M2} Mukai has obtained a mirror map for certain of the so-called loop potentials when $S$ is a group of prime order.

In the current paper, we build on these results, giving an isomorphism between the LG A-model state space for a pair $(W,G)$ and the LG B-model state space for the \emph{mirror pair} $(W^\vee,G^\vee)$ according to the rule described in \cite{Ebel1} and \cite{PWW}. This rule is defined only for symmetry groups of the form $G=S\rtimes H$, where $H$ is a group of so-called \emph{diagonal} symmetries, and $S$ is a group permuting the variables of $W$ (see Section~\ref{sec:group}). 

Unfortunately, it turns out that the two state spaces are not always isomorphic, as pointed out in the papers \cite{Ebel1}, \cite{EG2}, 
and \cite{PWW}. In fact, in \cite{EG2}, 
Ebeling and Gusein-Zade introduce a condition called the \emph{Parity Condition} (PC), which is a condition on the group $S$. They conjecture that PC is necessary in order to have a mirror isomorphism. In the current article, we give different (somewhat technical) conditions that ensure that such an isomorphism exists (see Sections~\ref{sec:strongPC} and \ref{sec:EV} and the examples of Section~\ref{sec:examples}). It is not clear how exactly our conditions are related to PC. And we conjecture that these conditions are also necessary.

Our main result is as follows:
\begin{result}
	Let $(W,G)$ be a pair consisting of an invertible polynomial $W$ and a group $G$ of the form $G=S\rtimes H$ such that $S$ satisfies Condition~\ref{EV}. Then there is a bidegree-preserving isomorphism between the state space of the LG A-model for $(W,G)$ and the LG B-model state space for the mirror pair $(W^\vee,G^\vee)$. 
\end{result}

In Section~\ref{sec:prime-order}, we show that these conditions are satisfied in many cases under a condition that is much easier to verify. 
\begin{corollary*}
	Let $(W,G)$ be a pair consisting of an invertible polynomial $W$ and a group $G=S\rtimes H$ with $S$ generated by a permutation of odd prime order. Then there is a bidegree-preserving isomorphism between the state space of the LG A-model for $(W,G)$ and the LG B-model state space for the mirror pair $(W^\vee,G^\vee)$. 
\end{corollary*}

Another interesting result in this direction was obtained by Filipazzi and Rota in \cite{FR}, where they give an example of an isomorphism of state spaces for complete intersections, rather than strictly for hypersurfaces, as we consider in this article. There is also a vast generalization of BHK mirror symmetry described by Berglund and H\"ubsch in \cite{BH2} and \cite{BH3}. It would be interesting to see if these constructions can be extended to the non-abelian case as well. Also in related work to appear in Advances in Math, in \cite{GRZ} Gr\"afnitz, Ruddat, and Zaslow obtain an open mirror map for del Pezzo surfaces via wall-crossing of the Hori--Vafa LG potential.

Finally, we would like to mention that in the diagonal case, much more is known beyond an isomorphism of state spaces. Both the LG A-model and LG B-model state spaces can be given the structure of a Frobenius algebra (see \cite{FJR},  \cite{KPABR}, or \cite{FJJS} for the A-model, and \cite{basalaev2} for the B-model), and in \cite{FJJS}, it was shown that for a large class of LG models, the two state spaces are also isomorphic as Frobenius algebras. This is conjectured to be true in general, but the general case remains open. And it turns out there is more structure to both models. See for example \cite{CIR:11}, \cite{ChiRu:09}, \cite{FJRW}, \cite{Guere}, \cite{HLSW}, and \cite{PrSh:13}.

The paper is organized as follows: In Section~\ref{sec:poly}, we will describe the classes of polynomials $W$ and symmetry groups $G$ that we will consider for mirror symmetry. In Section~\ref{sec:state-spaces}, we will define the building blocks of the respective state spaces, and give their definitions. In Section~\ref{sec:ms-diagonal}, we will describe BHK mirror symmetry for the case when $G$ is a group of diagonal symmetries. In Section~\ref{sec:ms-non-diag}, we give the rule for finding the mirror in the case that $G$ is nonabelian. We give an explicit example illustrating the isomorphism between the respective state spaces, and then describe the isomorphism in general and prove that such an isomorphism exists when $W$ and $G$ satisfy certain conditions. These conditions are quite technical, it turns out, so in Section~\ref{sec:prime-order}, we will give several stronger conditions which ensure that the isomorphism exists. Section~\ref{sec:examples} gives some examples and nonexamples that illustrate why the conditions are necessary.

Finally, let us mention that three of the authors of this article were undergraduates at the time of this writing. The discovery of the isomorphism is actually due to them. As such, this article is aimed at the level of advanced undergraduates. It is heavy with notation, but we provide many examples and have tried to make the paper as self-contained as possible. The aim is to hopefully ease the work of the reader and clarify the exposition. The cost is that the article turned out a bit longer than expected. 

Acknowledgements: The authors would like to thank Brigham Young University for its emphasis on undergraduate research, and for supporting the A.C., D.M., and C.B.W. so generously with space, resources, and funding during this project. The third author was supported by Collaboration Grant 586691 from the Simons Foundation.

\section{Invertible Polynomials and their Symmetries} \label{sec:poly}

The LG model state space is constructed from a potential function $W$ and a group $G$ of symmetries of $W$. In this section, we specify the class of functions and groups that are candidates for BHK mirror symmetry.

\subsection{The Polynomial \texorpdfstring{$W$}{W}}

\begin{definition} We say that a polynomial $W: \CC^N\rightarrow\CC$ is \textit{quasihomogeneous} if there exist rational numbers $q_1,\dots, q_N$
such that for every $c\in\CC^*$, we can write
\begin{align*}
    W(c^{q_1}x_1,\dots,c^{q_N}x_N)=cW(x_1,\dots,x_N)
\end{align*}
We call the numbers $q_1,\dots, q_N$ the \emph{weights}.

Furthermore, we say that a polynomial $W:\CC^N\rightarrow\CC$ is \textit{nondegenerate} if the rational weights $q_1,\dots,q_N$ of $W$ are uniquely determined and if $W$ has an isolated critical point at the origin.
\end{definition}

To each pair $(W,G)$, we can associate a geometric object, which is the quotient by $G$ of the hypersurface in weighted projective space defined by $W=0$. There is a condition on the pair $(W,G)$ called the Calabi--Yau condition (see, e.g., \cite{ChiRu:09}) so that when the pair $(W,G)$ satisfies this condition, then the corresponding geometric object will in fact be a Calabi--Yau orbifold. It has been shown in many cases that mirror symmetry for the pair $(W,G)$ agrees with the mirror symmetry we expect for the corresponding Calabi--Yau model (see, e.g., \cite{CIR:11}, \cite{PrSh:13}, \cite{PLS:14}). However, mirror symmetry for LG models does not require that $(W,G)$ satisfy any sort of Calabi--Yau type condition. In our generalization of the mirror map for LG models, we will also not need any such condition. 

\begin{definition} For a polynomial $W(x_1,\dots,x_N)$, we say that $W$ is \textit{invertible} if it satisfies the following conditions: 
\begin{enumerate}
    \item $W$ is nondegenerate and quasihomogeneous.
    \item $W$ has the same number of monomials as it has variables.
\end{enumerate}
\end{definition}

The following theorem, due to Kreuzer and Skarke, gives us a nice classification of invertible polynomials. 

\begin{theorem}[Kreuzer and Skarke \cite{KS}] \label{lem:kreu-ska} Every invertible polynomial is a disjoint sum of the following \textit{atomic types}:
\begin{enumerate}
\item A \textit{Fermat} type polynomial of the form $x_1^{a_1}$.
\item A \textit{loop} type polynomial of the form $x_1^{a_1}x_2+x_2^{a_2}x_3+\dots+x_N^{a_N}x_1$.
\item A \textit{chain} type polynomial of the form $x_1^{a_1}x_2+x_2^{a_2}x_3+\dots+x_{N-1}^{a_{N-1}}x_N+x_N^{a_N}$.
\end{enumerate}
In all cases, we require $a_i\geq 2$.
\end{theorem}

The term ``disjoint'' in this context means that the various types making up an invertible polynomial share no common variables. 

Because we can scale away any coefficient, an invertible polynomial is uniquely determined by the exponents of each monomial. Given an invertible polynomial $W$, we organize the exponents into the \emph{exponent matrix} $A_W$. Each row corresponds to a single monomial, and each column corresponds to one of the variables. 

By Lemma \ref{lem:kreu-ska}, every monomial of a nondegenerate invertible polynomial is of the form $x_i^a x_j$ or $x_i^a$. When forming the exponent matrix, it is useful to take the following important convention.

\begin{convention}
We order the monomials so that the $i$-th monomial is of the form $x_i^ax_j$ or $x_i^a$. This fact will be critical later. 
\end{convention} 

\begin{example}
Consider the polynomial $W=x_1^4+x_2^4+x_3^4+x_4^6$. It is a sum of Fermat type polynomials with weights $(\frac{1}{4},\frac{1}{4},\frac{1}{4}, \frac{1}{6})$.
Clearly, $W$ is quasihomogeneous and nondegenerate. The exponent matrix is 
\[
A_W=\begin{bmatrix}
4 & 0 & 0 & 0\\
0 & 4 & 0 & 0\\
0 & 0 & 4 & 0\\
0 & 0 & 0 & 6
\end{bmatrix}.\]

\end{example}

\begin{example}
The polynomial $W=x_1^4x_2+x_2^5x_3+x_3^3x_4+x_4^2$ is a chain potential. The weights for this polynomial are $(\frac{5}{24}, \frac{1}{6}, \frac{1}{6}, \frac{1}{2})$. The exponent matrix is
\[
\begin{bmatrix}
4 & 1 & 0 & 0\\
0 & 5 & 1 & 0\\
0 & 0 & 3 & 1\\
0 & 0 & 0 & 2
\end{bmatrix}.
\]
On the other hand, the polynomial $W=x_1^3x_2+x_2^3x_3+x_3^3x_4+x_4^3x_5+x_5^3x_6+x_6^3x_1$ is a loop potential. The weights for this polynomial are $(\frac{1}{4}, \frac{1}{4}, \frac{1}{4}, \frac{1}{4}, \frac{1}{4}, \frac{1}{4})$. The following matrix is the exponent matrix:
\[
\begin{bmatrix}
3 & 1 & 0 & 0 & 0 & 0 \\
0 & 3 & 1 & 0 & 0 & 0\\
0 & 0 & 3 & 1 & 0 & 0\\
0 & 0 & 0 & 3 & 1 & 0\\
0 & 0 & 0 & 0 & 3 & 1\\
1 & 0 & 0 & 0 & 0 & 3
\end{bmatrix}.
\]
\end{example}

\begin{remark} \label{rem:weightsAi}
The vector of weights $q_i$ of an invertible polynomial can be computed as $A_W^{-1} \mathbf 1$, where $\mathbf 1$ is a column vector of ones. This is equivalent to adding the rows of the inverse matrix. 
\end{remark}

\subsection{The Group \texorpdfstring{$G$}{G}}\label{sec:group}

Next we define the maximal symmetry group $G_W^{max}$ and two important subgroups.

\begin{definition}[Mukai \cite{M}] Let $W: \CC^N \rightarrow \CC$ be an invertible polynomial with weights $(q_1,\dots,q_N)$. Then the \textit{maximal symmetry group of W}, denoted $G_W^{max}$, is defined as follows:
\begin{equation}\label{def:Gmax}
G_W^{max}:=\{g\in \GL_N(\CC):W(g\cdot(x_1,\dots,x_N))=W(x_1,\dots,x_N)\text{ and } g_{ij}=0\text{ if } q_i\neq q_j\}.
\end{equation}

\end{definition}

\begin{definition}
Let $W(x_1,\dots,x_N)$ be a polynomial and $g\in G^{max}_{W}$. We define the \textit{fixed locus} of $g$
\[
\Fix(g)=\{\mathbf{x} \in \CC^N : g \cdot \mathbf{x}=\mathbf{x}\}.
\]
If we consider $g$ as a matrix, this is the eigenspace of $g$ with eigenvalue 1. The dimension of $\Fix(g)$ is denoted $N_g=\dim_\CC\Fix (g)$.
We will denote $W_g:=W|_{\Fix(g)}$. 
If $\Fix(g)=\{(0,\dots,0)\}$ then we say $g$ is \emph{narrow}. Otherwise, $g$ is \emph{broad}. 
\end{definition}

\begin{definition}[Mukai \cite{M}] Let $g \in G_W^{max}$. We define the $age$ of $g$ to be
\[
    \age g=\frac{1}{2\pi i}\sum_{j=1}^N \log(\gamma_j),
\]
 where the $\gamma_j$s represent the eigenvalues of $g$ and we choose a branch of the log such that $0 \leq \frac{\log{(\gamma_j)}}{2\pi i} < 1$. 
\end{definition}

The following lemma follows from simultaneously diagonalizing $g$ and $g^{-1}$, and then considering the eigenvalues. 

\begin{lemma}[Mukai \cite{M}] \label{lem:age-g+ginv}
Let $W(x_1,\dots,x_N)$ be a quasihomogeneous polynomial and $G$ be a finite subgroup of the symmetry group $G_W^{max}$. For every $g\in G$, we have
\begin{center}
    $\age(g)+\age(g^{-1})=N-N_g$.
\end{center}
\end{lemma}

There are two important subgroups of $G^{max}_W$ that will play a prominent role in our mirror map. We describe these next.

\subsubsection{Diagonal Symmetries}

\begin{definition} Let $W$ be a nondegenerate quasihomogeneous polynomial. We define the \textit{diagonal symmetry group of W}:
\begin{align*}
    G_W^{diag}=\{(\lambda_1,\dots,\lambda_N)\in(\CC^*)^N:W(\lambda_1x_1,\dots,\lambda_Nx_N)=W(x_1,\dots,x_N)\}.
\end{align*}
If we consider $\lambda\in G_W^{diag}$ as a matrix with (nonzero) entries along the diagonal, then we have a canonical inclusion map $G_W^{diag}\subseteq G_W^{max}$. 

For a nondegenerate quasihomogeneous polynomial $W(x_1,\dots,x_N)$ with weights $(q_1,\dots, q_N)$, we define the \textit{exponential grading operator} as $j_W=(e^{2\pi i q_1}, \dots, e^{2\pi iq_N})$. 
\end{definition}

As mentioned in \cite{M}, the condition on $g_{ij}$ in \eqref{def:Gmax} is equivalent to the requirement that $j_W$ lies in the center of $G^{max}_W$.

By \cite{ABS}, we know that for $\lambda=(\lambda_1,\dots,\lambda_N)\in G_W^{diag}$, the entries $\lambda_i$ are roots of unity. Thus we can represent any element $\lambda\in G^{diag}_W$ (additively) as an element of $(\QQ/\ZZ)^N$:
\begin{center}
$\lambda=(e^{2\pi i\theta_1},\dots,e^{2\pi i\theta_N})\leftrightarrow(\theta_1,\dots,\theta_N)$.  
\end{center}
As we will use both additive notation and multiplicative notation extensively, we will write $\lambda$ to represent a diagonal symmetry as an element of $(\CC^*)^N$, and we will write $[\lambda]$ to mean the corresponding element of $(\QQ/\ZZ)^N$. For example, we would write the exponential grading operator as $[j_W]=(q_1,\dots, q_N)$.

The following fact is useful, which we state without proof. 

\begin{lemma}[Krawitz \cite{Krawitz}] $G_W^{diag}$ is generated (additively) by the columns of $A_W^{-1}$. 
\end{lemma}

\begin{remark}
We will employ the notation used in \cite{Krawitz} as well. Let $\rho_i\in G^{diag}_W$ be the diagonal symmetry of $W$ corresponding to the $i$-th column of $A_W^{-1}$. 
\end{remark}

\begin{example}
Consider the chain potential $W=x_1^4x_2+x_2^5x_3+x_3^3x_4+x_4^2$ from earlier. The weights for this polynomial are $(\frac{5}{24}, \frac{1}{6}, \frac{1}{6}, \frac{1}{2})$. The inverse of the exponent matrix is
\[
\begin{bmatrix}
\tfrac{1}{4} & -\tfrac{1}{20} & \tfrac{1}{60} & -\tfrac{1}{120}\\
0 & \tfrac{1}{5} & -\tfrac{1}{15} & \tfrac{1}{30}\\
0 & 0 & \tfrac{1}{3} & -\tfrac{1}{6}\\
0 & 0 & 0 & \tfrac{1}{2}
\end{bmatrix}.
\]
One can readily see, again using the previous lemma, that the fourth column (i.e., $(e^{-2\pi i/120}, e^{2\pi i/30}, e^{-2\pi i/6}, e^{2\pi i/2})$) generates $G^{diag}_W$, which has order 120. The exponential grading operator $[j_W]$ is equal to the sum of all four columns. 

For the loop potential $W=x_1^3x_2+x_2^3x_3+x_3^3x_4+x_4^3x_5+x_5^3x_6+x_6^3x_1$, the weights are $(\frac{1}{4}, \frac{1}{4}, \frac{1}{4}, \frac{1}{4},\frac{1}{4}, \frac{1}{4})$, as we have seen. The inverse of the exponent matrix is
\[
\frac{1}{728} \begin{bmatrix}
243 & -81 & 27 & -9 & 3 & -1\\
-1 & 243 & -81 & 27 & -9 & 3\\
3 & -1 & 243 & -81 & 27 & -9\\
-9 & 3 & -1 & 243 & -81 & 27\\
27 & -9 & 3 & -1 & 243 & -81\\
-81 & 27 & -9 & 3 & -1 & 243
\end{bmatrix}.
\]
One can also check that any column of $A_W^{-1}$ generates $G_W^{diag}$. This is true for any loop. 
\end{example}

Finally, we want to mention the fixed locus for diagonal symmetries. Let $\lambda\in G^{diag}_W$. Since $\lambda$ is diagonal, it is clear that the fixed locus corresponds with diagonal entries of $1$. In additive notation, if $[\lambda]=(\theta_1,\dots, \theta_N)$ with $0\leq \theta_i<1$, then $\mathbf{e}_i\in \Fix(\lambda)$ if and only if $\theta_i=0$, where $\{\mathbf{e}_i\}$ denotes the standard basis of $\CC^N$.

In this notation, we also note that 
 \[
 \age \lambda=\sum_{i=1}^N \theta_i.
 \]

\subsubsection{Permutations}

The other subgroup that will play a pivotal role in this paper is the subgroup of $G^{max}_W$, which consists of symmetries that permute the variables of $W$. We denote this group by $\Sigma$. It can be constructed as follows.

Let's assume $S_N$ acts on $\{1, 2, \dots, N\}$ on the left. This gives us a left action on $\mathbb C^N$, via 
\[
\sigma \cdot \mathbf{e}_i = \mathbf{e}_{\sigma(i)}.
\]

For $\sigma\in S_N$, in coordinates, the (left) action of $S_N$ on $\CC^N$ is 
\[
\sigma \cdot (x_1, \dots, x_N)= (x_{\sigma^{-1}(1)}, \dots, x_{\sigma^{-1}(n)}). 
\]
We write $[\sigma]$ for the corresponding matrix. We have  
\[
[\sigma]_{ij}=\begin{cases}
1 & \text{if } i=\sigma(j)\\
0 & \text{otherwise}
\end{cases}.
\]
For example, if $n=3$, we would have 
\[
[(123)]=\begin{pmatrix}
0 & 0 & 1\\
1 & 0 & 0\\
0 & 1 & 0
\end{pmatrix}. 
\]

Note that $[\sigma^{-1}] = [\sigma]^T$.

One can check that this is a group homomorphism, if we take the convention on $S_N$ that we multiply from right to left. In other words, $S_N$ acts on $\{1, \dots, N\}$ from the left. For example, we have $(12)(23)=(123)$, and we have 
\[
[(12)][(23)]=\begin{pmatrix}
0 & 1 & 0\\
1 & 0 & 0\\
0 & 0 & 1
\end{pmatrix} \begin{pmatrix}
1 & 0 & 0\\
0 & 0 & 1\\
0 & 1 & 0
\end{pmatrix} =
\begin{pmatrix}
0 & 0 & 1\\
1 & 0 & 0\\
0 & 1 & 0
\end{pmatrix}=[(123)].
\]

The group $\Sigma = \Sigma_W$ contains those elements $\sigma\in S_N$ that preserve $W$.   
Alternatively, we can think of this as $\Sigma=G^{max}_W\cap [S_N]$. 
 Recalling our convention that the $i$-th monomial is of the form $x^a_ix_j$ or $x_i^a$, it follows that $\sigma$ permutes the order of the monomials the same way it permutes variables. Since multiplying a matrix on the right by $[\sigma]$ permutes columns and multiplying on the left by $[\sigma]^T$ permutes rows, we see that $\sigma$ is a symmetry of $W$ if and only if $A_W = [\sigma]^T A_W [\sigma]$, or in other words,
\begin{equation}
 [\sigma]A_W = A_W[\sigma]  \label{eq:sA=As}.
\end{equation}
The permutations that can appear in $\Sigma$ are restricted considerably by the polynomial $W$. We will describe the possible permutations shortly.

In terms of notation, we will write elements of $\Sigma$ in cycle notation, thinking of them as permuting the variables of $W$. In particular, the left action of $\Sigma$ on $\CC^N$ 
induces a right action on $\Hom(\CC^N,\CC)$. In coordinates, if we let $x_i$ be a dual basis, then the action is 
\[
(x_i.\sigma) (v) = x_i(\sigma \cdot v)
\]
for $v \in \mathbb C^N$. (We take the convention that left actions are written with a center dot, and right actions are written with a period.) 
Applying this to the basis, we have
\[
(x_i.\sigma) (\mathbf{e}_j) = x_i(\mathbf{e}_{\sigma(j)}) = \delta_{i\sigma(j)} = \delta_{\sigma^{-1}(i)j},
\]
so it follows that 
\begin{equation}
    x_i.\sigma = x_{\sigma^{-1}(i)}.
\end{equation}
This induces a right action of $\Sigma$ on $\CC[x_1,\dots, x_N]$ compatible with the representation mentioned above.

\begin{example}
If $W=x_1^4+x_2^4+x_3^4+x_4^6$, then we can consider the permutation $(123)\in \Sigma$, which permutes the first three variables. Then we have 
\[
W.\sigma=x_3^4+x_1^4+x_2^4+x_4^6. 
\]
So we can see that $W$ is preserved by $\sigma$.
\end{example}

\subsubsection{Permutations of variables in invertible polynomials} \label{sec:perms}
Most of this section will not be used directly, but we include it as a more precise description of the possible permutations. This is partially described by Ebeling and Gusein-Zade in \cite{EG2}. 
First, we will decompose $W$ into the sums of its atomic types, grouping together those that are ``the same''---meaning of the same type and with the same exponents. To this end, let 
\begin{equation}\label{eqn:WTypesum}
W=\sum_{i=1}^{n_F} \sum_{j=1}^{k_{F_i}} W_{F_{i,j}} + \sum_{i=1}^{n_C}\sum_{j=1}^{k_{C_i}} W_{C_{i,j}} + \sum_{i=1}^{n_L}\sum_{j=1}^{k_{L_i}} W_{L_{i,j}}.
\end{equation}

\noindent Here, the F refers to Fermat type, the C refers to chain type, and the L refers to loop type. Furthermore, we assume $W_{F_{i,j}}=W_{F_{i,\ell}}$ for all $j, \ell$ up to a relabelling of variables.

In other words, we can write $W$ in terms of its atomic types, grouping first the Fermat types together, then the chains, and then the loops. In each grouping, we further group those polynomials that are ``the same''.  

In order to be a symmetry of $W$, the permutation $\sigma$ must be a product of the following types of permutations:

\noindent\textbf{Case 1:} $\sigma_F$ permutes Fermat types. 

These are the easiest to describe. Since Fermat types consist of only one variable, several Fermat types can be permuted as long as they have the same exponent. For example, if we consider $W=x_1^4+x_2^4+x_3^4+x_4^6$, then we can consider the permutation $(123)\in \Sigma$, which permutes the first three variables.

\noindent\textbf{Case 2:} $\sigma_C$ permutes chain type polynomials. 

Similar to the Fermat case, these must permute entire chains. In other words, suppose $W_{C_j}=x_{j,1}^{a_1}x_{j,2} + x_{j,2}^{a_2}x_{j,3} +\dots + x_{j,m}^{a_m}$, for $1\leq j\leq k$ (here, we are using similar notation to Equation~\eqref{eqn:WTypesum}, but suppressing the first index $i$). If $\sigma_C$ cycles through these polynomials sending $W_{C_j}$ to $W_{C_{j+1}}$ with the indices taken modulo $k$, then we can express $\sigma_C$ as 
\[
(i_{1,1}i_{2,1}\dots i_{k,1})(i_{1,2}i_{2,2}\dots i_{k,2})\dots(i_{1,m}i_{2m}\dots i_{k,m}).
\]
Here, $i_{j,\ell}$ denotes the index of the variable $x_{j,\ell}$.

\noindent\textbf{Case 3:} $\sigma_L$ is a bit more complicated, due to the fact that there are two ways to permute the variables in a sum of loops. 
\begin{enumerate}
    \item We have permutations similar to Fermats and chains, which permute entire loops of the same length with the same exponents. 
    \item We have permutations that cycle through the variables in a single loop, e.g., if $W_L=x_1^{a_1}x_2+x_2^{a_2}x_3+x_3^{a_1}x_4+x_4^{a_2}x_1$, we could have the permutation $(13)(24)$, and if furthermore $a_1=a_2$, we would also have the permutation $(1234)$.  
\end{enumerate}

When we combine these two types of permutations for loops, we obtain elements of a wreath product, which we now describe. Suppose $W_{L_j}=x_{j,1}^{a_1}x_{j,2} + x_{j,2}^{a_2}x_{j,3} +\dots + x_{j,m}^{a_m}x_{j,1}$ for $1\leq j\leq k$ (again notated as in Equation~\ref{eqn:WTypesum}, but suppressing the first subscript) and denote by $Z_m$ the cyclic group of order $m$ acting on the variables $x_{j,1},\dots, x_{j,m}$. Then the possible permutations come from the wreath product $Z_m\wr S_k$. Notice that not all cycles in $Z_m$ will be symmetries of $W_{L_j}$. We have to restrict to elements of $Z_m$ that actually preserve the loop. 

Again, we can denote the index of the variable $x_{j,\ell}$ by $i_{j,\ell}$. For the sake of this description, we may assume $\sigma_L$ acts as a $k$-cycle on the $k$ loops $W_{L_j}$. Any other element of the wreath product can be written as a product of such permutations. And furthermore, if two such permutations act on distinct sets of loops, then they commute. 

We can represent $\sigma_L$ as a permutation of the entries of the following matrix (see \cite{frucht}):
\[
\begin{bmatrix}
i_{1,1} & i_{1,2} & \dots & i_{1,m}\\
i_{2,1} & i_{2,2} & \dots & i_{2,m}\\
\vdots &  & \ddots & \vdots\\
i_{k,1} & i_{k,2} & \dots & i_{k,m}
\end{bmatrix},
\]
where we first cycle the rows, followed by cycling the entries of each row (independently) by some cycle in $Z_m$. 

If we write $m_1,m_2,\dots m_k\in \{0,1,\dots, m-1\}$, where $m_j$ denotes the power of the generator of $Z_m$, thought of as an element of $Z_m$ of the cycle on the $j$-th row (i.e., $i_{j,\ell}\to i_{j,m_j+\ell}$), then $\sigma_L$ can be written as a product of $d=\gcd(\sum_{j=1}^k m_j, m)$ cycles of length $\frac{km}{d}$. 

Let us give a few examples: Suppose we have three loops of length 4, with all of the exponents equal (as described in Case 3 above). Then if we put $m_1=2, m_2=1, m_3=0$ (i.e., cycle the first row by 2, the second row by 1, and the third row not at all), we have 
\[
\begin{bmatrix}
i_{1,1} & i_{1,2} & i_{1,3} & i_{1,4}\\
i_{2,1} & i_{2,2} & i_{2,3} & i_{2,4}\\
i_{3,1} & i_{3,2} & i_{3,3} & i_{3,4}
\end{bmatrix}\to
\begin{bmatrix}
i_{2,2} & i_{2,3} & i_{2,4} & i_{2,1}\\
i_{3,1} & i_{3,2} & i_{3,3} & i_{3,4}\\
i_{1,3} & i_{1,4} & i_{1,1} & i_{1,2}
\end{bmatrix}.
\]
This would correspond to the permutation
\[
(i_{1,1}i_{2,2}i_{3,2}i_{1,4}i_{2,1}i_{3,1}i_{1,3}i_{2,4}i_{3,4}i_{1,2}i_{2,3}i_{3,3}).
\]

If we put $m_1=0, m_2=2, m_3=0$, we have 
\[
\begin{bmatrix}
i_{1,1} & i_{1,2} & i_{1,3} & i_{1,4}\\
i_{2,1} & i_{2,2} & i_{2,3} & i_{2,4}\\
i_{3,1} & i_{3,2} & i_{3,3} & i_{3,4}
\end{bmatrix}\to
\begin{bmatrix}
i_{2,3} & i_{2,4} & i_{2,1} & i_{2,2}\\
i_{3,1} & i_{3,2} & i_{3,3} & i_{3,4}\\
i_{1,1} & i_{1,2} & i_{1,3} & i_{1,4}
\end{bmatrix}.
\]
This would correspond to the permutation
\[
(i_{1,1}i_{2,3}i_{3,3}i_{1,3}i_{2,1}i_{3,1})(i_{1,2}i_{2,4}i_{3,4}i_{1,4}i_{2,2}i_{3,2}).
\]

Finally, if we put $m_1=0, m_2=1, m_3=1$, we have 
\[
\begin{bmatrix}
i_{1,1} & i_{1,2} & i_{1,3} & i_{1,4}\\
i_{2,1} & i_{2,2} & i_{2,3} & i_{2,4}\\
i_{3,1} & i_{3,2} & i_{3,3} & i_{3,4}
\end{bmatrix}\to
\begin{bmatrix}
i_{2,2} & i_{2,3} & i_{2,4} & i_{2,1}\\
i_{3,2} & i_{3,3} & i_{3,4} & i_{3,1}\\
i_{1,1} & i_{1,2} & i_{1,3} & i_{1,4}
\end{bmatrix}.
\]
This would correspond to the permutation
\[
(i_{1,1}i_{2,2}i_{3,3}i_{1,3}i_{2,4}i_{3,1})(i_{1,2}i_{2,3}i_{3,4}i_{1,4}i_{2,1}i_{3,2}).
\]

\begin{example}
Let $W=x_1^3x_2+x_2^3x_3+x_3^3x_4+x_4^3x_5+x_5^3x_6+x_6^3x_1$ and $\sigma= (14)(25)(36)$. We see that $W\sigma=x_4^3x_5+x_5^3x_6+x_6^3x_1+x_1^3x_2+x_2^3x_3+x_3^3x_4=W$. In the previous notation, this is $k=1$, $m=6$, and $m_1=3$, the reason why $\sigma$ is a product of three 2-cycles. 
\end{example}

\subsubsection*{The fixed locus for permutations}

We have already described $\Fix(\lambda)$ for diagonal symmetries. We will now describe $\Fix(\sigma)$ for permutations. We will generalize this in Section~\ref{sec:fixed-basis}, but it is important enough to make mention of it on its own. First, we will define following notation.

\begin{definition}\label{def:orbits}
Given $\sigma\in \Sigma$ and $k\in\{1,\dots,N\}$, we denote by $\orbit{\sigma}(k)$ the orbit of $k$ under the action of $\sigma$. If $\tau=(i_1,\dots,i_m)\in \Sigma$ is a cycle, then we denote the primary orbit of $\tau$ as $\orbit{\tau}=\{i_1,\dots,i_m\}$.
\end{definition}

\begin{lemma}\label{lem:eigenvector}
Let $\tau=(i_1,\dots,i_m)\in \Sigma$ be a cycle. 
Define the vector $\mathbf{f}_{\tau}=\sum_{j=1}^m e_{i_j}$. 
Then $\mathbf{f}_{\tau}\in \Fix(\sigma)$. 
Furthermore, if $\{\mathbf{e}_i\}_{i=1}^n$ is the standard basis of $\CC^N$, then the set 
\[
\{\mathbf{e}_j\}_{j\notin\orbit{\tau}}\cup\{\mathbf{f}_{\tau}\}
\]
is a basis for $\Fix(\tau)$. 
\end{lemma}

\begin{proof}
Recall that $\sigma \cdot \mathbf{e}_j=\mathbf{e}_{\sigma(j)}$. For $0\leq \ell\leq m-1$, let $\xi_m=e^{2\pi i/m}$, and define 
\[
\mathbf{f}^{(\ell)}_\tau=\sum_{j=1}^m\xi_m^{-j\ell}\mathbf{e}_{i_j},
\]

Taking the convention that $i_{m+1}=i_1$, we have 
\begin{align*}
\tau \cdot \mathbf{f}^{(\ell)}_\tau &=\sum_{j=1}^m\xi_m^{-j\ell}\mathbf{e}_{\tau(i_j)}\\
  &=\sum_{j=1}^m\xi_m^{-j\ell}\mathbf{e}_{i_{j+1}}\\
  &=\sum_{J=2}^{m+1}\xi_m^{(1-J)\ell}\mathbf{e}_{i_{J}}\\
  &=\xi_m^{\ell} \sum_{J=1}^m\xi_m^{-J\ell}\mathbf{e}_{i_{J}}.
\end{align*}
So the vector $\mathbf{f}_\tau^{(\ell)}$ is an eigenvector for $\tau$ with eigenvalue $\xi_m^{\ell}$. Notice there are exactly $m$ such eigenvectors, and $\mathbf{f}_{\tau}=\mathbf{f}_\tau^{(0)}$. Thus together with the set $\{\mathbf e_j\}_{j\notin\cO_\tau}$, we have a basis for $\CC^N$. 
\end{proof}

As shown in the proof, $\Fix(\tau)$ can be decomposed into two subspaces, the \emph{eigen fixed locus} $\Fix^E(\tau)=\Span\{\mathbf{f}_\tau\}$ and the \emph{tautological fixed locus} $\Fix^T(\tau)=\Span\{e_i\}_{i\notin \orbit{\tau}}$. 

\begin{corollary}
Write $\sigma=\tau_1\dots\tau_{N_\sigma}$ in its disjoint cycle notation, including cycles of length 1. Then the fixed locus $\Fix(\sigma)$ is spanned by $\{\mathbf{f}_{\tau_i}\}_{i=1}^k$, i.e., 
\[
\Fix(\sigma)=\bigoplus_{i=1}^{N_\sigma}\Fix^E(\tau_i).
\]
\end{corollary}

\begin{example}
Suppose $n=4$ and $\sigma=(123)(4)$. Then $\Fix(\sigma)=\{(1,1,1,0), (0,0,0,1)\}$. If, on the other hand $\sigma=(13)(24)$, then $\Fix(\sigma)=\{(1,0,1,0), (0,1,0,1)\}$. 
\end{example}

\begin{convention} When writing a permutation as a product of disjoint cycles, \textbf{we always include cycles of length 1}. 
\end{convention} 

With the eigenvalues of $\sigma$, we can now describe the age of any permutation. First, if $\tau$ is a cycle of length $m$, then we have 
\begin{equation}
\age \tau = \sum_{\ell=0}^{m-1}\frac{\ell}{m}=\frac{m-1}{2}. \label{eq:age-tau}
\end{equation}
Thus, for a general permutation, we as usual write it as a product of disjoint cycles $\sigma=\tau_1\dots\tau_{N_\sigma}$ with $\tau_j$ having length $m_j$. Then we have 
\begin{equation}
    \age \sigma = \sum_{j=1}^{N_\sigma} \frac{m_j-1}{2}=\frac{1}{2}\Big( \sum_{j=1}^{N_\sigma} m_j -N_\sigma \Big).
\end{equation}

\subsubsection{Semidirect Products} \label{sec:semi}

In this section, we will put the diagonal symmetries and the permutations together to form a nonabelian group of symmetries. 

There is a right action of $\Sigma$ on $G^{diag}$ given by \[
\lambda.\sigma = 
(\lambda_{\sigma(1)},\dots, \lambda_{\sigma(n)}). 
\]
This is just (right) conjugation by $\sigma$, as can be seen by the following computation:
\[
\sigma^{-1} \lambda \sigma \cdot \mathbf e_i = \sigma^{-1} \lambda \cdot \mathbf e_{\sigma(i)} = \sigma^{-1} \cdot \lambda_{\sigma(i)}  \mathbf e_{\sigma(i)} = \lambda_{\sigma(i)} \mathbf e_i.
\]
In additive notation, we have
\[
	[\lambda.\sigma] = [\sigma]^T [\lambda].
\]
We take the convention that the ``dot'' of this action has higher operational precedence than the group multiplication; however, we will try to use redundant parentheses when it makes reading easier.

Combining the action of $\lambda$ and $\sigma$, we get the action of $G$ on $\CC^N$ 
\begin{equation}\label{eq:RightActionGeneral}
\sigma \lambda \cdot(x_1,\dots, x_N)=(\lambda_{\sigma^{-1}(1)} x_{\sigma^{-1}(1)}, \dots, \lambda_{\sigma^{-1}(N)} x_{\sigma^{-1}(N)}). 
\end{equation}

We obtain a (right) semidirect product $\Sigma \ltimes G^{diag}$, with the product defined by
\[
(\sigma_1 \lambda_1) (\sigma_2 \lambda_2) = (\sigma_1\sigma_2) ((\lambda_1.\sigma_2) \lambda_2).
\]
We have the commutation rules
\[
 \sigma \lambda = (\lambda.\sigma^{-1}) \sigma
\]
and more importantly
\[
\lambda \sigma = \sigma (\lambda.\sigma).
\]

Let $S\subseteq \Sigma$ and $H\subseteq G^{diag}_W$ such that $S$ is in the normalizer of $H$. For the remainder of this work, we will assume that $G$ is a group of the form $G=S\ltimes H\subseteq G^{max}_W$. Thus, every $g\in G$ can be represented uniquely as $g=\sigma\lambda$, where $\lambda\in G^{diag}_W$ and $\sigma\in S$. This is just the interior product: $S\ltimes H \cong S\cdot H \le G^{max}$. 

There are a few subgroups of $G$ that will be important in the mirror map. First, given $g\in G$, we will denote the centralizer of $g$ by $C_G(g)$. In particular, if we consider only a permutation $\sigma\in S$, then we have also the centralizer $C_S(\sigma)\subset S$. We also define the subgroup $C_H(\sigma)=\{\lambda\in H : \lambda\sigma=\sigma\lambda\}$. The following lemma is obvious, but we will use it later. 
\begin{lemma} \label{lem:CHs}
$\lambda \in C_H(\sigma)$ if and only if $\lambda_i = \lambda_j$ whenever $\cO_\sigma(i) = \cO_\sigma(j)$.
\end{lemma}

Later, we will be interested in the orbits of $H$ acting on the right of $G$ by conjugation.  We define the \emph{commutator} $[H,\sigma]=\{(\mu^{-1}.\sigma)\mu: \mu \in H\}$. This is just the set of commutators $\sigma^{-1}\lambda^{-1}\sigma\lambda$ for a fixed $\sigma$.     

\begin{lemma}
 If $\sigma$ is in the normalizer of $H$, then $[H,\sigma]$ is a subgroup of $H$.
\end{lemma}
\begin{proof}
Since $\sigma$ is in the normalizer of $H$, we have $\lambda.\sigma=\sigma^{-1}\lambda\sigma \in H$ for any $\lambda\in H$. Let $(\lambda^{-1}.\sigma) \lambda $ and $(\mu^{-1}.\sigma)\mu$ be arbitrary elements of $[H,\sigma]$. Then
\begin{align*}
    ((\lambda^{-1}.\sigma) \lambda)^{-1} &= (\sigma^{-1}\lambda^{-1}\sigma\lambda)^{-1}\\
  &=    \lambda^{-1} \sigma^{-1}\lambda\sigma\\
  &= (\sigma^{-1}\lambda\sigma)\lambda^{-1}\\
    &=(\lambda.\sigma)\lambda^{-1}.
\end{align*}
The third equality is true because $H$ is commutative. Furthermore, 
\begin{align*}
(\lambda^{-1}.\sigma) \lambda \cdot (\mu^{-1}.\sigma) \mu &= 
  \sigma^{-1}\mu^{-1}\sigma\mu\sigma^{-1}\lambda^{-1}\sigma \lambda \\
   &=\sigma^{-1}\mu^{-1}\sigma\sigma^{-1}\lambda^{-1}\sigma \lambda \mu\\
   &=((\lambda \mu)^{-1} .\sigma) (\lambda \mu).    
\end{align*}
    
\end{proof}

\begin{corollary} \label{cor:coset-orbits}
    If we let $H$ act on $G$ by right conjugation, then $H$ orbit of $\sigma \lambda$ is the left coset of $[H, \sigma]$ represented by $\sigma \lambda$.
\end{corollary}
\begin{proof}
    One computes $(\sigma \lambda) .\mu =  \mu^{-1} \sigma \lambda \mu =  \sigma \lambda (\mu^{-1}. \sigma) \mu \in \sigma\lambda[H,\sigma]$.
\end{proof}

\subsubsection*{The fixed locus of a general group element}

As with diagonal symmetries and permutations, we describe $\Fix(\sigma\lambda)$ for  any element of $\Sigma\ltimes G^{diag}_W$. 

In order to do so, let us first introduce the following terminology. 

\begin{definition}
Suppose $\tau$ is a cycle of length $m$ and that $\lambda=(\lambda_1,\dots,\lambda_N)\in G^{diag}_W$ satisfies the property that $\lambda_i=1$ if $i\notin \orbit{\tau}$ (recall Definition~\ref{def:orbits}). Then we call the element $\tau\lambda$ a \emph{weighted cycle}. Furthermore, if $\det\lambda=1$, then $\tau\lambda$ is a \emph{special weighted cycle}. 
\end{definition}

We will first describe the fixed locus of weighted cycles, and then put them together to describe $\Fix(\sigma\lambda)$ for general elements of $G$. 

\begin{lemma}\label{lem:eigenvector1}
Let $\tau\lambda$ denote a weighted cycle of length $m$, and let $k$ be the smallest element of the orbit $\orbit{\tau}$. Let $\alpha=(\det\lambda)^{1/m}$ be some fixed $m$-th root of $\det\lambda$ and $\xi_m$ be a primitive $m$-th root of unity. For $0\leq \ell < m$, define the vector 
\[
\mathbf{f}_{\tau\lambda}^{(\ell)} = \sum_{j = 0}^{m-1} (\alpha\xi_m^\ell)^{-j}\left( \prod_{i=0}^{j-1} \lambda_{\tau^{i}(k)} \right) \mathbf{e}_{\tau^j(k)}.
\]
Then $\mathbf{f}_{\tau\lambda}^{(\ell)}$ is an eigenvector of $\tau\lambda$ with eigenvalue $\alpha\xi_m^{\ell}$.
Hence, $\Fix^E(\tau\lambda)$ is one-dimensional if $\det\lambda=1$, and zero-dimensional otherwise.
\end{lemma}

\begin{proof}
Consider the following computation:
\begin{align*}
\tau\lambda\cdot \mathbf{f}_{\tau\lambda}^{(\ell)} &= \sum_{j = 0}^{m-1} (\alpha\xi_m^{\ell})^{-j} \lambda_{\tau^j(k)} \left( \prod_{i=0}^{j-1} \lambda_{\tau^{i}(k)} \right) \mathbf{e}_{\tau^{j+1}(k)}\\
&=\sum_{j = 0}^{m-1} (\alpha\xi_m^\ell)^{-j}  \left( \prod_{i=0}^{j} \lambda_{\tau^{i}(k)} \right) \mathbf{e}_{\tau^{j+1}(k)}\\
&= \sum_{J = 1}^{m} (\alpha\xi_m^\ell)^{ (1-J)}  \left( \prod_{i=0}^{J-1} \lambda_{\tau^{i}(k)} \right) \mathbf{e}_{\tau^{J}(k)}\\
&= \alpha \xi_m^{\ell} \sum_{J = 1}^{m} (\alpha\xi_m^\ell)^{-J}  \left( \prod_{i=0}^{J-1} \lambda_{\tau^{i}(k)} \right) \mathbf{e}_{\tau^{J}(k)}.\\
\end{align*}
Finally, one compares the $J=m$ term of the sum in the last line above with the (missing) $J=0$ term. Of course, $\mathbf{e}_{\tau^m(k)} = \mathbf{e}_k$. The coefficient is $(\alpha\xi_m^\ell)^{-m} \prod_{i=0}^{m-1} \lambda_{\tau^i(k)}$ in the former, and $1$ in the latter. These are the same.

Together with $\Fix^T(\tau\lambda)$, we obtain a basis for $\CC^N$. 

Finally, $\alpha\xi_m^\ell = 1$ for some $\ell$ is equivalent to $\alpha^m=1$, which is equivalent to $\det \lambda = 1$. 
\end{proof}

Given $\sigma\lambda\in \Sigma\rtimes G^{diag}_W$, we can decompose it into a product of weighted cycles
\begin{equation}\label{eq:cycledecomp}
\sigma\lambda = \prod_{i=1}^{N_\sigma} \tau_i\lambda^{(i)},    
\end{equation}
where $\sigma=\tau_1\dots\tau_{N_\sigma}$ is the cycle decomposition of $\sigma$ (as usual including cycles of length 1).

As before, we obtain the following description of the fixed locus of $\sigma\lambda$. First, we decompose into weighted cycle decomposition (including cycles of length 1) as in \eqref{eq:cycledecomp}. Let $F_{\sigma\lambda}$ be the set of $i$ so that $\tau_i\lambda^{(i)}$ is special  (that is, $\det\lambda^{(i)}=1$). 

\begin{corollary}
If $i\notin F_{\sigma\lambda}$, then $\Fix^E(\tau_i\lambda^{(i)})$ is zero-dimensional. For $i \in F_{\sigma\lambda}$, $\Fix^E(\tau_i\lambda^{(i)})$ is one-dimensional. Furthermore, 
\[
\Fix(\sigma\lambda)=\bigoplus_{i \in F_{\sigma\lambda}} \Fix^E(\tau_i\lambda^{(i)}). 
\]
In particular, $\dim\Fix(\sigma\lambda)$ is the number of special weighted cycles in the decomposition \eqref{eq:cycledecomp}.    
\end{corollary} 
We will pick a particular basis for $\Fix(\sigma\lambda)$ in Section~\ref{sec:fixed-basis}.

\begin{example}
Let us again consider the example $W=x_1^4+x_2^4+x_3^4+x_4^6$. We had $[j_W]=(\tfrac{1}{4}, \tfrac{1}{4}, \tfrac{1}{4}, \tfrac{1}{6})$ and the permutation $(123)(4)$, which permutes the variables of $W$. From the corollary, we see that $\Fix(\sigma j_w)$ is trivial. 

On the other hand, if we consider the diagonal symmetry $\lambda=(\tfrac{1}{4}, 0, \tfrac{3}{4},0)\in G^{diag}_W$, then we have $\tau_1=(123)$ and $\tau_2=(4)$. Then $\Fix^E(\tau_1\lambda^{(1)})$ is generated by the vector $\mathbf{f}_{\tau_1\lambda^{(1)}}=(i, i, 1, 0)$, whereas $\Fix^E(\tau_2\lambda^{(2)})$ is generated by the vector $\mathbf{f}_{\tau_2\lambda^{(2)}}=(0,0,0,1)$. So $\Fix(\sigma\lambda)$ is two-dimensional. 
\end{example}

With the eigenvalues of $\sigma\lambda$ in hand, we can compute the age. As with $\age(\sigma)$, we write $\age(\sigma\lambda)$ as a sum of ages of its weighted cycles. 

Therefore, assume first that $\tau\lambda$ is a weighted cycle. Then the eigenvalues of $\tau\lambda$ are $\alpha\xi_m^{\ell}$, with $\alpha=(\det\lambda)^{1/m}$, i.e., the $m$ roots of $\det\lambda$. 
Suppose $\lambda=(\Theta_1,\dots,\Theta_N)$ with $\Theta_i\in\QQ\cap [0,1)$. In particular, $\Theta_i=0$ unless $i\in\orbit{\tau}$, since $\tau\lambda$ is a weighted cycle. Let $a_\ell\in\QQ\cap[0,1)$ and $t_\ell \in \{0,1\}$ such that 
\[
\frac{\ell}{m}+\frac{1}{m}\sum_{i=1}^N \Theta_i =a_\ell +t_\ell.
\]
Then 
\begin{align}
\begin{split}
    \age(\tau\lambda) &=\sum_{\ell=0}^{m-1} a_\ell \\
      &= \frac{1}{m}\sum_{\ell=0}^{m-1}\left(\ell-t_\ell+\sum_{i=1}^N\Theta_i\right)\\
      &= \frac{m-1}{2} + \sum_{i=1}^N\Theta_i -\sum_{\ell=0}^{m-1} t_\ell \\
      &= \age(\tau) + \age(\lambda)  -\sum_{\ell=0}^{m-1} t_\ell. 
\end{split} \label{eq:age-ts}
\end{align}
In particular, 
\[
\age\tau\lambda \equiv \age \tau + \age \lambda \pmod{\ZZ}.
\]
In fact, we get the same fact for any general $\sigma\lambda$, namely
\[
\age\sigma\lambda \equiv \age \sigma + \age \lambda \pmod{\ZZ}.
\]

\subsubsection{A basis for $\Fix(\sigma\lambda)$} \label{sec:fixed-basis}
 We will now describe a nice basis for $\Fix(\sigma\lambda)$.

We are now in the unfortunate position of having two different uses for the natural numbers. On the one hand, we are using them to index the variables of $W$, and these indices show up in the permutations. On the other hand, we are using them to index the cycles in the cycle decomposition of $\sigma$. So we need a way to go from one set of indices to the other. For the following, we fix a cycle decomposition $\sigma=\tau_1\tau_2\dots\tau_{N_\sigma}$ for each permutation. 

For each $i \in \{1, \dots, N\}$ the set of indices on the variables $x_i$, let $r(i) \in {1, \dots, N_\sigma}$ be the unique number so that $i$ is in the cycle $\tau_{r(i)}$. For example, if $\sigma=(123)(4)$, then we can index $\tau_1=(123)$ and $\tau_2=(4)$. Then we have $r(1)=r(2)=r(3)=1$ and $r(4)=2$. 

For $k \in \{1,\dots, N_\sigma\}$ the set of indices on cycles, let $s(k)$ be the smallest number in the cycle $\tau_k$. As in our previous example, we have $s(1)=1$ and $s(2)=4$. Of course, $s$ and $r$ depend on $\sigma$, so when necessary, we will write $r_\sigma(i)$ and $s_\sigma(k)$.

Recall $F_{\sigma\lambda}=\{k\in \{1,\dots,N_\sigma\}: \tau_i\lambda^{(i)} \text{ is special}\}$ is the set of  indexes of the \emph{special cycles} of $\sigma\lambda$.

For $k \in F_{\sigma\lambda}$, we set
 \begin{equation}
     \mathbf f_{k,\sigma \lambda} = \sum_{j=0}^{m-1} \left(\prod_{i=0}^{j-1}  \lambda_{\sigma^i(s(k))}\right) \mathbf e_{\sigma^j(s(k))}.    \label{def:f}
\end{equation}

In fact, we have seen this vector before. In our previous notation, notice that $\mathbf f_{k,\sigma\lambda} = \mathbf{f}_{\tau_k\lambda^{(k)}}^{(\ell)}$ when $\ell$, $\xi_m$, and $\alpha$ are chosen so that $\alpha\xi_m^\ell = 1$. 
By the discussion in the previous section, $\{\mathbf{f}_{i,\sigma\lambda}\}_{i \in F_{\sigma\lambda}}$ is a basis for $\Fix(\sigma\lambda)$.

\begin{definition} \label{def:z}
We let $z_{k,\sigma\lambda}$ be the basis for the dual of $\Fix(\sigma\lambda)$, dual to the basis $\mathbf{f}_{k,\sigma\lambda}$.
\end{definition}

We will generalize this in Lemma~\ref{lem:act-f}, but for now let us mention the following result.  
\begin{lemma} \label{lem:mu-f}
	Let $\mu \in C_H(\sigma)$ and $i \in F_{\sigma\lambda}$. We have
	\[
	\mu\cdot\mathbf{f}_{i,\sigma\lambda} = \mu_{s(i)}\mathbf{f}_{i,\sigma\lambda}
	\]
	and hence
	\[
	 z_{i, \sigma \lambda}.\mu = \mu_{s(i)} z_{i, \sigma \lambda}.
	\]
	In other words, the action of $\mu$ on our preferred basis is the same as scaling by $\mu_{s(k)}$. 
\end{lemma}

\begin{proof}
	Recall that $\mathbf{f}_{i,\sigma\lambda}$ is a linear combination of $\mathbf e_{s(i)}, \mathbf{e}_{\tau(s(i))}, \dots, \mathbf{e}_{\tau^{m-2}(s(i))}$. The property that $\mu \in C_H(\sigma)$ implies that $\mu_{s(i)} = \mu_{ \tau^k(s(i))}$ for all $k$ (see Lemma \ref{lem:CHs}).
\end{proof}

\section{State Spaces} \label{sec:state-spaces}

In this section, we will define the state spaces for both the A- and B-models. We begin by defining the various pieces that we will need. 

\subsection{Milnor Rings} 

\begin{definition}
Let $W(x_1,\dots,x_N)$ be a nondegenerate quasihomogeneous polynomial with weights $q_1,\dots,q_N$. We define the \textit{Milnor ring} of $W$ as the quotient 
\begin{align*}
    \cQ_W=\frac{\CC[x_1,\dots,x_N]}{(\frac{\partial W}{\partial x_1},\dots,\frac{\partial W}{\partial x_N})}.
\end{align*}
\end{definition}
 
It is a well-known fact that $\cQ_W$ is a finite-dimensional vector space (over $\CC$) as long as $W$ is nondegenerate. In that case, the dimension of $\cQ_W$, which we denote by $\mu$, can be computed using the following formula (see, e.g., \cite[Section 12.3, Cor. 3]{arnold}):
\[
\mu=\prod_{i=1}^N\left(\frac{1}{q_i}-1\right).
\]

In fact, the definition of $\cQ_W$ is independent of coordinates (see, for example \cite[Section 4.2, Remark 3]{arnold}).  We will use this fact in the next section, where we define the state space.

For a chain or loop, we can make a tuple of the exponents $\beta=(a_1,\dots, a_N)$, which we call the \emph{vector of exponents}. The following property will help us define a basis for the Milnor ring for a chain potential.

\begin{definition}(see, e.g., \cite{FJJS})
Let $W(x_1,\dots,x_N)$ be a chain potential with vector of exponents $\beta=(a_1,\dots,a_N)$. We say that $\beta$ has the \emph{Chain Property} for $W$ if there exists a $k\in \ZZ$ with $0\leq k\leq \frac{N+1}{2}$ 
such that 
\[
\textbf{$\beta$}=(a_1-1,0,a_3-1,0,\dots,a_{2k-1}-1,0,b_{2k+1},b_{2k+2},\dots,b_N),
\]
where $b_{2k+1}<a_{2k+1}-1$, and $b_i\leq a_i-1$ for all $i>2k$. 
\end{definition}

This allows us to describe a particularly nice basis for $\cQ_W$. 

\begin{theorem}[Kreuzer \cite{Kreuzer}] \label{thm:kreu-basis} Let $W$ be an atomic type polynomial in $\CC[x_1,...,x_N]$.
\begin{enumerate}
    \item If $W=x^{a}$ is Fermat type, a basis for $\cQ_W$ is given by $\{x^r:0\leq r\leq a-2\}$.
    \item If $W=x_1^{a_1}x_2+x_2^{a_2}x_3+\dots+x_N^{a_N}x_1$ is a loop, a basis for $\cQ_W$ is given by
    \begin{align*}
        \left\{\prod_{i=1}^n x_i^{r_i}:0\leq r_i\leq a_i-1\right\}.
    \end{align*}

    \item If $W=x_1^{a_1}x_2+x_2^{a_2}x_3+\dots+x_N^{a_N}$ is a chain, a basis for $\cQ_W$ is given by
    \begin{align*}
        \Bigg\{\prod_{i=1}^N x_i^{r_i}:
        \text{The vector of exponents \textbf{r} has the Chain Property}
        \Bigg\}.
    \end{align*}
\end{enumerate}
\end{theorem}

\begin{example}
We will demonstrate how to find a basis for $\cQ_W$ given the chain potential $W=x_1^4x_2+x_2^5x_3+x_3^3x_4+x_4^2$. The weights for this polynomial are $(\frac{5}{24}, \frac{1}{6}, \frac{1}{6}, \frac{1}{2})$. We can choose a basis of monomials of the form $x_1^{b_1}x_2^{b_2}x_3^{b_3}x_4^{b_4}$, and by Theorem~\ref{thm:kreu-basis}, we see that all the choices for exponent vectors are:
\[
(b_1,b_2, b_3,b_4)\text{ satisfying } b_1\leq 2,b_2\leq4, b_3\leq2,b_4\leq1;
\]
\[
(3,0,b_3,b_4)\text{ satisfying } b_3\leq1, b_4\leq1;
\]
\[
(3,0,2,0).
\]
One can check there are 95 such monomials, and the dimension of the Milnor ring is 
\[
\mu=\frac{19}{5}\cdot 5\cdot 5\cdot 1=95. 
\]
This gives us the basis for the Milnor ring of $W$. 
\end{example}

\subsection{The Definition of State Spaces} 

We now have all of the ingredients needed to define our state spaces. Recall that $W_g:=W|_{\Fix(g)}$.
\begin{lemma}[Mukai \cite{M}] If $N_g$ is greater than zero, then $W_g$ is a nondegenerate quasihomogeneous polynomial.
\end{lemma}

As mentioned in the previous section, we can define $\cQ_{W_g}$ without the use of coordinates. The definition of the state space combines this Milnor ring with the volume form on $\Fix(g)$, which we denote by $\omega_g$. As a vector space, $\cQ_{W_g} \cdot \omega_g$ is the same as $\cQ_{W_g}$, but the volume form is important when considering the $G$ action and the bidegree, which we describe in this section.

\begin{example}
If $\lambda = (\frac13,0,0)$, then the (dual of the) fixed locus of $\lambda$ is spanned by $x_2$ and $x_3$. So we can write the volume form as $\omega_\lambda = dx_2 \wedge dx_3$. One can check, for example, that $\omega_\lambda .(1/3,1/4,0) = i \omega_\lambda$ and $\omega_\lambda.(23) = - \omega_\lambda$.  
\end{example}

\begin{lemma} \label{lem:ss-act}
Let $g\in G$ and let $\conj(g)$ denote the conjugacy class of $g$ in $G$. Then $G$ acts on 
\[
\bigoplus_{h\in\conj(g)}\cQ_{W_h}\cdot\omega_h
\]
on the right, satisfying
\[
\mathcal{Q}_{W_g}\cdot\omega_g . h\subset \mathcal{Q}_{W_{h^{-1}gh}}\cdot\omega_{h^{-1}gh}.
\]

Summing these actions, we get an action of $G$ on 
\[
\bigoplus_{h\in G}\cQ_{W_h}\cdot\omega_h.
\]
\end{lemma}
\begin{proof}
Notice that if $h \in G$ and $\mathbf{x}\in \Fix(g)$, then 
\[
\mathbf x \mapsto h \cdot \mathbf{x} \in \Fix(h g h^{-1})
\]
gives an isomorphism between $\Fix(g)$ and $\Fix(hgh^{-1})$.

Taking duals, an isomorphism of vector spaces gives an isomorphism of volume forms and of Milnor rings:
\[
    h^*: \cQ_{W_g \circ h^{-1}}\cdot \omega_{hgh^{-1}} \to \cQ_{W_{g}} \cdot \omega_g 
\]
One can check that $W_g \circ h^{-1} = W_{hgh^{-1}}$.

Renaming $hgh^{-1}$ as $g$, we get the desired result.
\end{proof}

\begin{definition}

Let $W$ be a quasihomogeneous polynomial and $G \le G_W^{max}$ a group of symmetries. The \emph{unprojected state space} is defined as
\[
    \HH_{W,G}^{un} :=\bigoplus_{g\in G}\cQ_{W_g}\cdot\omega_g
\]

The \emph{state space}, or \emph{Landau--Ginzburg orbifold}, is obtained by taking $G$ invariants (of the action in Lemma \ref{lem:ss-act}):
\begin{equation}\label{eq:unprojected}
    \HH_{W,G} := \left(\HH_{W,G}^{un}\right)^G.
\end{equation}
Sometimes we say the \emph{projected state space} to emphasize the distinction from the unprojected state space.
\end{definition}

We will use the notation $\lfloor p; g\rceil$ to denote a polynomial $p\in\cQ_{W_g}\cdot\omega_g \subseteq \HH^{un}_{W,G}$. This notation puts the group element on equal footing with the polynomial rather than relegating it to a subscript. This is natural since the element of our group and the monomial ``trade places'' in the mirror map (see Section~\ref{sec:ms-diagonal}). 

When writing $p$ explicitly, we typically do not write the volume form, since it is determined by $g$. When using the basis $\{ z_{\sigma\lambda,k} \}$ to express $p$, the $\sigma \lambda$ in the subscript may be omitted as well. 

Next, we want to define the bidegree. First, for $p\in \cQ_{W_g}\cdot \omega_g$, we define $\deg p$ to be the weight of the action of $j_W$ on $p$. It is important to remember that the form $\omega_g$ contributes to $\deg p$. So, for example, if $y_{i_1},\dots,y_{i_m}$ is a basis for $\Fix(g)$, then 
\[
\deg y_{i_1}^{a_1}\dots y_{i_m}^{a_m}\cdot\omega_g = \sum_{j=1}^m (a_j+1)q_{i_j}.\] 
We are now ready to define the bidegree.

\begin{definition}[see also Krawitz \cite{Krawitz}] \label{def:bideg} The \textit{A-model bidegree} of an element $\lfloor p;g\rceil\in \cQ_{W_g}\cdot \omega_g$ is defined to be the ordered pair 
\begin{align*}
    \bideg_A(\lfloor p;g\rceil):= (\deg p+\age g -\age j_W,\; N_g-\deg p+\age g-\age j_W).
\end{align*}

The \textit{B-model bidegree} of an element $\lfloor p;g\rceil$ is defined to be the ordered pair
\begin{align*}
    \bideg_B(\lfloor p;g\rceil):= (\deg p+\age g-\age j_W,\; \deg p+\age g^{-1}-\age j_W).
\end{align*}
\end{definition}

\begin{definition}\label{def:AandBmodels}
    If we use the A-model bidegree, then we write
    \[
        \mathcal A_{W,G} := \HH_{W,G}
    \]
    and call it the \emph{A-model state space}.
    
    If we use the B-model bidegree, then we write
    \[
        \mathcal B_{W,G} := \HH_{W,G}
    \]
    and call it the \emph{B-model state space}.
\end{definition}

\begin{remark}
The original definition of the A-model state space is given in terms of a relative Chen--Ruan cohomology of the stack quotient of the fixed locus (see \cite{FJR:GLSM}). This perspective is necessary when one is interested in the full Frobenius manifold and cohomological field theory of the A-model. However, the A-model state space can be shown to be isomorphic, as a $G$ module, to the definition here (see, e.g., \cite{CIR:11}). In order for this Chen--Ruan cohomology to make sense, the group $G$ must be what is known as A-admissible---which essentially means that it contains the element $j_W$. This condition is not needed for the desired mirror map, so we will not comment on it further. 
The definition of the state space we use here in terms of Milnor rings is more convenient for explicit calculations. 

The B-model state space, sometimes called the orbifolded Milnor ring, was defined as in this section as a candidate for mirror symmetry. There is also a B-admissible condition, but only for groups of diagonal symmetries. 
\end{remark}

The following was proved in \cite{PWW}. It allows us to choose a basis of homogeneous elements (with respect to the bidegrees) for the state space.  
\begin{lemma}[see \cite{PWW}] 
Both the bidegrees of Definition~\ref{def:bideg} are invariant under the action of $G$. That is, for $h \in G$,
\[
	\bideg_A(\lfloor p; g \rceil) = \bideg_A(\lfloor p; g \rceil. h)
\]
if we consider the state space as an A-model and 
\[
	\bideg_B(\lfloor p; g \rceil) = \bideg_B(\lfloor p; g \rceil. h)
\] 
if we use the B-model. 
\end{lemma}

\subsection{State space computation example} \label{sec:state-space-ex}

Let $W =x_1^4+x_2^4+x_3^4+x_4^6$. Let $H$ be the group generated by three elements 
$\alpha = (\frac14,0,-\frac14,0)$, $\beta = (0,\frac14,-\frac14,0)$, and $\gamma = (0,0,\frac12,\frac12)$. Since these are in echelon form, we easily see that this group  has order 32. Let $S = \left<(123)\right>$. One can check that $S$ is in the normalizer of $H$. Let $G = S \ltimes H$.

We give a basis for  $\HH_{W,G}^{un}$, but omit those basis elements $\lfloor p; g \rceil$ that satify $\lfloor p; g \rceil.h = a\lfloor p; g \rceil$ for some $h \in C_G(g)$ and $a \in \CC$ with $a \neq 1$. This is because such elements cannot contribute to the projected state space $\HH_{W,G}$, as we will see. We also compute the B-model bidegree for these basis elements.

\textbf{Case 1: $g=(0,0,0,0)$} Using Theorem~\ref{thm:kreu-basis}, a basis for the Milnor ring is given by 
\[
\{ x_1^ax_2^bx_3^cx_4^d : 0 \le a \le 3, 0 \le b \le 3, 0 \le c \le 3, 0 \le d \le 5\}.
\] However, only those below are fixed by $H$. (Don't forget the volume form $dx_1 dx_2 dx_3 dx_4$.) Each of these is $S$ invariant.

\begin{center}
\begin{tabular}{c|c}
     \textbf{Basis Element} &  \textbf{Bidegree}\\
     \hline
     $\lfloor 1; \text{id}\rceil$ & $(0,0)$\\
     $\lfloor x_1x_2x_3x_4; \text{id}\rceil$ & $(\frac{11}{12},\frac{11}{12})$\\
     $\lfloor x_1^2x_2^2x_3^2x_4^2; \text{id}\rceil$ & $(\frac{11}{6},\frac{11}{6})$\\
     $\lfloor x_1x_2x_3x_4^3; \text{id}\rceil$ & $(\frac{4}{3},\frac{4}{3})$\\
     $\lfloor x_1^2x_2^2x_3^2x_4^4; \text{id}\rceil$ & $(\frac{13}{12},\frac{13}{12})$\\
     $\lfloor x_1^2x_2^2x_3^2; \text{id}\rceil$ & $(\frac{3}{2},\frac{3}{2})$\\
     $\lfloor x_4^2; \text{id}\rceil$ & $(\frac{1}{3},\frac{1}{3})$\\
     $\lfloor x_4^4; \text{id}\rceil$ & $(\frac{2}{3},\frac{2}{3})$\\
\end{tabular}
\end{center}

\textbf{Case 2: $g = \lambda \in H$ has trivial fixed locus}. Then the Milnor ring is one-dimensional, and the $H$ action is trivial. $S$ acts by permuting the entries of $\lambda$. The elements here are grouped into three $S$ orbits, shown here by adjacent rows with the same shading.

\begin{center}
\begin{tabular}{c|c}
     \textbf{Basis Element} &  \textbf{Bidegree}\\
     \hline
     $\lfloor 1; (\frac{1}{2},\frac{1}{2},\frac{1}{2},\frac{1}{2})\rceil$ & $(\frac{13}{12},\frac{13}{12})$\\
     \rowcolor{TableShade}
     $\lfloor 1; (\frac{1}{2},\frac{3}{4},\frac{1}{4},\frac{1}{2})\rceil$ & $(\frac{13}{12},\frac{13}{12})$ \\
     \rowcolor{TableShade}
     $\lfloor 1; (\frac{1}{4},\frac{1}{2},\frac{3}{4},\frac{1}{2})\rceil$ & \\
     \rowcolor{TableShade}
     $\lfloor 1; (\frac{3}{4},\frac{1}{4},\frac{1}{2},\frac{1}{2})\rceil$ & \\

     $\lfloor 1; (\frac{3}{4},\frac{1}{2},\frac{1}{4},\frac{1}{2})\rceil$ & $(\frac{13}{12},\frac{13}{12})$ \\
     $\lfloor 1; (\frac{1}{4},\frac{3}{4},\frac{1}{2},\frac{1}{2})\rceil$ &\\
     $\lfloor 1; (\frac{1}{2},\frac{1}{4},\frac{3}{4},\frac{1}{2})\rceil$ &\\
\end{tabular}
\end{center}

\textbf{Case 3: $g = \lambda \in H$ fixes $x_4$.}

\begin{center}
\begin{tabular}{c|c}
     \textbf{Basis Element} &  \textbf{Bidegree}\\
     \hline
     $\lfloor x_4; (\frac{1}{4},\frac{1}{4},\frac{1}{2},0)\rceil$&$(\frac{5}{12},\frac{7}{12})$\\
     $\lfloor x_4; (\frac{1}{4},\frac{1}{2},\frac{1}{4},0)\rceil$&\\
     $\lfloor x_4; (\frac{1}{2},\frac{1}{4},\frac{1}{4},0)\rceil$&\\ 
     \rowcolor{TableShade}
     $\lfloor x_4; (\frac{3}{4},\frac{3}{4},\frac{1}{2},0)\rceil$ & $(\frac{7}{12},\frac{5}{12})$ \\
     \rowcolor{TableShade}
     $\lfloor x_4; (\frac{3}{4},\frac{1}{2},\frac{3}{4},0)\rceil$ & \\
     \rowcolor{TableShade}
     $\lfloor x_4; (\frac{1}{2},\frac{3}{4},\frac{3}{4},0)\rceil$ &\\
     $\lfloor x_4^3; (\frac{1}{2},\frac{1}{4},\frac{1}{4},0)\rceil$ & $(\frac{3}{4},\frac{7}{4})$\\
     $\lfloor x_4^3; (\frac{1}{4},\frac{1}{2},\frac{1}{4},0)\rceil$ &\\
     $\lfloor x_4^3; (\frac{1}{4},\frac{1}{4},\frac{1}{2},0)\rceil$ & \\
     \rowcolor{TableShade}
     $\lfloor x_4^3; (\frac{3}{4},\frac{3}{4},\frac{1}{2},0)\rceil$  &$(\frac{7}{4},\frac{3}{4})$\\
     \rowcolor{TableShade}
     $\lfloor x_4^3; (\frac{1}{2},\frac{3}{4},\frac{3}{4},0)\rceil$ &\\
     \rowcolor{TableShade}
     $\lfloor x_4^3; (\frac{3}{4},\frac{1}{2},\frac{3}{4},0)\rceil$ &
\end{tabular}
\end{center}

One can check that if $g \in H$, $g \neq \text{id}$, and $g$ fixes any of $x_1, x_2$ or $x_3$, then there are no elements fixed by $H$.

\textbf{Case 4: $g$ is conjugate to (123) or (132).}
\noindent The (dual of the) fixed locus is spanned by $z_{1,g}$ and $z_{4,g}$. One can also check that in this basis, we have $W_{g} = z_{1,g}^4 + z_{4,g}^6$ up to rescaling the monomials. Recall from Corollary~\ref{cor:coset-orbits} that the conjugacy class of $\sigma$ in $G$ is $\sigma[H,\sigma]$. One can check that $[H,\sigma]$ is the subgroup generated by the first two generators $\alpha = (\frac14,0,-\frac14,0)$ and $\beta = (0,\frac14,-\frac14,0)$ and is of order 16. 
  In the table below, each row lists elements with the ``same'' monomial in conjugate sectors. Up to scalars, each of these is an $H$ orbit. Each of the $H$ orbits here are invariant (up to scalars) under $S$ (as a set, not pointwise).

\begin{center}
\begin{tabular}{c|c}
     \textbf{Basis Element} &  \textbf{Bidegree} \\
     \hline
     $\lfloor 1; (123)\rceil, \lfloor 1; (123) \alpha \rceil, \dots (\text{14 more})$ & $(\frac{1}{2},\frac{1}{2})$  \\
     $\lfloor z_1z_4;(123)\rceil,\lfloor z_1z_4; (123) \alpha \rceil, \dots (\text{14 more})$ & $(\frac{11}{12},\frac{11}{12})$ \\
     $\lfloor z_1z_4^3;(123)\rceil,\lfloor z_1z_4^3; (123) \alpha \rceil, \dots (\text{14 more})$ & $(\frac{5}{4},\frac{5}{4})$\\
     $\lfloor z_1^2;(123)\rceil,\lfloor z_1^2; (123) \alpha \rceil, \dots (\text{14 more})$ & $(1,1)$\\
     $\lfloor z_1^2z_4^2;(123)\rceil,\lfloor z_1^2z_4^2;(123)\alpha\rceil, \dots (\text{14 more})$ & $(\frac{4}{3},\frac{4}{3})$\\
     $\lfloor z_1^2z_4^4;(123)\rceil,\lfloor z_1^2z_4^4;(123)\alpha\rceil, \dots (\text{14 more})$ & $(\frac{5}{3},\frac{5}{3})$\\
     $\lfloor z_4^2;(123)\rceil,\lfloor z_4^2;(123)\alpha\rceil, \dots (\text{14 more})$ & $(\frac{5}{6},\frac{5}{6})$\\
     $\lfloor z_4^4;(123)\rceil,\lfloor z_4^4;(123)\alpha\rceil, \dots (\text{14 more})$ & $(\frac{7}{6},\frac{7}{6})$\\
     $\lfloor 1; (132)\rceil,\lfloor 1; (132)\alpha\rceil , \dots (\text{14 more})$ & $(\frac{1}{2};\frac{1}{2})$\\
     $\lfloor z_1z_4;(132)\rceil,\lfloor z_1z_4;(132)\alpha\rceil, \dots (\text{14 more})$ & $(\frac{11}{12},\frac{11}{12})$\\
     $\lfloor z_1z_4^3;(132)\rceil,\lfloor z_1z_4^3;(132)\alpha\rceil, \dots (\text{14 more})$ & $(\frac{5}{4},\frac{5}{4})$\\
     $\lfloor z_1^2;(132)\rceil,\lfloor z_1^2;(132)\alpha\rceil, \dots (\text{14 more})$ & $(1,1)$\\
     $\lfloor z_1^2z_4^2;(132)\rceil,\lfloor z_1^2z_4^2;(132)\alpha\rceil, \dots (\text{14 more})$ & $(\frac{4}{3},\frac{4}{3})$\\
     $\lfloor z_1^2z_4^4;(132)\rceil,\lfloor z_1^2z_4^4;(132)\alpha\rceil, \dots (\text{14 more})$ & $(\frac{5}{3},\frac{5}{3})$\\
     $\lfloor z_4^2;(132)\rceil,\lfloor z_4^2;(132)\alpha\rceil, \dots (\text{14 more})$ & $(\frac{5}{6},\frac{5}{6})$\\
     $\lfloor z_4^4;(132)\rceil,\lfloor z_4^4;(132)\alpha\rceil, \dots (\text{14 more})$ & $(\frac{7}{6},\frac{7}{6})$\\
\end{tabular}
\end{center}

\textbf{Case 5: $g$ is conjugate to $(123) \gamma$ or $(132)\gamma$}. Recall $\gamma = (0,0,\frac12,\frac12)$ is one of the generators listed above. These sectors are all narrow. Each row is an $H$ orbit, and is $S$ invariant.

\begin{center}
\begin{tabular}{c|c}
     \textbf{Basis Element} &  \textbf{Bidegree} \\
     \hline
     $\lfloor 1; (123) \gamma \rceil, \lfloor 1; (123) \gamma \alpha \rceil, \dots (\text{14 others})\rceil$ & $(\frac{13}{12},\frac{13}{12})$\\
     $\lfloor 1; (132) \gamma \rceil, \lfloor 1; (132) \gamma \alpha \rceil, \dots (\text{14 others})\rceil$ & $(\frac{13}{12},\frac{13}{12})$ 
\end{tabular}
\end{center}

To pass to the projected state space, one can sum over $G$ orbits.

\subsection{Relation to other work} In \cite{M}, the state space has a different definition. We will show that the definition agrees with our definition of the B-model. 

\begin{definition}[see \cite{M}] Let $W(x_1,...,x_N)$ be a nondegenerate quasihomogeneous polynomial, $G$ a finite subgroup of the symmetry group $G_W^{max}$, and $\cK\subset G$ a set of representatives of the conjugacy classes of $G$. The \textit{Landau--Ginzburg orbifold }$\HH^{alt} _{W,G}$ for the pair $(W,G)$ is defined by

\[
    \HH_{W,G}^{alt}=\bigoplus_{g\in \cK}(\cQ_{W_g}\cdot\omega_g)^{C_G(g)}.
\]
\end{definition}

\begin{lemma}[see \cite{M}] If $N_g$ is greater than zero, then centralizer $C_G(g)$ of $G$ is a finite subgroup of the symmetry group $G_{W_g}$.
\end{lemma}

\begin{lemma}[see \cite{M}] \label{lem:conj-inv} An unprojected  sector does not depend on a choice of a representative of a conjugacy class. More precisely, if $g$ and $h$ are conjugate in $G$, then we have $\cQ_{W,g}\simeq\cQ_{W,h}$ and $\cQ_{W,g}\cdot\omega_g\simeq\cQ_{W,h}\cdot\omega_h$ as bigraded vector spaces.
\end{lemma}

\begin{proposition}\label{prop:statespaceiso}
There exists an isomorphism of vector spaces between the Landau--Ginzburg orbifold $\HH_{W,G}^{alt}$ defined in \cite{M} and the Landau--Ginzburg orbifold $\mathcal{B}_{W,G}$ defined in Equation~\eqref{eq:unprojected} and Definition~\ref{def:AandBmodels}.
\end{proposition}

\begin{proof}
Fix $g \in S$. By Lemma~\ref{lem:ss-act}, we have a decomposition
\[
\mathcal{B}_{W,G}=\left(\bigoplus_{g\in G}\cQ_{W_g}\cdot\omega_g\right)^G =\bigoplus_{g\in S}\left(\bigoplus_{h \in \conj(g)} \cQ_{W_h}\cdot\omega_h\right)^G.
\]
Therefore, it is sufficient to provide an isomorphism 
\[
\phi_g:\left(\cQ_{W_g}\cdot\omega_g\right)^{C_G(g)} \rightarrow \left( \bigoplus_{h \in \conj(g)} \mathcal{Q}_{W_h}\cdot \omega_h \right)^G.
\]
We will denote elements of the left-hand side simply by $p$ and elements of the right-hand side using the notation $\lfloor p;g\rceil$.

For $p\in \left(\cQ_{W_g}\cdot\omega_g\right)^{C_G(g)}$, define 
\[
\phi_g(p)=\frac{1}{|C_G(g)|}\sum_{h\in G} \lfloor p;g \rceil . h.
\]

It is clear that $\phi_g(p)$ is $G$ invariant, so $\phi_g(p)\in \left( \bigoplus_{h \in [g]} \mathcal{Q}_{W_h}\cdot \omega_h \right)^G$.

We need to show that $\phi_g$ is a bijection.
Assume $p\in\ker(\phi_g)$, i.e., $\phi_g(p)=0$. 
Then we have
\begin{align*}
\phi_g(p)&=\frac{1}{|C_{G}(g)|}\sum_{h\in G}\lfloor p;g \rceil. h\\
&=\Big\lfloor\big(\sum p_{i,g_1});g_1\Big\rceil+
\Big\lfloor\big(\sum p_{i,g_2});g_2\Big\rceil+\dots+\Big\lfloor\big(\sum p_{i,g_n});g_n\Big\rceil,
\end{align*}
where $g_i$ are the elements of $\conj(g)$ and $p_{i,g_i}$ are result of the action of $h$ on $p$, for each $h$ satisfying $h^{-1}gh=g_1$. Each of these summands must be equal to 0 since each they are linearly independent and $\phi_g(p)=0$. Since $g\in \conj(g)$, we may assume $g_1=g$. Furthermore, all $h\in G$ that contribute to the summand $\lfloor \sum p_{i,g_1};g_1 \rceil \in \cQ_{W_{g}}\cdot \omega_{g}$ must satisfy $h \in C_G(g)$. However, $p$ is invariant under the action of $C_G(g)$. Thus, $\lfloor \sum p_{i,g_1};g_1 \rceil= |C_G(g)|\lfloor p;g\rceil$. Hence, $p = 0$.
Therefore, $\ker(\phi_g)$ is trivial, and $\phi$ is injective.

To prove surjectivity, we choose $y\in \left( \bigoplus_{g_i \in \conj(g)} \cQ_{W_h}\cdot \omega_h \right)^G$. Again, $g\in \conj(g)$, so we can write
\begin{align*}
    y=\lfloor q;g\rceil+\sum_{g_i \in \conj(g), g_i \neq g} \lfloor q_i;g_i\rceil.
\end{align*}
Since $y$ is $G$ invariant, for any $h \in C_G(g)$ we have $y. h =y$. But also 
\begin{align*}
y.h &=\lfloor p;g\rceil.h+ \sum_{g_i \in \conj(g), g_i \neq g} \lfloor q_i;g_i\rceil. h\\
&=\lfloor  q. h;g\rceil+\sum_{g_i \in \conj(g), g_i \neq g} \lfloor q_i. h;h^{-1}g_i h\rceil\\
&=y.
\end{align*}
So $\lfloor q;g\rceil.h=\lfloor q;g\rceil$. 
We need to show now that $\phi_g(q)=y$.

To this end, take any element $h\in G\setminus C_G(g)$, and act on $y$. Since $y$ is invariant, we still have $y. h = y$. But now we have $h^{-1}g h\neq g$. So $h^{-1} g h= g_i\in \conj(g)$ for some $i$. But then we have
\begin{align*}
    \lfloor q;g \rceil .h = \lfloor q_i;g_i \rceil
\end{align*}
for some $i$. 

We can glue together the various $\phi_g$ for each representative in $S$.  Therefore, as vector spaces, we can conclude that $\HH_{W,G}\cong\mathcal{B}_{W,G}$.
\end{proof}

In \cite{Bas}, Basalaev and Ionov construct an \textit{a priori} different candidate for the B-model state space. It is constructed out of Milnor rings, as we do here, but the $G$ action is determined by braided commutativity of a certain Hochschild cohomology ring. However, they show that the state space is isomorphic to $\HH^{alt}_{W,G}$ defined above. 

\section{BHK Mirror Symmetry for Diagonal Groups}\label{sec:ms-diagonal}

We are now prepared to define mirror symmetry. First, we will explain how to start with a pair $(W,G)$ and obtain the mirror LG pair $(W^\vee,G^\vee)$. This will require us to first recall BHK mirror symmetry from \cite{Krawitz}, \cite{BH}, and \cite{BergHenn}. 

We will use two different notations for the mirror LG model. In the case that $G$ is a group of diagonal symmetries, we will denote the BHK mirror model by $(W^T,G^T)$ as is common. In the next section, we will denote the mirror model for nonabelian LG models by $(W^\vee, G^\vee)$.

\begin{definition}[see, e.g., \cite{BH}] Let $W(x_1,\dots,x_N)$ be an invertible polynomial. The \textit{BHK dual polynomial} $W^T$ is the invertible polynomial determined by the exponent matrix $(A_W)^T$.
\end{definition}

\begin{definition}
Let $W$ be an invertible polynomial, and $H\subseteq G^{diag}_W$. The diagonal \textit{dual group} of $H$ is the group
\begin{align*}
    H^T=\{\lambda\in G_{W^T}^{diag}: [\lambda]^T A_W [\mu]\in\ZZ \text{ for all } \mu\in H\}.
\end{align*}

\end{definition}

\begin{lemma}\label{lem:dual-props}
	\mbox{}
	\begin{enumerate}
	    \item $\age(j_W) = \age(j_{W^T})$.
		\item $(H^T)^T = H$.
		\item $|H|\cdot |H^T| = |G_W^{diag}|$.
		\item If $H_1 \le H_2$, then $H_2^T \le H_1^T$ and $H_1/H_2 \cong H_2^T/H_1^T$.
	\end{enumerate}
\end{lemma}

Part 1 of the previous lemma follows from \cite[Equation (3)]{ChiRu:09} by Chiodo and Ruan. The remainder are proven by Artebani, Boisierre and Sarti in \cite[Proposition 3.5]{ABS}.

\begin{example}
We return to the previous example, namely $W=x_1^4+x_2^4+x_3^4+x_4^6$ with $G=\langle j_W\rangle$.

The BHK dual $W^T$ is defined by the exponent matrix 
\begin{align*}
 A_{W}^T&=\begin{bmatrix}
4 & 0 & 0 & 0\\
0 & 4 & 0 & 0 \\
0 & 0 & 4 & 0 \\
0 & 0 & 0 & 6 \\
\end{bmatrix}.
\end{align*}
Notice $A_W=(A_W)^T$, so $W^T=W=x_1^4+x_2^4+x_3^4+x_4^6$. The dual group $G^T$ is $\SL_W$ generated by $\langle (\frac{1}{4},\frac{1}{4},\frac{1}{2},0),(\frac{1}{4},\frac{1}{2},\frac{1}{4},0),(\frac{1}{2},0,0,\frac{1}{2})\rangle$.
\end{example}

\begin{example}
If we set $W=x_1^3x_2+x_2^3x_3+x_3^3x_4+x_4^3x_5+x_5^3x_6+x_6^3x_1$, our exponent matrices are
\begin{align*}
 A_W&=\begin{bmatrix}
3 & 1 & 0 & 0 & 0 & 0\\
0 & 3 & 1 & 0 & 0 & 0\\
0 & 0 & 3 & 1 & 0 & 0 \\
0 & 0 & 0 & 3 & 1 & 0 \\
0 & 0 & 0 & 0 & 3 & 1 \\
1 & 0 & 0 & 0 & 0 & 3
\end{bmatrix}
&
(A_W)^T&=\begin{bmatrix}
3 & 0 & 0 & 0 & 0 & 1 \\
1 & 3 & 0 & 0 & 0 & 0\\
0 & 1 & 3 & 0 & 0 & 0 \\
0 & 0 & 1 & 3 & 0 & 0 \\
0 & 0 & 0 & 1 & 3 & 0\\
0 & 0 & 0 & 0 & 1 & 3\\
\end{bmatrix},
\end{align*}
and thus $W^T=x_1^3x_6+x_2^3x_1+x_3^3x_2+x_4^3x_3+x_5^3x_4+x_6^3x_5$. Recall that for loops, the group $G^{diag}_W$ is cyclic, generated by any column of the inverse matrix $A_W^{-1}$. In this case, $|G^{diag}_W| = 728$. If we let $G=\langle (\frac{7}{8}, \frac{3}{8}, \frac{7}{8}, \frac{3}{8}, \frac{7}{8}, \frac{3}{8}) \rangle$, we can notice this is the 91st power of the last column. Furthermore, the square of this generator is equal to $j_W$. 

By Lemma~\ref{lem:dual-props}, we see that the dual group $G^T$ must be cyclic of order 91. This is the group $G^T=\langle (\frac{10}{91}, \frac{27}{91}, \frac{82}{91}, \frac{3}{91}, \frac{90}{91}, \frac{61}{91}) \rangle$ generated by the 8th power of the last row of $A_W^{-1}$. 
\end{example}

\begin{definition}[Mirror Map, see \cite{Krawitz}] Let $W$ be an invertible potential and $H\subseteq G^{diag}_W$. We will use variables $x_i$ for $W$ and $y_i$ for $W^T$. We write $\rho_i$ for the symmetries of $W$ corresponding to the columns of the exponent matrix $A_W$ and $\bar \rho_i$ for the columns of $A_{W^T}$. For $\lambda \in H$, we write $F_\lambda = \{ i : x_i.\lambda = x_i \}$.

Given $\displaystyle \left\lfloor \prod_{j \in F_\mu} y_i^{\alpha_i} dy_i; \mu \right\rceil$, we set $\displaystyle \lambda = \prod_{j \in F_\mu} \rho_j^{\alpha_j + 1}$. By \cite[Lemma 2.2]{Krawitz}, there is (usually, see Remark below) a unique element $\displaystyle \prod_{i \in F_\lambda} x_i^{r_i}dx_i$ so that $\mu= \displaystyle\prod_{i \in F_\lambda} \bar \rho_i^{r_i +1}$. Then we can define the map 
\begin{gather} 
    \HH_{W^T,G^T}^{un} \rightarrow \HH_{W,G}^{un}
    \\
    \left \lfloor \prod_{j \in F_\mu} y_j^{\alpha_j}dy_j;\prod_{i \in F_\lambda} \bar{\rho}_i^{r_i+1} \right \rceil
    \mapsto 
    \left\lfloor \prod_{i \in F_\lambda} x_i^{r_i}dx_i;\prod_{j \in F_\mu} \rho_j^{\alpha_j+1} \right\rceil \label{eq:kmm}.
\end{gather}
A key observation is that $\prod_{j \in F_\mu} y_j^{\alpha_j}dy_j$ is invariant under $H^T$ if and only if $\prod_{j \in F_\mu} \rho_j^{\alpha_j+1} \in H$, and analogously for the $x_i$ and $\bar \rho_i$.

In \cite{Krawitz}, Krawitz showed that this map descends to a bidegree-preserving isomorphism between $\mathcal B_{W^T,H^T}$ and $\mathcal A_{W,H}$.

\begin{remark}
In the case of a loop polynomial with an even number of variables, the identity group element can be written in two different ways. In this case, the mirror map of Krawitz in \cite{Krawitz} matches up two-dimensional subspaces from each side (see also \cite{Kreuzer}). 
\end{remark}

\end{definition}

\begin{example}
For example, we again consider the example above with $W=x_1^4+x_2^4+x_3^4+x_4^6$ and $G=\langle j_W \rangle$. In order to be precise, we will include the volume form in this example, to emphasize the mirror map. In order to make the presentation more compact, we will use the convention of Krawitz in \cite{Krawitz}, writing the wedge product simply as a product, but the reader should be aware that we do mean the anticommutative wedge product. 

One of the invariant elements from the state space is $\lfloor 1; j_W\rceil$. We can write $j_W=\rho_1\rho_2\rho_3\rho_4$. So Krawitz mirror map would give us
\[
\lfloor 1; j_W\rceil \leftrightarrow \lfloor dx_1 dx_2 dx_3 dx_4;\id \rceil.
\]
On the other hand, $j_W^5=\rho_1\rho_2\rho_3\rho_4^5$, so we obtain
\[
\lfloor 1; j_W^5\rceil \leftrightarrow \lfloor x_4^4dx_1 dx_2 dx_3 dx_4;\id \rceil.
\]
On the other hand, one can check that $\lfloor x_1^2x_2x_3dx_1 dx_2 dx_3 dx_4; \id\rceil $ is invariant, and by Krawitz's mirror map we obtain
\[
\lfloor x_1^2x_2x_3dx_1 dx_2 dx_3 dx_4; \id\rceil\leftrightarrow \lfloor 1; \rho_1^3\rho_2^2\rho_3^2\rho_4 \rceil,
\]
and the last element is $\rho_1^3\rho_2^2\rho_3^2\rho_4=(\tfrac{3}{4},\tfrac{1}{2}, \tfrac{1}{2}, \tfrac{1}{4})$, which we can check is an element of $G^T=\SL_W$. 
\end{example}

\begin{example}
For another example, consider the LG model defined by the pair $W=x_1^3x_2+x_2^3x_3+x_3^3x_4+x_4^3x_5+x_5^3x_6+x_6^3x_1$ and $G=\langle (\frac{7}{8}, \frac{3}{8}, \frac{7}{8}, \frac{3}{8}, \frac{7}{8}, \frac{3}{8}) \rangle$. If we let $g=(\frac{7}{8}, \frac{3}{8}, \frac{7}{8}, \frac{3}{8}, \frac{7}{8}, \frac{3}{8})$, then we recall that $g$ is the 91st power of $\rho_4$. In order to apply Krawitz's mirror map, we first use the relations $\rho_i^3=\rho_{i-1}^{-1}$ to write
\[
\rho_4^{91}=\rho_2\rho_4\rho_6= \rho_1^3\rho_2^2\rho_3^3\rho_4^2\rho_5^3\rho_6^2. 
\]
Then Krawitz's mirror map gives us
\[
\lfloor 1; \rho_1^3\rho_2^2\rho_3^3\rho_4^2\rho_5^3\rho_6^2\rceil \leftrightarrow \lfloor x_1^2 x_2 x_3^2 x_4 x_5^2 x_6 dx_1 dx_2 dx_3 dx_4 dx_5 dx_6 ; \id\rceil,
\]
which we can check is invariant under $G^T=\langle (\frac{10}{91}, \frac{27}{91}, \frac{82}{91}, \frac{3}{91}, \frac{90}{91}, \frac{61}{91}) \rangle$.
\end{example}

\subsection{Compatibility of diagonal mirror map with permutations} \label{sec:compat}

Recall from \eqref{eq:sA=As} that $\sigma$ is a symmetry of $W$ if and only if $A_W = [\sigma]^T A_W [\sigma]$. Taking transposes, we see that $\sigma$ is a symmetry of $W$ if and only if it is a symmetry of $W^T$. If we apply this to the inverse matrix, we obtain $A_W^{-1} = [\sigma]^T A_W^{-1} [\sigma]$, from which we also learn that 
\[
	[\sigma]^T\Big[[\rho_1] \cdots [\rho_N]\Big] =  \Big[[\rho_1] \cdots [\rho_N]\Big] [\sigma]^T = \Big[[\rho_{\sigma^{-1}(1)}] \cdots [\rho_{\sigma^{-1}(n)}]\Big],
\]
where $[\rho_i]$ is the $i$-th column of the inverse matrix, thought of as a diagonal symmetry of $W$ (written additively). This computation shows that $\rho_i.\sigma = \rho_{\sigma^{-1}(i)}$, i.e., $\sigma$ acts on the standard generators of $G_W^{diag}$ the same way it acts on variables.

We will use the following lemma in Section~\ref{sec:reducing}, not applied to $W$, $H$, and $S$, but to the restriction of $W$ to the fixed locus of $\sigma$ (see also Lemma~\ref{lem:permutations}).

\begin{lemma}\label{lem:kraw-comp} Let $H$ be a group of diagonal symmetries of $W$ and $S$ be a group of permutation symmetries so that $S$ is in the normalizer of $H$. Then there is a right action of $S$ on the state space $\HH_{W,H}$. This action is compatible with the Krawitz mirror map, up to the sign of $S$ acting on the volume form.
\end{lemma}
\begin{proof}
	Of course, the right action is defined by
	\[
		\lfloor p; \mu \rceil.\sigma = \lfloor p.\sigma; \mu.\sigma \rceil.
	\]
	
Acting by $\sigma$ on the left-hand side of \eqref{eq:kmm}, we get
\[
\Big\lfloor \prod_{j \in F_\mu} y_{\sigma^{-1}(j)}^{\alpha_j} dy_{\sigma^{-1}(j)}  ;\prod_{j \in F_{\lambda}} \bar{\rho}_{\sigma^{-1}(j)}^{r_j+1}\Big\rceil.
\]
Notice $j \in F_{\mu.\sigma}$ iff $\mu_{\sigma(j)} = 0$, so $\sigma^{-1}(j) \in F_{\mu.\sigma}$ iff $\mu_j = 0$ iff $j \in F_\mu$. 

So the above is equal to 
\[
\Big\lfloor \prod_{\sigma^{-1}(j) \in F_{\mu.\sigma} } y_{\sigma^{-1}(j)}^{\alpha_j} dy_{\sigma^{-1}(j)}  , \prod_{\sigma^{-1}(j) \in F_{\lambda.\sigma}} \bar{\rho}_{\sigma^{-1}(j)}^{r_j+1}\Big\rceil.
\]
Of course, we can replace $\sigma^{-1}(j)$ with another letter $k$, so we have
\[
\Big\lfloor \prod_{k \in F_{\mu.\sigma} } y_{k}^{\alpha_{\sigma(k)}} dy_{k} \sign(\sigma|_{F_\mu}) , \prod_{k \in F_{\lambda.\sigma}} \bar{\rho}_{k}^{r_{\sigma(k)}+1}\Big\rceil.
\]
Here, $\sign(\sigma|_{F_\mu})$ is $\pm 1$ depending on the parity of the number of transpositions needed to put the list $(\sigma(i))_{i \in F_\mu}$ in order.

The action on the right-hand side is similar, so one sees that the Krawitz map is compatible with the action, up to the signs.
\end{proof}

\section{Mirror Symmetry for Nondiagonal Groups} \label{sec:ms-non-diag}

\subsection{Extending the definition of the transpose group}
Now we are ready to extend BHK mirror symmetry to nonabelian groups. As before, we take an invertible polynomial $W$ and a group $G=S\ltimes H$ as in Section~\ref{sec:semi}. Recall $H\subseteq G^{diag}_W$ and $S\subseteq \Sigma$ is a group permuting the variables with $S$ a subgroup of the normalizer of $H$.

We observed in Section~\ref{sec:compat} that a permutation $\sigma$ is a symmetry of $W$ if and only if it is a symmetry of $W^T$. Furthermore, one can check that if $S$ is in the normalizer of $H$, then we can also regard $S$ as a subgroup of the normalizer of $H^T$. Hence we may define the \textit{nonabelian dual group}
\[
G^\vee:=S\ltimes H^T.
\]
In this definition, $H^T$ refers to the BHK dual of $H$, as discussed previously.

\begin{example}
As an example, consider again the example above with $W=x_1^4+x_2^4+x_3^4+x_4^6$, but this time with $H=\langle j_W \rangle$ and $S=\langle (123)\rangle$. As before, we set $G=S\rtimes H$. So for a mirror dual model, we have $H^T=\SL_W$. So we obtain the pair $(W^\vee,G^\vee)$, with $W^\vee=W$ and $G^\vee=\langle (123)\rangle\rtimes \SL_W$. 
\end{example}

The first indication that this does indeed describe mirror symmetry would be to show that we have a bidegree-preserving isomorphism
\[
\mathcal{A}_{W,G}\simeq \mathcal{B}_{W^T,G^T}. 
\]
However, as seen in \cite{EG2} and \cite{PWW}, this is not true for every choice of symmetry group. In fact, we can see from \cite{M} that sometimes there is no bidegree-preserving isomorphism. In order to have a bidegree-preserving isomorphism, we must impose some conditions on $S$, which we discuss in Sections~\ref{sec:strongPC} and \ref{sec:EV}. 

Some work has been done in exhibiting a mirror map between state spaces, as in the following theorem. Notice that this theorem only gives a partial mirror map, and only for a restricted class of polynomials. Furthermore, the theorem requires that $S$ contain only even permutations.

\begin{theorem}[see \cite{PWW}] Let $W$ be an invertible Fermat polynomial and $G\subseteq G_W^{max}$ be an admissible group of the form $S \ltimes H$, where $S\subseteq \Sigma$ is the subgroup of even permutations and $H\leq G^{diag}_W$. Define $\mathcal{A}_0\subseteq\mathcal{A}_{W,G}$ and $\mathcal{B}_0\subseteq\mathcal{B}_{W^\vee,G^\vee}$ to be the identity sectors for the A- and B-models, respectively. Let $nar'\leq H$ be the set of narrow diagonal symmetries. We also denote $nar'\leq H^T$ to be the corresponding set on the B side. Then there exist bigraded vector space isomorphisms
	\begin{align*}
		\mathcal{A}_0\xrightarrow{\sim}\mathcal{B}_{nar'} \text{ and } \mathcal{A}_{nar'}\xrightarrow{\sim}\mathcal{B}_0.
	\end{align*}
\end{theorem}

\subsection{The polynomial \texorpdfstring{$W^\sigma$}{W sigma} and its Milnor ring} \label{sec:Ws}
Our mirror map for nonabelian LG models makes use of Krawitz's mirror map \eqref{eq:kmm}, but for a polynomial $W^\sigma$ in place of $W$. The polynomial $W^\sigma$ will be (up to a diagonal rescaling) equal to $W_\sigma$. More generally, we will explain how we can view $W_{\sigma\lambda}$ as the restriction of $W^\sigma$ to the fixed locus of a diagonal symmetry $\beta_\sigma(\lambda)$. 

For a vector space $V$, we write $V_g$ for the fixed locus of the group element $g$ acting on $V$. We denote by $\diag(V, \{\mathbf{b}_i\})$ the group of diagonal (with respect to the basis $\{\mathbf b_i\}$) matrices acting on $V$. We write $V^\vee$ for $\Hom(V, \CC)$.

For each $\sigma$, we decompose it into disjoint cycles $\sigma = \tau_1 \cdots \tau_{N_\sigma}$. There is a map $\beta_\sigma:\diag(\CC^N, \{\mathbf e_i\}) \rightarrow \diag(\CC^N_\sigma, \{\mathbf f_{k,\sigma}\})$, defined by
\begin{equation}
    \beta_\sigma(\lambda).\mathbf{f}_{k,\sigma} = \left(\prod_{i \in \mathcal O_{\tau_k}} \lambda_i\right) \mathbf{f}_{k,\sigma} \label{def-beta}
\end{equation}
and a map $\gamma_\sigma:  \diag(\CC^N, \{\mathbf e_i\}) \rightarrow \diag(\CC^N_\sigma, \{\mathbf f_{k,\sigma}\})$, defined by
\begin{equation}
    \gamma_\sigma(\lambda).\mathbf f_{k,\sigma} = \lambda_{s(k)} \mathbf f_{k,\sigma}.
\end{equation}
Recall $s(k)$ is the smallest index in the $k$-th cycle of $\sigma$ (see Section~\ref{sec:fixed-basis}).  

Notice that $\mathbf f_{k, \sigma}$ is fixed by $\beta(\lambda)$ if and only if  $k \in F_{\sigma\lambda}$ (that is, $\tau_k \lambda^{(k)}$ is a special weighted cycle). 
Hence there is an isomorphism $\psi_{\sigma\lambda}:(\CC^N_{\sigma})_{\beta(\lambda)}\to \CC^N_{\sigma\lambda}$, via
\begin{equation} 
 \psi_{\sigma\lambda}:  \mathbf f_{k, \sigma} \mapsto \mathbf f_{k,\sigma \lambda} \text{ for } k \in F_{\sigma\lambda}. \label{eq:fixed_iso}
\end{equation}
We will drop the indices on $\beta_\sigma$, $\gamma_\sigma$, and $\psi_{\sigma\lambda}$ when they are clear from the context.

It follows from Lemma \ref{lem:mu-f} that $\psi$ is $C_H(\sigma)$ equivariant with the action induced by $\gamma$. That is, for $\mathbf{x}\in (\CC^N_\sigma)_{\beta(\lambda)}$ and $\mu \in C_H(\sigma)$, we have 
\begin{equation}
	\mu.\psi(\mathbf{x}) = \psi(\gamma(\mu).\mathbf{x}),
\end{equation}
or, rewriting in terms of the dual, if $y \in (\CC^N_{\sigma \lambda})^\vee$, we have
\begin{equation}
    (\psi^* y) . \gamma(\mu) = \psi^*(y.\mu). \label{eq:dual-j-ev}
\end{equation}

Now, $\psi$ induces an isomorphism 
\begin{equation}
	\psi^*:\mathcal Q_{W_{\sigma\lambda}} \cdot \det((\CC^N_{\sigma\lambda})^\vee) \cong \mathcal Q_{\psi^*W_{\sigma\lambda}} \cdot \det((\CC^N_\sigma)^\vee_{\beta(\lambda)}).  \label{eq:j-milnor}
\end{equation}

As an intermediate step to our mirror map, we need to express $\psi^* W_{\sigma \lambda}$ in our coordinates $\{z_{k, \sigma}\}$ for $(\CC^N_\sigma)_{\beta(\lambda)}$.
We will write $W = \sum_{i=1}^N M_i$, where $M_i$ are the monomials of $W$. Recall that since $W$ is invertible, every monomial of $W$ is of the form $x_i^a x_j$ or $x_i^a$, and we have taken the convention that $M_i=x_i^a x_j$ or $M_i=x_i^a$. Thus, $M_i.\sigma = M_{\sigma^{-1}(i)}$.

Define a polynomial on $\CC^N_\sigma$:
\begin{equation}
    W^{\sigma} = \sum_{k =1}^{N_\sigma} M_{s(k)}( x_i \mapsto z_{r(i), \sigma} ). \label{eq:def-W-sig}
\end{equation}
Here we make use of some nonstandard notation: the parenthesized expression means to make the indicated substitution into the stated monomial. We see that $W^\sigma$ is an invertible polynomial. 

The following lemma provides an important fact about the symmetries of $W^\sigma$. 

\begin{proposition} \label{prop:beta-sym}
   If $\lambda$ and $\sigma$ are symmetries of $W$, then $\beta(\lambda)$ is a symmetry of $W^\sigma$.
\end{proposition}
\begin{proof}
We see by \eqref{def-beta} that the coefficient of the (diagonal) action of $\beta(\lambda)$ on $z_{k,\sigma}$ is 
\begin{equation}
 \prod_{i \in \cO_{\tau_k}} \lambda_i. \label{eq:coefofaction}
\end{equation}
Let $m = |\cO_{\tau_{k}}|$, the size of the $\sigma$ orbit of $k$. Then \eqref{eq:coefofaction} is the same as the coefficient of the (diagonal) action of $\prod_{i=0}^{m-1} \lambda.\sigma^i$ on any $x_j$ with $j \in \cO_{\tau_k}$.
The monomials of $W^\sigma$ are all monomials $W$ with each $x_j$ replaced with $z_{k,\sigma}$, with $j \in \cO_{\tau_k}$. Since the monomial of $W$ is invariant under $H$, which includes $\prod_{i=0}^{m-1} \lambda.\sigma^i$, it follows that the corresponding monomial of $W^\sigma$ is invariant under $\beta(\lambda)$.
\end{proof}

\begin{lemma} \label{lem:resW}
The functions $\psi^* W_{\sigma\lambda}$ and $(W^\sigma)_{\beta(\lambda)}$ on $(\CC^N_\sigma)_{\beta(\lambda)}$ are the same up to a diagonal (with respect to the basis  $\{\mathbf f_{k,\sigma} \}_{k\in F_{\sigma \lambda}}$) rescaling.

Hence, we get an isomorphism
\[
    \cQ_{\psi^*W_{\sigma\lambda}} \cdot \det((\CC^N_\sigma)^\vee_{\beta(\lambda)}) \cong \cQ_{(W^\sigma)_{\beta(\lambda)}} \cdot \det((\CC^N_\sigma)^\vee_{\beta(\lambda)}).
\]

\end{lemma}

\begin{proof}
    
   We will evaluate both sides on an arbitrary element of $(\CC^N_\sigma)_{\beta(\lambda)}$,  
\[
z=\sum_{k \in F_{\sigma\lambda}} z_k \mathbf f_{k,\sigma}.
\]
The results will be polynomials in the $z_k$. We will see that the two polynomials have the same monomials, differing only in their coefficients. They also have the same number of monomials as variables, so the result will follow.

 We have $z_{k,\sigma}(z) = z_k$ if $k \in F_{\sigma\lambda}$, and $z_{k,\sigma}(z) = 0$ if $k \notin F_{\sigma\lambda}$. So we have
 \[
W^{\sigma}(z) = \sum_{k=1}^{N_\sigma} M_{s(k)}\left(x_i \mapsto \begin{cases} z_{r(i)} & r(i) \in F_{\sigma\lambda} \\
        0 & \text{ otherwise} \end{cases}\right).
 \] 
 
  For the other side, one can check from \eqref{def:f} that $x_i(\mathbf f_{k,\sigma\lambda}) = 0$ if and only if $r(i) =k$. 
  Hence, there are constants $c_{r(i),\sigma\lambda}$ so that
  \[
    x_i \circ \psi (z) = x_i\left( \sum_{k \in F_{\sigma \lambda}} z_k \mathbf f_{k,\sigma \lambda} \right) = \begin{cases} c_{r(i),\sigma\lambda} z_{r(i)} & r(i) \in F_{\sigma \lambda} \\
    0 & \text{ otherwise} \end{cases}.
  \]
 
Next we note that if $M_i$ is the $i$-th monomial of $W$, then $M_i.\sigma\lambda = M_{\sigma^{-1}(i)}$. Also, $\psi(z)$ is $\sigma\lambda$ invariant. It follows that $M_i \circ \psi(z) = M_i(\sigma\lambda\cdot\psi(z)) = M_{\sigma^{-1}(i)}\circ \psi(z)$. Hence, we have
\begin{align*}
   \psi^*W(z) = W \circ \psi (z) &= \sum_{i=1}^N M_i \circ \psi (z) \\
                  &= \sum_{k =1}^{N_\sigma} \sum_{\ell=0}^{|\tau_k|-1} M_{\sigma^{-\ell}(s(k))} \circ \psi(z) \\
                  &= \sum_{k =1}^{N_\sigma} |\tau_k| M_{s(k)} \circ \psi(z) \\
                  &= \sum_{k =1}^{N_\sigma} |\tau_k| M_{s(k)}\left( x_i \mapsto \begin{cases} c_{r(i),\sigma\lambda} z_{r(i)} & r(i) \in F_{\sigma \lambda} \\
    0 & \text{else} \end{cases}\right),
\end{align*}
so it is the same as $W^\sigma(z)$, up to the scalar on each monomial.
\end{proof}
Combining \eqref{eq:j-milnor} and Lemma \ref{lem:resW}, we obtain the following.

\begin{corollary} \label{cor:psi}
$\psi^*$ gives a $C_H(\sigma)$ equivariant (via $\gamma$) isomorphism
    \begin{equation}
	\psi^*: \mathcal Q_{W_{\sigma\lambda}} \cdot \det(\CC^N_{\sigma\lambda}) \cong \mathcal Q_{(W^\sigma)_{\beta(\lambda)}} \cdot \det((\CC^N_\sigma)_{\beta(\lambda)}).  \label{eq:milnor-iso}
\end{equation}
\end{corollary}

\subsection{\texorpdfstring{$W^\sigma$}{W sigma} and exponent matrices}

We will define matrices $B_\sigma$ and $C_\sigma$ so that
\begin{equation}
 [\beta_\sigma(\lambda)] = B_\sigma [\lambda], \;\;\; [\gamma_\sigma(\lambda)] = C_\sigma[\lambda].
\end{equation}
As before, we will drop the subscripts when $\sigma$ is clear from context.

Let $\sigma = \tau_1 \cdot \cdots \cdot \tau_{N_\sigma}$ be the cycle decomposition.  Define $N_\sigma \times n$ matrices
\[
 B_{ij} = \begin{cases}
1 & \text{if } j \in \cO_{\tau_i}\\
0 & \text{ otherwise }
\end{cases}
\]
and 
\[
C_{ij} = \begin{cases}
1 & \text{ if } j = \min \cO_{\tau_i}\\
0 & \text{ otherwise}
\end{cases}.
\]
\begin{example}
 Let $\sigma = (12)(345)$. Then 
 \[
    B = \begin{bmatrix}
    1 & 1 & 0 & 0 & 0 \\
    0 & 0 & 1 & 1 & 1
    \end{bmatrix}
 \]
 and
 \[
    C = \begin{bmatrix}
    1 & 0 & 0 & 0 & 0 \\
    0 & 0 & 1 & 0 & 0
    \end{bmatrix}.
 \]
\end{example}

Notice that multiplying $B$ on the left of a (column) vector has the effect of adding together the entries of the vector corresponding to the indices in each cycle. Furthermore, mulitplication by $C$ just picks out the entry corresponding to the smallest number in each cycle.

We also observe that Lemma~\ref{lem:CHs} can be restated as
\begin{equation}
    \lambda\sigma = \sigma \lambda \Leftrightarrow [\lambda] = B^T\ell \text{ for some } \ell \in (\QQ/\ZZ)^k.
\end{equation}

The following properties will be useful.

\begin{lemma} \label{lem:matrix-prop}
\mbox{}
\begin{enumerate}
    \item $BC^T =  I = CB^T$.
    \item $A_{W^\sigma}B = B A_W$.
    \item $B^T A_{W^\sigma} = A_W B^T$.
    
\end{enumerate}
\end{lemma}

\begin{corollary}
    The exponent matrix of $W^\sigma$ is $A_{W^\sigma} = B A_W C^T = C A_W B^T$.
\end{corollary}

\begin{corollary} \label{cor:Tsig}
   $(W^T)^\sigma = (W^\sigma)^T$. 
\end{corollary} 

\begin{proof}[Proof of Lemma \ref{lem:matrix-prop}]
    For (1), we have
    \[
        (BC^T)_{ij} = \sum_k B_{ik}C_{jk}. 
    \]
    A product in the sum above is nonzero if and only if $k \in \cO_{\tau_i}$ and $k = \min \cO_{\tau_j}$. So we must have $i=j$, and there is exactly one value of $k$ with $k=\min \cO_{\tau_i}$. So $(BC^T)_{ij} = 1$ if $i=j$ and $0$ otherwise, as desired. The right-hand side is just the transpose of the left-hand side.

     When one multiplies an exponent matrix $A_W$ on the right by a matrix $D$ with integer entries, the result is the exponent matrix of the polynomial obtained by the substitution $x_i \mapsto \prod_j z_j^{D_{ij}}$.
     
     When one multiplies an exponent matrix $A_W$ on the left by $D$, the result is the exponent matrix of the polynomial $\sum_i \prod_j M_j^{D_{ij}}$, where $M_j$ is the $j$th monomial of $W$. 
     
     We make one more observation: if $x_i^a x_j$ is a monomial of $W$, then $|\cO_\sigma(i)| = |\cO_\sigma(j)|$. This follows from the fact that there cannot also be a monomial $x_k^b x_j$ for $k\neq i$.
     
     Hence the left-hand side of (2) corresponds to the polynomial 
     \begin{align*}
         &\sum_{k=1}^{N_\sigma} M_{s(k)}\left( x_i \mapsto z_{r(i),\sigma} \mapsto \prod_{j \in \cO_{\tau_{r(i)}}} x_{j}\right)\\
        = &\sum_{k=1}^{N_\sigma} M_{s(k)}\left( x_i \mapsto \prod_{\ell=0}^{ |\cO_{\sigma}(i)|-1} x_{\sigma^{\ell}(i)}\right)\\
        =& \sum_{k=1}^{N_\sigma} \prod_{\ell = 0}^{|\cO_\sigma(k)| - 1} M_{s(k)}(x_i \mapsto x_{\sigma^{\ell}(i)}) \\
        =& \sum_{k=1}^{N_\sigma} \prod_{\ell = 0}^{|\cO_\sigma(k)| - 1} M_{\sigma^{\ell}(s(k))} \\
        =&\sum_{k=1}^{N_\sigma} \prod_{j \in \cO_{\tau_k}} M_j,
     \end{align*}
    which is precisely the polynomial corresponding to the right-hand side of (2). The second equality is by the observation above and a property of monomials.
    
    The left-hand side of (3) corresponds to the polynomial
    \begin{align*}
        \sum_{k=1}^N \prod_{j=1}^{N_\sigma} M_{s(j)}(x_i \mapsto z_{r(i),\sigma})^{B_{jk}}.
    \end{align*}
    For fixed $k$, we have $B_{jk} \neq 0$ for exactly one value of $k$, and in this case $r(k) = j$, so this is equal to 
    \[
        \sum_{k=1}^N M_{s(r(k))}(x_i \mapsto z_{r(i),\sigma}).
    \]
    The substitution for the right-hand side of (3) is
    \[
        x_k \mapsto \prod_{j=1}^{N_\sigma} z_{j,\sigma}^{B_{jk}},
    \]
    which, for similar reasons, is
    \[
        x_k \mapsto z_{r(k),\sigma}.
    \]
    
    Recall that for each $i \in \{1, \dots, N\}$, $s(r(i))$ is the smallest number in the $\sigma$ orbit of $i$. Let $\ell(i)$ be the number of times we must apply $\sigma$ to the smallest number in the cycle to obtain $i$: 
\begin{equation}
    i = \sigma^{\ell(i)}(s(r(i))). \label{eq:elli}
\end{equation} 
For example, if $\sigma=(123)(4)$, then $s(r(2))=1$ and $\ell(2)=1$, but $\ell(3)=2$. 
    
    Then the right-hand side of (3) is 
    \begin{align*}
        \sum_{k=1}^N M_k( x_i \mapsto z_{r(i),\sigma} ) &= \sum_{k=1}^N M_{\sigma^{\ell(k)}(s(r(k)))} ( x_i \mapsto z_{r(i),\sigma} )\\
        &=  \sum_{k=1}^N M_{s(r(k))} ( x_i \mapsto x_{\sigma^{\ell(k)}(i)} \mapsto z_{r(\sigma^{\ell(k)}(i)),\sigma} )\\
        &= \sum_{k=1}^N M_{s(r(k))} (x_i \mapsto z_{r(i), \sigma}).
    \end{align*}

\end{proof}

\subsection{An example illustrating the idea of the new mirror map}
With all of the pieces now in place, we can describe the mirror map. We will first do so through an illustrative example. Our strategy is roughly as follows. 

We use the maps $\psi^*_{\sigma\lambda}$ to send elements indexed by $\lambda\sigma$ (for a given $\sigma\in S$) of the unprojected state space $\HH_{W,G}^{un}$ to the state space $\HH_{W^\sigma, G^{diag}_{W^\sigma}}^{un}$. Similarly, we map elements of $\HH_{W^T,G^T}^{un}$ indexed by $\lambda'\sigma$ (the same $\sigma$) to $\HH_{(W^T)^\sigma, G^{diag}_{(W^T)^\sigma}}^{un}$. We know that $\HH_{W^\sigma, G^{diag}_{W^\sigma}}^{un}$ and $\HH_{(W^T)^\sigma, G^{diag}_{(W^T)^\sigma}}^{un}$ are related by the Krawitz mirror map, so we will use this to match up our elements. One should note that this does not provide a one-to-one correspondence between basis elements but rather many elements of $\HH_{W,G}^{un}$ will map to a single element of $\HH_{W^\sigma, G^{diag}_{W^\sigma}}^{un}$. We will show later how elements that have the same image are related. So under certain circumstances, we will get an isomorphism on the invariant subspace $\HH_{W,G}$. 
We will spell the details out more precisely later, but for now the following example may illustrate the idea.

Let $W=x_1^4+x_2^4+x_3^4+x_4^6$. Then $j_W=(\frac{1}{4},\frac{1}{4},\frac{1}{4},\frac{1}{6})$. We define $G=\langle j_W,(123)\rangle$. We will take $(W,G)$ as the data for our $A$-model.

Then $W=W^T$, and one can check that $H^T$ is the group generated by three elements 
$\alpha = (\frac14,0,-\frac14,0)$, $\beta = (0,\frac14,-\frac14,0)$, and $\gamma = (0,0,\frac12,\frac12)$, as in the example of Section~\ref{sec:state-space-ex}. We will take this example as our B-model. We can construct the A-model state space in a similar manner, which we will not describe in detail. 

In Table~\ref{tab:ms-example}, we list a basis for $\HH_{W,G}^{un}$ and $\HH_{W^\vee,G^\vee}^{un}$ together and show how they correspond under the mirror map. Elements in adjacent rows of the same color are in the same $S$ orbit, and elements in a horizontal list are in the same $H$ orbit  (up to scalars). The bidgree pair is the A-model bidegree for the A-model and the B-model bidegree for the B-model. We write the images of A- and B-model elements under $\psi$ in the columns labeled $\psi_A$ and $\psi_B$. We don't write the image of $\psi$ if $\psi$ is trivial. A- and B-model elements $x$ and $y$ in the same row have $\psi_A(x)$ and $\psi_B(y)$ corresponding under the Krawitz isomorphism. 
\begin{table}[] 
\begin{center}
\begin{tabular}{c|c|c|c||c}
    \textbf{A-model basis element}  & $\psi_A$ & $\psi_B$ & \textbf{B-model basis element} & \textbf{Bidegree}\\
    \hline
    $\lfloor x_1x_2x_3x_4^2; \text{id}\rceil$ & & & $\lfloor 1; (\frac{1}{2},\frac{1}{2},\frac{1}{2},\frac{1}{2})\rceil$& $(\frac{13}{12},\frac{13}{12})$\\
    \rowcolor{TableShade}
     $\lfloor x_1x_2^2x_4^2; \text{id}\rceil$  &&& $\lfloor 1; (\frac{1}{2},\frac{3}{4},\frac{1}{4},\frac{1}{2})\rceil$& $(\frac{13}{12},\frac{13}{12})$\\ 
     \rowcolor{TableShade}
     $\lfloor x_2x_3^2x_4^2; \text{id}\rceil$ &  && $\lfloor 1; (\frac{1}{4},\frac{1}{2},\frac{3}{4},\frac{1}{2})\rceil$ &\\
     \rowcolor{TableShade}
     $\lfloor x_3x_1^2x_2^2; \text{id}\rceil$   &&& $\lfloor 1; (\frac{3}{4},\frac{1}{4},\frac{1}{2},\frac{1}{2})\rceil$ &\\
     $\lfloor x_1^2x_2x_4^2; \text{id}\rceil$  &&& $\lfloor 1; (\frac{3}{4},\frac{1}{2},\frac{1}{4},\frac{1}{2})\rceil$& $(\frac{13}{12},\frac{13}{12})$\\ 
     $\lfloor x_2^2x_3x_4^2; \text{id}\rceil$ &&& $\lfloor 1; (\frac{1}{4},\frac{3}{4},\frac{1}{2},\frac{1}{2})\rceil$ &\\
     $\lfloor x_3^2x_1x_4^2; \text{id}\rceil$ &&& $\lfloor 1; (\frac{1}{2},\frac{1}{4},\frac{3}{4},\frac{1}{2})\rceil$ &\\
     \rowcolor{TableShade}
     $\lfloor 1; j_W\rceil$  &&& $\lfloor 1; \text{id}\rceil$& $(0,0)$\\
     $\lfloor 1; j_W^2\rceil$  &&& $\lfloor x_1x_2x_3x_4; \text{id}\rceil$& $(\frac{11}{12},\frac{11}{12})$\\
    \rowcolor{TableShade}
     $\lfloor 1; j_W^3\rceil$  &&& $\lfloor x_1^2x_2^2x_3^2x_4^2; \text{id}\rceil$& $(\frac{11}{6},\frac{11}{6})$\\
     $\lfloor 1; j_W^5\rceil$ &&& $\lfloor x_4^4; \text{id}\rceil$ & $(\frac{2}{3},\frac{2}{3})$\\
    \rowcolor{TableShade}
     $\lfloor 1; j_W^7\rceil$ &&&  $\lfloor x_1^2x_2^2x_3^2; \text{id}\rceil$ & $(\frac{3}{2},\frac{3}{2})$\\
     $\lfloor 1; j_W^9\rceil$ &&& $\lfloor x_4^2; \text{id}\rceil$ & $(\frac{1}{3},\frac{1}{3})$\\
    \rowcolor{TableShade}
     $\lfloor 1; j_W^{10}\rceil$ &&& $\lfloor x_1x_2x_3x_4^3; \text{id}\rceil$ & $(\frac{4}{3},\frac{4}{3})$\\
     $\lfloor 1; j_W^{11}\rceil$ &&& $\lfloor x_1^2x_2^2x_3^2x_4^4; \text{id}\rceil$ & $(\frac{13}{12},\frac{13}{12})$\\
     \rowcolor{TableShade}
     $\lfloor x_1; j_W^4\rceil$ &&& $\lfloor x_4^3; (\frac{1}{2},\frac{1}{4},\frac{1}{4},0)\rceil$ & $(\frac{3}{4},\frac{7}{4})$\\
     \rowcolor{TableShade}
     $\lfloor x_2; j_W^4\rceil$ &&& $\lfloor x_4^3; (\frac{1}{4},\frac{1}{2},\frac{1}{4},0)\rceil$ &\\
     \rowcolor{TableShade}
     $\lfloor x_3; j_W^4\rceil$ &&& $\lfloor x_4^3; (\frac{1}{4},\frac{1}{4},\frac{1}{2},0)\rceil$ &\\
     $\lfloor x_1^2x_2^2x_3; j_W^4\rceil$ &&& $\lfloor x_4^3; (\frac{3}{4},\frac{3}{4},\frac{1}{2},0)\rceil$ & $(\frac{7}{4},\frac{3}{4})$\\
     $\lfloor x_1^2x_2x_3^2; j_W^4\rceil$ &&& $\lfloor x_4^3; (\frac{1}{2},\frac{3}{4},\frac{3}{4},0)\rceil$ &\\
     $\lfloor x_1x_2^2x_3^3; j_W^4\rceil$ &&& $\lfloor x_4^3; (\frac{3}{4},\frac{1}{2},\frac{3}{4},0)\rceil$ &\\
     \rowcolor{TableShade}
     $\lfloor x_1; j_W^8\rceil$ &&& $\lfloor x_4; (\frac{1}{4},\frac{1}{4},\frac{1}{2},0)\rceil$  & $(\frac{5}{12},\frac{7}{12})$\\
     \rowcolor{TableShade}
     $\lfloor x_2; j_W^8\rceil$ &&& $\lfloor x_4; (\frac{1}{4},\frac{1}{2},\frac{1}{4},0)\rceil$ &\\
     \rowcolor{TableShade}
     $\lfloor x_3; j_W^8\rceil$ &&& $\lfloor x_4; (\frac{1}{2},\frac{1}{4},\frac{1}{4},0)\rceil$ &\\
     $\lfloor x_1^2x_2^2x_3; j_W^8\rceil$&&& $\lfloor x_4; (\frac{3}{4},\frac{3}{4},\frac{1}{2},0)\rceil$  & $(\frac{7}{12},\frac{5}{12})$  \\
     $\lfloor x_1^2x_2x_3^2; j_W^8\rceil$ &&&$\lfloor x_4; (\frac{3}{4},\frac{1}{2},\frac{3}{4},0)\rceil$ &\\
     $\lfloor x_1x_2^2x_3^3; j_W^8\rceil$ &&& $\lfloor x_4; (\frac{1}{2},\frac{3}{4},\frac{3}{4},0)\rceil$ &\\
     \rowcolor{TableShade}
     $\lfloor 1; (123)j_W\rceil$ & $\lfloor 1; (\frac34,\frac16)\rceil$ & $\lfloor z_1^2; \text{id} \rceil$ & $\lfloor z_1^2;(123)\rceil, \dots (\text{15 others})$ & $(1,1)$\\
     $\lfloor 1; (123)j_W^2\rceil$ & $\lfloor 1; (\frac24, \frac26)\rceil$ & $\lfloor z_1z_4; \text{id} \rceil$ & $\lfloor z_1z_4;(123)\rceil, \dots (\text{15 others})$ & $(\frac{11}{12},\frac{11}{12})$\\
     \rowcolor{TableShade}
     $\lfloor 1; (123)j_W^3\rceil$ & $\lfloor 1; (\frac14,\frac36)\rceil$ & $\lfloor z_4^2; \text{id} \rceil$ & $\lfloor z_4^2;(123)\rceil, \dots (\text{15 others})$ & $(\frac{5}{6},\frac{5}{6})$\\
     $\lfloor 1; (123)j_W^5\rceil$ &$\lfloor 1; (\frac34,\frac56)\rceil$ & $\lfloor z_1^2z_4^4; \text{id} \rceil$ & $\lfloor z_1^2z_4^4;(123)\rceil, \dots (\text{15 others})$  & $(\frac{5}{3},\frac{5}{3})$\\
     \rowcolor{TableShade}
     $\lfloor 1; (123)j_W^7\rceil$ &$\lfloor 1; (\frac14,\frac16)\rceil$ & $\lfloor 1; \text{id} \rceil$ & $\lfloor1; (123)\rceil, \dots (\text{15 others})$ & $(\frac{1}{2},\frac{1}{2})$\\
     $\lfloor 1; (123)j_W^9\rceil$ &$\lfloor 1; (\frac34,\frac36)\rceil$ & $\lfloor z_1^2z_4^2; \text{id} \rceil$ & $\lfloor z_1^2 z_4^2;(123)\rceil, \dots (\text{15 others})$ & $(\frac{4}{3},\frac{4}{3})$\\
     \rowcolor{TableShade}
     $\lfloor 1; (123)j_W^{10}\rceil$ &$\lfloor 1; (\frac24,\frac46)\rceil$ & $\lfloor z_1z_4^3; \text{id} \rceil$ & $\lfloor z_1z_4^3;(123)\rceil, \dots (\text{15 others})$ & $(\frac{5}{4},\frac{5}{4})$\\
     $\lfloor 1; (123)j_W^{11}\rceil$ &$\lfloor 1; (\frac14,\frac56)\rceil$ & $\lfloor z_4^4; \text{id} \rceil$ & $\lfloor z_4^4;(123)\rceil, \dots (\text{15 others})$ & $(\frac{7}{6},\frac{7}{6})$\\
     \rowcolor{TableShade}
     $\lfloor 1; (132)j_W\rceil$ & $\lfloor 1; (\frac34,\frac16)\rceil$ & $\lfloor z_1^2; \text{id} \rceil$ & $\lfloor z_1^2;(132)\rceil, \dots (\text{15 others})$ & $(1,1)$\\
     $\lfloor 1; (132)j_W^2\rceil$ & $\lfloor 1; (\frac24, \frac26)\rceil$ & $\lfloor z_1z_4; \text{id} \rceil$ & $\lfloor z_1z_4;(132)\rceil, \dots (\text{15 others})$ & $(\frac{11}{12},\frac{11}{12})$\\
     \rowcolor{TableShade}
     $\lfloor 1; (132)j_W^3\rceil$ & $\lfloor 1; (\frac14,\frac36)\rceil$ & $\lfloor z_4^2; \text{id} \rceil$ & $\lfloor z_4^2;(132)\rceil, \dots (\text{15 others})$ & $(\frac{5}{6},\frac{5}{6})$\\
     $\lfloor 1; (132)j_W^5\rceil$ &$\lfloor 1; (\frac34,\frac56)\rceil$ & $\lfloor z_1^2z_4^4; \text{id} \rceil$ & $\lfloor z_1^2z_4^4;(132)\rceil, \dots (\text{15 others})$ & $(\frac{5}{3},\frac{5}{3})$\\
     \rowcolor{TableShade}
     $\lfloor 1; (132)j_W^7\rceil$ &$\lfloor 1; (\frac14,\frac16)\rceil$ & $\lfloor 1; \text{id} \rceil$ & $\lfloor1; (132)\rceil, \dots (\text{15 others})$ & $(\frac{1}{2},\frac{1}{2})$\\
     $\lfloor 1; (132)j_W^9\rceil$ &$\lfloor 1; (\frac34,\frac36)\rceil$ & $\lfloor z_1^2z_4^2; \text{id} \rceil$ & $\lfloor z_1^2 z_4^2;(132)\rceil, \dots (\text{15 others})$ & $(\frac{4}{3},\frac{4}{3})$\\
     \rowcolor{TableShade}
     $\lfloor 1; (132)j_W^{10}\rceil$ &$\lfloor 1; (\frac24,\frac46)\rceil$ & $\lfloor z_1z_4^3; \text{id} \rceil$ & $\lfloor z_1z_4^3;(132)\rceil, \dots (\text{15 others})$ & $(\frac{5}{4},\frac{5}{4})$\\
     $\lfloor 1; (132)j_W^{11}\rceil$ &$\lfloor 1; (\frac14,\frac56)\rceil$ & $\lfloor z_4^4; \text{id} \rceil$ & $\lfloor z_4^4;(132)\rceil, \dots (\text{15 others})$ & $(\frac{7}{6},\frac{7}{6})$\\
     \rowcolor{TableShade}
     $\lfloor z_1z_4^2; (123)\rceil$ & $\lfloor z_1z_4^2; \text{id} \rceil$ & $\lfloor1; (\frac12, \frac36) \rceil$ &$\lfloor 1;(\frac{1}{2},0,0,\frac{1}{2})(123)\rceil, \dots (\text{15 others})$ & $(\frac{13}{12},\frac{13}{12})$\\
     $\lfloor z_1z_4^2; (132)\rceil$ & $\lfloor z_1z_4; \text{id} \rceil$ &$\lfloor1; (\frac12, \frac36) \rceil$ &  $\lfloor 1;(\frac{1}{2},0,0,\frac{1}{2})(132)\rceil, \dots (\text{15 others})$ & $(\frac{13}{12},\frac{13}{12})$\\
\end{tabular}
\end{center}
    \caption{An example of the mirror map for a nonabelian group.}
    \label{tab:ms-example}
\end{table}

For example, consider the first row with an entry in the second and third columns. On the A side there is a single basis element indexed by $g=(123)j_w$. Recall $j_W=(\tfrac{1}{4},\tfrac{1}{4},\tfrac{1}{4},\tfrac{1}{6})$. In the notation of the previous section, we have $W^{\sigma}=z_1^4+z_2^6$. Then $G_{W^\sigma}^{diag}$ is generated by $(\tfrac{1}{4}, 0)$ and $(0,\tfrac{1}{6})$. Furthermore, we have $B_\sigma=\begin{bmatrix}
    1 & 1 & 1 & 0\\ 0 & 0 & 0 & 1\end{bmatrix}$, and so $\beta_\sigma(j_w)=(\tfrac{3}{4}, \tfrac{1}{6})$. 

Both domain and codomain are trivial, so the map $\psi_g:\Fix g\to \Fix(\beta_\sigma(j_w))$ is trivial. Hence we have $\psi_A(\lfloor 1; (123)j_w\rceil)=\lfloor 1; (\tfrac{3}{4}, \tfrac{1}{6})\rceil$. This completes the description of the first two columns of this row. 

On the B side in the same row, we can consider $\lfloor z_1^2; (123)\rceil$, i.e., $\lambda$ is the identity.  As in Corollary~\ref{cor:Tsig}, we have $(W^\sigma)^T=(W^T)^\sigma=z_1^4+z_2^6$. The matrix $B_\sigma$ is the same as on the A side, and so $\beta_\sigma(\lambda)=(0,0)$. This time, we have the map $\psi_g:\Fix g\to \Fix((0,0))$ given by $\psi_g(1,1,1,0)= (1,0)$ and $\psi_g(0,0,0,1)=(0,1)$, which induces the isomorphism of corresponding Milnor rings of Corollary~\ref{cor:psi}.

So we obtain $\psi_B(\lfloor z_1^2; (123)\rceil)=\lfloor z_1^2; (0,0)\rceil$. Again, we have suppressed the volume form.

In the same row, we also have elements of the form $\lfloor z_1^2; g\rceil$, where $g=(123)\lambda$ and $\lambda=(\theta_1,\theta_2,\theta_3,0)$ with $\sum_{i=1}^3\theta_i\equiv 0\pmod{\ZZ}$. One can check that there are 16 such elements (one of them where $\lambda$ is the identity), and that all such elements of $G^\vee$ are conjugate by an element of $H^T=\SL_W$. A similar analysis shows us that $\psi_B(\lfloor z_1^2; g\rceil)=\lfloor z_1^2; (0,0)\rceil$ as in Table~\ref{tab:ms-example}. 

The last step in this example is to notice that $\lfloor 1; (\tfrac{3}{4}, \tfrac{1}{6})\rceil$ and $\lfloor z_1^2; (0,0)\rceil$ correspond under Krawitz's mirror map.  

To pass to the (projected) state space, one needs to identify $G$ invariants with $G$ orbits, and the same for $G^\vee$. We will describe this in the next section. In our example in Table~\ref{tab:ms-example}, elements on the A side sharing a shaded region form a $G$ orbit, and each $G$ orbit corresponds to an element of a basis for the state space $\HH_{W,G}$. A similar statement holds for $G^\vee$. Assuming this fact, Table~\ref{tab:ms-example} demonstrates the mirror map for this first example.

\subsection{A basis for the unprojected state space}
For each $\sigma \in S$, consider the unprojected maximal state space (for diagonal symmetries)
\[
    \HH_{W^\sigma, G^{diag}_{W^\sigma}}^{un} \cong \bigoplus_{\alpha \in G^{diag}_{W^\sigma}} \cQ_{W^\sigma_\alpha} \cdot \omega_\alpha. 
\]
Pick a basis of monomials for each $\cQ_{W^\sigma_\alpha}$, and use these to form a basis $\Ca_\sigma$ of $\HH_{W^\sigma, G^{diag}_{W^\sigma}}^{un}$. This can be done using tensor products of the bases of Theorem~\ref{thm:kreu-basis}, for example. 

For each $c=\lfloor m; \alpha \rceil \in \Ca_\sigma$, let  
\[
 \Ba_c := \left\{\lfloor (\psi_{\sigma\lambda}^{-1})^*m; \sigma \lambda \rceil : \lambda \in H, \beta_\sigma(\lambda) = \alpha\right\}.
\]

Using the isomorphisms of Corollary~\ref{cor:psi}, one can see that 
\begin{equation} \label{eq:defB}
\Ba := \bigcup_{\sigma \in S} \bigcup_{c\in \Ca_\sigma} \Ba_c
\end{equation}
is a basis for the unprojected state space $\HH_{W,G}^{un}$. If we return to the example illustrated in Table~\ref{tab:ms-example}, we might consider $\sigma=(123)$ and $c=\lfloor z_1^2; \id\rceil$, found in the third column. Then we have 
\[
\Ba_c=\{\lfloor z_1^2; (123)\lambda\rceil \mid \lambda \in \left< (1/4,0,-1/4,0), (0,1/4,-1/4,0) \right>\}.
\]
Notice that variable $z_1$ in $c$ is in $\cQ_{W^\sigma}$, while the $z_1$'s in $\Ba_c$ are in $\cQ_{W_{\sigma\lambda}}$. 

One can do the same thing replacing $\sigma$ with $\sigma^2=(132)$ (although notice the dependence of $z_1$ on $\sigma$ is hidden in the notation). Indeed, the element $\lfloor z_1^2; \id\rceil$ appears again later in the table, but this time as an element of $\HH_{W^{\sigma^2}, G^{diag}_{W^{\sigma^2}}}^{un}$.

Using the same constructions, let $\Ca^T_\sigma$ be the basis for $\HH_{(W^T)^\sigma,G_{(W^T)^\sigma}^{diag}}^{un}$ and $\Ba^T$ the basis for $\HH_{W^T,G^T}^{un}$.

By Krawitz \cite{Krawitz}, we have a bidegree-preserving isomorphism 
\[
 \HH_{W^\sigma,G^{diag}_{W^\sigma}}^{un} \cong \HH_{(W^\sigma)^T,G^{diag}_{(W^\sigma)^T}}^{un},
\]
and using Corollary \ref{cor:Tsig}, we get an isomorphism 
\begin{equation} \label{eq:Ciso}
\HH_{(W^\sigma)^T,G^{diag}_{(W^\sigma)^T}}^{un} \cong \HH_{(W^T)^\sigma,G^{diag}_{(W^T)^\sigma}}^{un}.
\end{equation}
In fact, the Krawitz mirror map takes elements of $\Ca_\sigma$ to elements of $\Ca^T_\sigma$.

\subsection{Bidegrees}
The constructions in the previous sections interact nicely with the bidegrees. 

\begin{lemma}
 If $b \in \Ba_c$, with $c \in \Ca_\sigma$, then
 \[
    \bideg(b) = \bideg(c) + (\age \sigma - \age j_W + \age j_{W^\sigma},\age \sigma - \age j_W + \age j_{W^\sigma})
 \]
 for both the A- and B-model bidegrees. In particular, $\bideg(b)-\bideg(c)$ is a constant depending only on $\sigma$. 
\end{lemma}
\begin{proof}
Let $b=\lfloor p; \sigma\lambda \rceil$. Then its A-model bidegree is 
\[
(\deg p + \age \lambda \sigma -\age j_W, N_{\sigma\lambda}-\deg p + \age \sigma\lambda-\age j_W)
\]
and the A-model bidegree of $c=\lfloor \psi(p) ;\beta_\sigma(\lambda) \rceil$ is 
\[
(\deg \psi(p) + \age \beta(\lambda) -\age j_{W^\sigma}, N_{\beta(\lambda)}-\deg \psi(p) + \age \beta(\lambda)-\age j_{W^\sigma}).
\]

Recall that the vector of weights $q_i$ is given by $A_W^{-1} \mathbf 1$. If $\mathbf p$ is the vector of exponents of $p$, then we have
\[
	\deg p = \mathbf p^T A_W^{-1} \mathbf 1,
\]
and the vector of $\psi(p)$ is $B\mathbf p$. Hence, using Lemma~\ref{lem:matrix-prop},
\[
    \deg \psi(p) = \mathbf p^T B^T A_{W^\sigma}^{-1} \mathbf 1 = \mathbf p^T A_W^{-1} B^T \mathbf 1 = \mathbf p^T A_W^{-1} \mathbf 1 = \deg p.
\]
    
By \eqref{eq:fixed_iso}, we have $N_{\sigma \lambda} = N_{\beta(\lambda)}$. 
    
Finally, we claim that 
\begin{equation}
    \age \sigma \lambda = \age \sigma + \age\beta_\sigma(\lambda). \label{eq:age-beta}
\end{equation} 

Indeed, recall from \eqref{eq:age-ts} that for a weighted cycle $\tau \lambda_\tau$ of length $m$ we have that 
\[
\age \tau\lambda = \age\tau + \age\lambda - \sum_{\ell=0}^{m-1} t_\ell,
\]
where $t_\ell = 1$ if $\frac{\ell}{m} + \frac1m \sum_{i=1}^N \Theta^\lambda_i \ge 1$ and $t_\ell=0$ otherwise (recall that since $\tau\lambda$ is a weighted cycle, $\Theta^\lambda_i=0$ unless $i\in \cO_\tau$). This condition simplifies to $\sum \Theta^\lambda_i \ge m- \ell$, so we see that $\sum t_\ell = \lfloor \sum \Theta^\lambda_i \rfloor$. But $\sum \Theta_i - \lfloor \sum \Theta^\lambda_i \rfloor$ is precisely the age of $\beta(\lambda_\tau)$, so we obtain 
\[
\age \tau\lambda =\age \tau +\age \beta(\lambda)
\]
for weighted cycles. 

The age is additive for ``disjoint'' weighted cycles, so we obtain for a general element of $\sigma\lambda\in G$
\[
\age \sigma\lambda =\age \sigma +\age \beta(\lambda).
\]

Making these three substitutions, we obtain the formula for the A-model bidegree in the lemma. 
    
The B-model bidegree of $b$ is
\[
(\deg p + \age \sigma\lambda - \age j_W, \deg p + \age (\sigma \lambda)^{-1} - \age j_W )
\]
while the B-model bidegree of $c$ is
\[
(\deg \psi(p) + \age \beta(\lambda) - \age j_{W^\sigma}, \deg \psi (p)  + \age \beta(\lambda)^{-1} - \age j_{W^\sigma})
\]
Hence, it remains only to check that $\age(\sigma \lambda)^{-1} - \age \beta_\sigma(\lambda)^{-1} = \age \sigma$. As before, it suffices to check this for weighted cycles.
Recall that, in general, we have
\[
    \age g  + \age g^{-1} + N_g = n.
\]

Let $\tau \lambda$ be  a weighted cycle of length $m$. Then we have
\[
\age (\tau \lambda)^{-1}  - \age \beta_\tau(\lambda)^{-1} = (- \age \tau\lambda  - N_{\tau\lambda} + n) - (-\age \beta_\tau(\lambda) - N_{\beta_\tau(\lambda)} + \dim \Fix(\tau)).
\]
The dimension of  $\Fix(\tau)$ is $n-m+1$. Hence, using \eqref{eq:age-beta}, $N_{\tau\lambda} = N_{\beta_\tau(\lambda)}$, and \eqref{eq:age-tau}, the difference simplifies to 
\[
-\age \tau + m-1 = \age \tau.
\]
\end{proof}

\begin{corollary} \label{cor:bidegrees}
If $c \in \Ca$ and $c'\in \Ca^T$ correspond under the Krawitz isomorphism, then the A-model bidegree of elements of $\Ba_c$ is the same as the B-model bidegree of elements of $\Ba^T_{c'}$.   
\end{corollary}
\begin{proof}
    Since the Krawitz isomorphism interchanges the bidegrees, it only remains to recall that $\age j_W = \age j_{W^T}$.
\end{proof}

\subsection{The actions of \texorpdfstring{$G$}{G} on \texorpdfstring{$\Ba$}{B} and \texorpdfstring{$\Ca$}{C}}
We want to define an action of $G$ on $\Ba$ (not on its span), denoted by $\star$, that captures the important information. We first introduce some notation.

Let $\pi \in S$. We define $P_{\sigma, \pi}\in S_{N_\sigma}$ to be the permutation of cycle decompositions induced by (left) conjugation by $\pi$. More precisely, let $\bar \sigma = \pi \sigma \pi^{-1}$. If $\sigma$ has cycle decomposition $\tau_1 \cdots \tau_{N_\sigma}$ (we assume we have already fixed the order of the cycles for each $\sigma \in S$) and  $\bar \sigma$ has cycle decomposition $\bar \tau_1 \dots \bar \tau_{N_\sigma}$ , then for each $n=1, \dots, N_\sigma$, we have $\pi \tau_n \pi^{-1} = \bar \tau_{P_{\sigma,\pi}(n)}$.

Equivalently, we could define the permutation $P_{\sigma,\pi}$ by the property
\begin{equation}
    P_{\sigma,\pi}(r_\sigma(i)) = r_{\bar \sigma}(\pi(i)) \label{eq:Yr=rpi}
\end{equation}
or by
\begin{equation}
    i \in \cO_{\tau_k} \Leftrightarrow \pi(i) \in \cO_{\bar \tau_{P_{\sigma,\pi}(k)}}. \label{eq:ik-piiYk}
\end{equation}

The next lemmas are technical in nature, but they allow us to define the $\star$ action, ensure it is well defined, and show that it is compatible with our setup. 

\begin{lemma} \label{lem:YBpi} Let $\pi \in S$ and $\bar \sigma = \pi \sigma \pi^{-1}$. Then we have $[P_{\sigma, \pi}^{-1}] B_{\bar \sigma} [\pi] =B_\sigma$.
	\end{lemma}
\begin{proof}
We have 
\[
([P_{\sigma, \pi}^{-1}] B_{\bar \sigma} [\pi])_{ij} = (B_{\bar \sigma})_{P_{\sigma,\pi}(i), \pi(j)} = \left\{\begin{aligned} 1 & \text{ if } \pi(j) \in \cO_{\bar \tau_{P_{\sigma,\pi}(i)}} \\ 0 & \text{ otherwise} \end{aligned}\right\} = \left\{\begin{aligned}1 & \text{ if } j \in \cO_{\tau_i} \\ 0 & \text{ otherwise} \end{aligned}\right\} = B_{ij}.
\]
\end{proof}

The next lemma shows that the bases $\{\mathbf f_{k,\sigma \lambda}\}$ of Section~\ref{sec:fixed-basis} are particularly nice with respect to the $G$ action---namely, if we act on any element of a basis by any element $g\in G$, then we obtain a scalar multiple of an element of a basis. 

\begin{lemma} \label{lem:act-f}
	Let $n \in F_{\sigma\lambda}$, and $\pi \in S$ and $\mu \in H$.   Then
	\[
	\pi\mu\cdot \mathbf f_{n,\sigma\lambda} = c\, \mathbf f_{ P_{\sigma, \pi}(n), \pi\mu \sigma \lambda (\pi\mu)^{-1}}
	\]
	for some $c \in \mathbb C$.
\end{lemma}
\begin{proof}
	
We first observe
\[
(\pi \mu)  \sigma \lambda (\pi \mu)^{-1}  = \pi \sigma \pi^{-1} \pi \sigma^{-1} \mu \sigma \lambda \mu^{-1} \pi^{-1} = \bar \sigma (\mu.\sigma  \lambda \mu^{-1}).\pi^{-1}.
\] 
Let $k = P_{\sigma, \pi}(n)$ (i.e., conjugating by $\pi$ transforms $\tau_n$ to $\bar \tau_k$). Recall $s_\sigma(n)$ gives the smallest index in the cycle $\tau_n$, so $\pi(s_\sigma(n))$ is in the same $\bar \sigma$ orbit as $s_{\bar \sigma}(k)$, so we may pick $\ell$ so that
\begin{equation}
\pi(s_\sigma(n)) = \bar\sigma^\ell (s_{\bar\sigma}(k)). \label{eq:pick-l}
\end{equation}
Then
\begin{align*}
	\mathbf f_{k,\pi\mu \sigma \lambda(\pi\mu)^{-1}} =& \sum_{j=0}^{m-1} \left(\prod_{i=0}^{j-1} ((\mu .\sigma \lambda \mu^{-1}) .\pi^{-1})_{\bar\sigma^i(s_{\bar \sigma}(k))}\right) \mathbf e_{\bar\sigma^j( s_{\hat\sigma}(k))} \\
	=& \sum_{j=0}^{m-1} \left(\prod_{i=0}^{j-1} ((\mu .\sigma \lambda \mu^{-1}) .\pi^{-1})_{ \bar \sigma^{i-\ell} \pi (s_{\sigma}(n))}\right) \mathbf e_{ \bar \sigma^{j-\ell} \pi( s_{\sigma}(n))} \\
	=& \sum_{j=0}^{m-1} \left(\prod_{i=0}^{j-1} ((\mu .\sigma \lambda \mu^{-1}) .\pi^{-1})_{ \pi \sigma^{i-\ell} (s_{\sigma}(n))}\right) \mathbf e_{ \pi \sigma^{j-\ell} (s_{\sigma}(n))} \\
	=& \pi \cdot \sum_{j=0}^{m-1} \left(\prod_{i=0}^{j-1} (\mu .\sigma \lambda \mu^{-1})_{\sigma^{i-\ell} (s_{\sigma}(n))}\right) \mathbf e_{\sigma^{j-\ell} (s_{\sigma}(n))} \\
	=& \pi\cdot \sum_{j=0}^{m-1} \left(\prod_{i=0}^{j-1} \mu_{\sigma^{i-\ell+1}(s_\sigma(n))}\mu^{-1}_{\sigma^{i-\ell}(s_\sigma(n))}   \lambda_{\sigma^{i-\ell} (s_{\sigma}(n))}\right) \mathbf e_{\sigma^{j-\ell} (s_{\sigma}(n))} \\
	=& \pi\cdot \sum_{j=0}^{m-1}  \mu_{\sigma^{j-\ell}(s_\sigma(n))}\mu^{-1}_{\sigma^{-\ell}(s_\sigma(n))}  \left(\prod_{i=0}^{j-1} \lambda_{\sigma^{i-\ell} (s_{\sigma}(n))}\right) \mathbf e_{\sigma^{j-\ell} (s_{\sigma}(n))} \\
	=& \pi\cdot \sum_{J=-\ell}^{m-\ell-1}  \mu_{\sigma^{J}(s_\sigma(n))}\mu^{-1}_{\sigma^{-\ell}(s_\sigma(n))}  \left( \prod_{I=-\ell}^{J-1} \lambda_{\sigma^{I} (s_{\sigma}(n))}\right) \mathbf e_{\sigma^{J} (s_{\sigma}(n))} \\
	=& \mu^{-1}_{\sigma^{-\ell}(s_\sigma(n))} \pi \mu\cdot \sum_{J=-\ell}^{m-\ell-1}   \left(\prod_{I=-\ell}^{J-1} \lambda_{\sigma^{I} (s_{\sigma}(n))}\right)
	\mathbf e_{\sigma^{J} (s_{\sigma}(n))}. 
\end{align*}

Next, we note that, of course, $\sigma^J = \sigma^{J+m}$, and also $\prod_{I = J }^{J + m -1} \lambda_{\sigma^{i-\ell}(s_\sigma(n))} = 1$ since $n \in F_{\sigma \lambda}$. It follows that we may replace the negative values of $J$ with $J + m$. So we obtain
\begin{align*}
	\mathbf f_{k,\pi\mu \sigma \lambda(\pi\mu)^{-1}} =& \mu^{-1}_{\sigma^{-\ell}(s_\sigma(n))} \pi \mu\cdot \sum_{J=0}^{m-1}   \left(\prod_{I=-\ell}^{J-1} \lambda_{\sigma^{I} (s_{\sigma}(n))}\right) \mathbf e_{\sigma^{J} (s_{\sigma}(n))} \\
	=& \mu^{-1}_{\sigma^{-\ell}(s_\sigma(n))} \left(\prod_{I=-\ell}^{-1} \lambda_{\sigma^{I} (s_{\sigma}(n))} \right) \pi \mu \cdot  \sum_{J=0}^{m-1}   \left(\prod_{I=0}^{J-1} \lambda_{\sigma^{I} (s_{\sigma}(n))}\right) \mathbf e_{\sigma^{J} (s_{\sigma}(n))} \\
	=& c^{-1}\, \pi\mu \cdot \mathbf f_{n,\sigma\lambda}.
	\end{align*}
\end{proof}

Translating Lemma \ref{lem:act-f} to the dual, and recording the constant $c$ from the proof, we have the following.

\begin{corollary} \label{cor:act-z} For $n\in F_{\sigma\lambda}$, let $\ell$ satisfy $\pi(s_\sigma(n)) = \bar\sigma^\ell(s_{\bar \sigma}(P_{\sigma, \pi}(n)))$. Then
\[
    z_{n,\sigma \lambda}.(\pi\mu)^{-1} = \left( \mu^{-1}_{\sigma^{-\ell}(s_\sigma(n))} \prod_{i=1}^{\ell} \lambda_{\sigma^{-i}(s_\sigma(n))} \right)z_{P_{\sigma,\pi}(n), (\pi\mu)\sigma \lambda(\pi\mu)^{-1}},
\]
and in particular
\[
    z_{n, \sigma \lambda} .\mu = \mu_{s_\sigma(n)} z_{n, \mu^{-1} \sigma \lambda \mu}.
\]
\end{corollary}

We define maps $(-) \star (\pi\mu)^{-1}:\mathcal Q_{W_{\sigma\lambda }} \rightarrow \mathcal Q_{W_{(\pi\mu)\lambda\sigma(\pi\mu)^{-1}}}$ that ignore the scalars
\begin{equation}
 z_{n, \sigma\lambda} \star (\pi\mu)^{-1} := z_{P_{\sigma,\pi\mu}(n),(\pi\mu)\lambda\sigma(\pi\mu)^{-1}}.\label{eq:Yact} 
\end{equation}
These maps are compatible with $\psi$, that is,
\[
  \psi_{\pi \sigma \lambda \pi^{-1}}^*(z_{n,\sigma\lambda} \star \pi^{-1}) = \psi^*_{\sigma\lambda}(z_{n,\sigma\lambda}) \star \pi^{-1}.
\]
These maps also interact nicely with $W^{\sigma}$, as illustrated in the next lemma.
\begin{lemma}\label{lem:permutations}
    For any $\pi \in S$, we have $ W^{\sigma} \star \pi^{-1} = W^{\bar\sigma}$. In particular, $C_S(\sigma)$ is a group of permutation symmetries of $W^\sigma$ (via $P_{\sigma, \pi}$).
    \label{lem:W-sigpi}
\end{lemma}
\begin{proof} Write for the permutation $P:= P_{\sigma,\pi}$.
For each $k \in \{1, \dots, N_\sigma\}$, pick $\ell$ as in \eqref{eq:pick-l} so  that  
\[
\pi(s_{\sigma}(k)) =  \sigma^{\ell} s_{\bar\sigma}(P(k)).
\]

Recall
\[
    W^\sigma = \sum_{k=1}^{N_\sigma} M_{k} (x_i \mapsto z_{r_\sigma(i),\sigma}).
\]
So, 
\begin{align*}
    W^\sigma \star \pi^{-1} &= \sum_{k=1}^{N_\sigma} M_{s_\sigma(k)}(x_i \mapsto z_{P(r_{\sigma}(i)),\bar \sigma}) \\
    &= \sum_{k=1}^{N_\sigma} M_{s_\sigma(k)}(x_i \mapsto z_{r_{\bar\sigma}(\pi(i)),\bar \sigma}) \\
    &= \sum_{k=1}^{N_\sigma} (M_{s_\sigma(k)}.\pi^{-1})(x_i \mapsto z_{r_{\bar\sigma}(i),\bar \sigma}) \\
    &= \sum_{k=1}^{N_\sigma} M_{\pi(s_\sigma(k))}(x_i \mapsto z_{r_{\bar\sigma}(i),\bar \sigma})\\
    &= \sum_{k=1}^{N_\sigma} M_{\sigma^{\ell} s_{\bar\sigma}(P(k))} (x_i \mapsto z_{r_{\bar\sigma}(i),\bar\sigma})\\
    &= \sum_{k=1}^{N_\sigma} M_{s_{\bar\sigma}(P(k))} (x_i \mapsto z_{r_{\bar\sigma}{(\sigma^{\ell}(i))},\bar\sigma})\\
   &= \sum_{k=1}^{N_\sigma} M_{s_{\bar\sigma}(P(k))} (x_i \mapsto z_{r_{\bar\sigma}(i),\bar\sigma})\\
   &= \sum_{k=1}^{N_\sigma} M_{s_{\bar\sigma}(k)} (x_i \mapsto z_{r_{\bar\sigma}(i),\bar\sigma})\\
   &= W^{\bar\sigma}.
\end{align*}
\end{proof}

Of course, if $W'$ is obtained from $W$ by a renaming of variables, then the bases of Theorem~\ref{thm:kreu-basis} for $W$ and $W'$ will also correspond under the substitution. It follows then that the maps in Equation~\eqref{eq:Yact} can be assembled into a $G$-$\star$ action on $\Ba$: 
\[
    \lfloor p; \sigma \lambda \rceil \star (\pi \mu)^{-1} = \lfloor p \star (\pi\mu)^{-1};(\pi\mu)\sigma\lambda(\pi\mu)^{-1} \rceil.
\]
The $\star$ action on $\Ba$ differs from the dot action of Lemma~\ref{lem:ss-act} on $\mathcal H_{W,G}^{un} = \Span(\Ba)$ by a scalar.

If $\alpha = \beta_\sigma(\lambda)$ for some $\lambda \in H$, then we have $\beta_{\bar \sigma}(\lambda.\pi^{-1})$ corresponding to $B_{\bar \sigma} [\pi] [\lambda]$, which by Lemma~\ref{lem:YBpi} is equal to $[P_{\sigma,\pi}] B_\sigma [\lambda]$. So we have an $S$-$\star$ action on $\Ca$:
\begin{equation}
    \lfloor m;\alpha \rceil \star \pi^{-1} = \lfloor m \star \pi^{-1}; [P_{\sigma,\pi}][\alpha] \rceil. \label{eq:CstarS}
\end{equation}
These actions are compatible---that is, if $b \in B_c$, then $b \star \pi \in B_{c \star \pi}$.

For use in Section~\ref{sec:basis-for-inv}, we observe the following.
\begin{lemma} \label{lem:H-acts-Bc}
	The $\star$ action restricted to $H$ preserves $\Ba_c$.
\end{lemma}
\begin{proof}
	Looking at \eqref{eq:Yact}, we see that the $H$-$\star$ action is trivial on monomials.
	
	On group elements, we have $(\sigma \lambda).\mu = \mu^{-1} \sigma \lambda \mu = \sigma ((\mu^{-1}.\sigma) \mu \lambda)$. It is easy to see that $(\mu^{-1}.\sigma) \mu$ is in the kernel of $\beta_\sigma$, that is, $\beta_\sigma((\mu^{-1}.\sigma)\mu\lambda) = \beta_\sigma(\lambda) $.
\end{proof}

Of course, in the same way, we have a $\star$ action of $G^T$ on $\Ba^T$ and of $S$ on $\Ca^T$. 

\subsection{A basis for the invariants} \label{sec:basis-for-inv}
We now have two potential issues that we should resolve. First, we have a nice basis $\Ba$ for the \emph{unprojected} state space, but we need to get rid of the elements of $\Ba$ that are not involved in a basis for the state space, i.e., not every element of $\Ba$ will correspond to a nonzero sum as in Proposition~\ref{prop:statespaceiso}.

The second issue is that we have two different actions of $G$: The natural action of $G$ on the unprojected state space, which we denote by ``.'', and the $\star$ action of $G$ on $\Ba$. The $\star$ action is more convenient when working with a basis; we don't have to keep track of scalars, so we get an honest action on our preferred basis. 

In this section, we will resolve both of these issues, showing that the orbits of the $\star$ action correspond to a basis for the ($G$ invariant) state space. 

We first need to get rid of some of the elements of $\Ba$ that are not involved in a basis for the state space. Let $\Ea$ be the set of elements $\lfloor p;g  \rceil \in \Ba$ so that for some $h \in G$ and $a\in \CC$, we have $\lfloor p;g  \rceil.h = a \cdot \lfloor p;g  \rceil$, with $a \neq 1$. Set $\Ra=\Ba\setminus \Ea$. Of course, $\Ra$ is preserved by the $\star$ action. We also write $\Ra_c := \Ra \cap \Ba_c$ and $\Ea_c := \Ea \cap \Ba_c$. Furthermore, for $b \in \Ba$, let $\mathcal O_G(b)$ be its $\star$ orbit in $\Ba$, and let $\Ba/G$ denote the set of orbits.

\begin{lemma} \label{lem:orb-inv}
	The $G$-$\star$ orbits in $\Ra$ correspond to a basis for $\mathcal H_{W,G}$.
\end{lemma}

\begin{proof}
Recall that $\mathcal H_{W,G} = \Span(\Ba)^G$. Define a  linear map $p: \Span(\Ba) \rightarrow \mathcal H_{W,G}$ by $x \mapsto \frac{1}{|G|} \sum_{g \in G} x.g$. It is easy to check that $p(x.g) = p(x)$ and if $x$ is $G$ invariant, then $p(x) = x$. 

Next, define a set map 
\begin{align*}
P:&\Ba/G\to \HH_{W,G}\\
&\mathcal O_G(b) \mapsto p(b).
\end{align*}
This is well defined up to scalars. The images of different orbits have disjoint support (with respect to basis $\Ba$), so they are linearly independent if none of them are zero.
	
For $b \in \Ba$, we claim that $p(b) = 0$ if and only if $b \in \Ea$. First, assume $p(b)=0$. Let $a_g$ be the coefficient of $b$ in $b.g$. So of course $\sum_{g \in G} a_g = 0$. But $a_{id} = 1$, so there must be a $g$ with $a_g \neq 0$, and then $b.g = a_g b$. Conversely, assume that $b \notin \Ra$, that is, $b.g = a \cdot b$ for some $a \neq 1$. We know $g$ has finite order, say $k$. It follows that $a$ is a $k$-th root of unity, and since $a \neq 1$, we have $\sum_{i=0}^{k-1} a^i = 0$. Hence $\sum_{i=0}^{k-1} b.g^i = 0$. We also know that $p(b) = p(b.g)$. Hence $0 = p(\sum_{i=0}^{k-1} b.g^i) = k p(b)$, so $p(b) = 0$.
	
So the image $P(\Ra/G)$ is a linearly independent set.
	
Any $y \in \mathcal H_{W,G}$ can be written $y = \sum_i c_i r_i + \sum_j d_j s_j$, where $r_i \in \Ra$ and $s_j \in \Ea$. Since $y$ is $G$ invariant, we see that 
\[
y = p(y) = \sum_i c_i p(r_i) + \sum_j d_j p(s_j) = \sum_i c_i P(\mathcal O_G(r_i)) 
\]
for some scalars $c_i$ and $d_j$. Hence $\mathcal H_{W,G} = \Span(P(\Ra/G))$.
\end{proof}

Lemma \ref{lem:orb-inv} shows us that we need to understand the $\star$ orbits of the $G$ action.  Indeed, in order to give an isomorphism between the A- and B-model state spaces, we will only need to give a bijection between the $G$-$\star$ orbits in $\Ra$ and the $G^\vee$-$\star$ orbits in $\Ra^T$. We observe that 
\begin{equation}\label{eq:orbit-union}
\cO_G(x) = \coprod_{\sigma \in S} \cO_H(x \star \sigma) .
\end{equation}
Indeed, $x \star \sigma \lambda = (x \star \sigma) \star \lambda$. 
Similarly,  $\cO_{G^\vee}(y) = \coprod_{\sigma \in S} \cO_{H^T}(y)\star\sigma$.  Hence, we can understand the $G$-$\star$ orbits by first understanding the $H$-$\star$ orbits and then studying the $S$ action (see Lemma \ref{lem:orbit-action}) on the $H$-$\star$ orbits. 

So far, we have an isomorphism $c \leftrightarrow c'$ between $\Ca_\sigma$ and $\Ca_\sigma^T$, which leads us to examine the relationship between  $\Ba_c$ and $\Ba^T_{c'}$.  By Lemma~\ref{lem:H-acts-Bc}, the sets $\Ba_c$ are unions of $H$-$\star$ orbits.  We will first tackle the question: How many $H$-$\star$ orbits are there in a nonempty $\Ba_c$? As noted above, the action of $H$ doesn't affect the monomials, so we need only worry about the action (by conjugation) on the group elements indexing the sectors. Notice that the set of group elements of sectors in $\Ba_c$ is a coset of $\ker \beta_\sigma \cap H$ in $H$. It is easy to check that $[H,\sigma] \subseteq\ker \beta_\sigma \cap H$. By Corollary \ref{cor:coset-orbits}, we see that the number of $H$ orbits in $\Ba_c$ is
\begin{equation}
    \frac{|H \cap \ker \beta_\sigma|}{[H,\sigma]}. \label{eq:num-orbits}
\end{equation}

%We use Burnside's lemma: the number of orbits is equal to the average number of fixed points of an element of $H$. Notice that $\mu \in H$ fixes $\lambda \sigma$ if and only if $\mu$ and $\sigma$ commute.  Hence the average number of fixed points (and the number of orbits) is \[|C_{H}(\sigma)| \cdot |H\cap \ker\beta| / |H|\]. 

The next lemma will show that the number of $H$-$\star$ orbits in $\Ba_c$ is the same as the number of $H^T$-$\star$ orbits in $\Ba^T_{c'}$ (if both are nonempty). 

\begin{lemma} \label{lem:count-fibers}
	Let $\sigma \in S$ and $\beta=\beta_\sigma$ be the map defined in Equation~\eqref{def-beta}. Then 
	
	\begin{equation}
\frac{|H \cap \ker\beta|}{|[H,\sigma]|} = \frac{|H^T \cap \ker\beta|}{|[H^T,\sigma]|} \label{eq:reduced-to}
\end{equation} 
\end{lemma}

\begin{proof}
 Let $J_i$ be the symmetry of $W$ corresponding to the $i$-th column of $A_W^{-1}B^T$. Then we can describe $\ker\beta$ as
\begin{align*}
    \ker\beta &= \{ \lambda \in \diag({\CC^N}) : B[\lambda] \in \ZZ^k\} \\
    		  &=  \{ \lambda \in \diag({\CC^N}) : [\lambda]^T A_W A_W^{-1} B^T \in \ZZ^k\} \\
              &=  \{ \lambda \in \diag({\CC^N}) :  [\lambda]^T A_W [J_i] \in \ZZ \text{ for } i =1, \dots, k \} \\
              &= \left< J_1, \dots, J_k\right>^T.
\end{align*}

Now, by the index-preserving property of the transpose (see Proposition~\ref{lem:dual-props}), we have
\begin{equation}
\frac{|H^T \cap \ker\beta|}{|[H^T,\sigma]|} = \frac{[H^T,\sigma]^T}{|(H^T \cap \ker\beta)^T|}. \label{eq:by-index-preserving}
\end{equation}

In general, it is true that $(H_1^T \cap H_2^T) = \left<H_1, H_2\right>^T$, so we see that 
\begin{equation}
(H^T \cap \ker\beta)^T = \left< H, J_1, \dots, J_k \right>. \label{eq:transpose-intersection}
\end{equation}

We have (setting $A=A_W$):
\begin{align*}
 [H^T,\sigma]^T &= \{ \mu \in G_W^{diag} : ([\lambda]-[\sigma][\lambda])^TA[\mu] \in \ZZ, \forall \lambda \in H^T \}    \\
                &= \{ \mu \in G_W^{diag} : ([\lambda]^TA-[\lambda]^T[\sigma]^TA)[\mu] \in \ZZ, \forall \lambda \in H^T \} \\
                &= \{ \mu \in G_W^{diag} : ([\lambda]^TA-[\lambda]^TA[\sigma]^T)[\mu] \in \ZZ, \forall \lambda \in H^T \} \\ 
                &= \{ \mu \in G_W^{diag} : [\lambda]^TA([\mu]-[\sigma]^T[\mu]) \in \ZZ, \forall \lambda \in H^T \} \\
                &= \{ \mu \in G_W^{diag} : [\mu]-[\sigma]^T[\mu] \in H\} \\
                &= \{ \mu \in G_W^{diag} : [\mu]-[\sigma][\mu] \in H\}.
\end{align*}

In other words, if we have $\phi: G_W^{diag} \rightarrow \ker\beta$ defined by $[\mu] \mapsto [\mu]-[\sigma][\mu]$, then $[H^T,\sigma]^T = \phi^{-1} (H \cap \ker\beta)$. Lemma \ref{lem:psi-sur-ker} below shows that $\phi$ is surjective.  It follows (from, e.g., the third isomorphism theorem for groups) that 
\begin{equation}
\frac{|G_W^{diag}|}{|[H^T,\sigma]^T|} = \frac{|\ker\beta|}{|H \cap \ker\beta|} \label{eq:by-third}.
\end{equation}
Using \eqref{eq:by-index-preserving},\eqref{eq:transpose-intersection}, and \eqref{eq:by-third}, we obtain
\[
 \frac{|H^T \cap \ker\beta'|}{|[H^T,\sigma]|} = \frac{|G_W^{diag}||H \cap \ker\beta|}{|\ker\beta|\; |\left<H, J_1, \dots, J_k\right>|}.
\]
Comparing to \eqref{eq:reduced-to}, we see that we need to show
\begin{equation*}
    \frac{|H \cap \ker\beta|}{|[H,\sigma]|} = \frac{|G_W^{diag}||H \cap \ker\beta|}{|\ker\beta|\; |\left<H, J_1, \dots, J_k\right>|},
    \end{equation*}
or equivalently
\begin{equation}
    \frac{|\left<H, J_1, \dots, J_k\right>|}{|[H,\sigma]|} = \frac{|G_W^{diag}|}{|\ker\beta| }.
\end{equation}
The right-hand side is, of course, equal to $|\ker \phi|$. Notice that $\ker \phi$ is precisely the $\lambda \in G_W^{diag}$ so that $\lambda_i = \lambda_j$ whenever $\cO(i) = \cO(j)$. Also, Lemma \ref{lem:psi-sur-ker} tells us that $\ker \phi = \left< J_1, \dots, J_k\right>$. 

If one restricts the domain of $\phi$ to $\left<H, J_1, \dots, J_k\right>$, one sees that the image is $[H,\sigma]$. It follows that the quotient on the left is equal to the size of the kernel of the restricted $\phi$. But the whole kernel of $\phi$ is contained in the restricted domain, so the kernel of $\phi$ restricted is the same as the kernel of $\phi$ (unrestricted), whose size is given by the right-hand side.
 \end{proof}

\begin{lemma} \label{lem:psi-sur-ker}
	Consider the map $\phi: G_{W}^{diag} \rightarrow \ker\beta$ given by $[\mu] \mapsto [\mu] - [\sigma][\mu]$. Then the following hold:
	\begin{enumerate}
	    \item $\phi$ is surjective.
	    \item The kernel of $\phi$ is $\left<J_1, J_2, \dots, J_k \right>$ (where $J_i$ are represented by the columns of $A^{-1}B^T$).
	\end{enumerate}
\end{lemma}
\begin{proof}
    (1) Assume $\lambda \in \ker \beta$. Let $y$ be a representative of $[\lambda]$ in $\QQ^N$ so that $B y = 0$. It is easy to see that the nullspace of $(I-[\sigma]^T)$ (in $\QQ^N$) is precisely the row space of $B$, and hence the range of $(I - [\sigma])$ is the kernel of $B$. It follows that there is a solution $x \in \QQ^N$ for $(I-[\sigma])x = y$. Let $m$ be the order of $\sigma$, and set $x' = x - \frac1m \sum_{i=1}^m [\sigma]^i x$. One can check that $(I-[\sigma])x' = (I-[\sigma])x = y$.
    
    We claim that $Ax' \in \ZZ^N$, so that $x'$ represents an element of $G_W^{diag}$. To see this, we first note that since $x = [\sigma] x + y$ and $Ay \in \ZZ^N$, we have that $Ax \equiv A[\sigma]^\ell x \mod \ZZ^N$ for any $\ell$. It follows that 
    \[
    Ax' = Ax - \frac1m \sum_{i=1}^m A[\sigma^i]x \equiv Ax -\frac1m \sum_{i=1}^m Ax = 0 \mod \ZZ^N. 
    \]
    
    (2) Now, assume $\lambda \in \ker \phi$. Then $\lambda$ is represented by $y \in \QQ^N$, and we may assume that $[\sigma]y = y$. Let $z = Ay$, which is in $\ZZ^N$. It follows that $[\sigma]z = z$ as well. Hence, $z = B^Tx$ for some $x \in \QQ^k$. But since $B^T$ has only one nonzero entry per row, and it is a one, we see that $x\in \ZZ^k$. So $y=A^{-1}Bx$, as desired.
\end{proof}

\begin{corollary} \label{cor:count-non-empty}
	If $c \in \Ca_\sigma$ and $c' \in \Ca^T_\sigma$ correspond under the Krawitz isomorphism \eqref{eq:kmm} and both $\Ba_c$ and $\Ba^T_{c'}$ are nonempty, then the number of $H$-$\star$ orbits in $\Ba_c$ is the same as the number of $H^T$-$\star$ orbits in $\Ba_{c'}$.
	
\end{corollary}

From this corollary, we see that we get the same number of $H$-$\star$ orbits as $H^T$-$\star$ orbits, as long as $\Ba_c$ is nonempty. So the next question is: When is \texorpdfstring{$\Ba_c$}{Bc} empty?  In the next section we define a condition that  must be satisfied in order to ensure we can define a mirror map between the $\star$ orbits on the A-model and the B-model.

\subsection{The Diagonal Scaling Condition} \label{sec:strongPC}

As mentioned in the introduction, we do not expect to get a mirror isomorphism unless the group $G$ satisfies certain conditions. In \cite{EG2}, Ebeling and Gusein-Zade defined a condition, called the \emph{Parity Condition} (PC), with the idea that PC was necessary in order to define a mirror map. We have found, on the other hand, that there are certain other conditions that should be met, which we will now describe. It is not clear what connection, if any, there is between our conditions and PC.

The following lemma gives us necessary and sufficient conditions to tell when $\Ba_c$ is nonempty in terms of the action of $H^T$ on $\Ba^T_c$.

\begin{lemma} \label{lem:empty-fibers}
Suppose $c\in \Ca_\sigma$ and $c'\in \Ca^T_\sigma$ correspond under the Krawitz isomorphism \eqref{eq:kmm} and $\Ba_c$ is nonempty. Then the following are equivalent:
    \begin{enumerate}
        \item $\Ba_{c'}^T$ is empty.
        \item For all $b \in \Ba_c$,  there exists $\mu \in C_H(\sigma)$ and $a \in \mathbb C$, $a \neq 1$ with $b.\mu = a \cdot b$. 
        \item For some $b \in \Ba_c$, there exists $\mu \in C_H(\sigma)$ and $a \in \mathbb C$, $a \neq 1$ with $b.\mu = a \cdot b$.
    \end{enumerate}
\end{lemma}

In other words, $\Ba_{c'}^{T}$ is empty if and only if there is some element in the fiber $\Ba_c$ that is scaled by a diagonal symmetry. Note the subtle difference between this condition and the condition that $b\in \Ea$. Furthermore, conditions (2) and (3) say that if one element $b\in \Ba_c$ has the given property, then they all do.

\begin{proof}An element $b \in \Ba_c$ is of the form $b = \lfloor (j^{-1})^* m ; \sigma \lambda \rceil$, where $B_\sigma [\lambda] = [\alpha]$. Then, for $\mu \in C_H(\sigma)$, we have $b.\mu = \lfloor ((j^{-1})^* m).\mu ; \mu^{-1}\sigma \lambda \mu \rceil = \lfloor  (j^{-1})^*(m.\gamma(\mu)) ; \sigma \lambda \rceil$. So we see that $b.\mu = a \cdot b$ if and only if $ m.\gamma(\mu) = a \cdot m$. This is independent of $\lambda$, so we see that (2) and (3) are equivalent.

 Let $c=\lfloor m; \alpha \rceil$ and $c':=\lfloor p;\epsilon \rceil$, and assume (3), so we have $\mu \in C_H(\sigma)$ with $ m.\gamma(\mu) = a m$, with $a \neq 1$. By \eqref{eq:kmm}---esp. the discussion following---this  implies that $\epsilon$ is not in $\left< \gamma(\mu) \right>^T$, that is, $[\epsilon]^T A_{W^\sigma} C [\mu] \notin \ZZ$. Assume, by way of contradiction, $\lfloor (j^{-1})^*(p); \sigma \nu\rceil \in \Ba_{c'}^T$, so $B[\nu] = [\epsilon]$ and $\nu \in H^T$. Since $\mu \in C_H(\sigma)$, there is an $\ell \in \QQ^{N_\sigma}$ with $B^T\ell = [\mu]$. So we have $[\nu]^TB^T A_{W^{\sigma}}CB^T\ell \notin \ZZ$. By Lemma~\ref{lem:matrix-prop}, 
	\[
    	B^T A_{W^{\sigma}}CB^T\ell = B^T A_{W^\sigma}\ell = A_W B^T \ell = A_W [\mu],
	\] so $[\nu]^T A_W [\mu] \notin \ZZ$, so $\nu \notin H^T$, a contradiction. We have shown (1). 

    On the other hand, assume (1), that is $\Ba_{c'}^T$ is empty. This happens exactly when $\epsilon \notin \beta(H^T)$. It follows that $m$ is not invariant under some $\mu' \in \beta(H^T)^T$. Saying $\mu' \in \beta(H^T)^T$ is saying that
    \begin{align*}
        [\nu]^T B^T A_{W^{\sigma}} [\mu'] \in \mathbb Z,\; \forall \nu \in H^T\\
        [\nu]^T A_{W} B^T[\mu'] \in \mathbb Z,\; \forall \nu \in H^T \\
        B^T[\mu'] \in H.
    \end{align*} 
    Let $\mu = B^T \mu'$. So $\mu \in C_H(\sigma)$. And since $CB^T = I$, we see $\gamma(\mu) = \mu'$. It follows that $(j^{-1})^* m$ is not invariant under $\mu$, so we obtain (2).
\end{proof}

With this lemma in mind, in order to have a mirror map, we will need the following condition.

\begin{condition} \label{PC}
We say the pair $(W,G)$ satisfies the \emph{diagonal scaling condition} (DSC) if the following are true:
\begin{enumerate}
    \item For each $c$ where $\Ea_c$ is nonempty (that is, there exists $b \in \Ba_c$ and $g$ so that $b.g = a b$ with $a \neq 1$), the (equivalent) conditions of Lemma \ref{lem:empty-fibers} are satisfied.
    \item For each $c'$ where $\Ea_{c'}^T$ is nonempty,  the (equivalent) conditions of Lemma \ref{lem:empty-fibers} are satisfied (with $c$ and $c'$ interchanged and transposes applied).
\end{enumerate}
\end{condition}

In other words, this condition says that a pair satisfies DSC if whenever there is a ``bad'' basis element in $\Ba_c$, then it is because of an element of $H$ (or rather $C_H(\sigma)$). And then every element of $\Ba_c$ is ``bad'' and $\Ba_{c'}^T$ is empty. When this condition holds, we obtain the following corollary. 

\begin{corollary} \label{cor:count-fibers-good}
	If the DSC \ref{PC} is satisfied and $c \in \Ca$ and $c' \in \Ca^T$ correspond under the Krawitz isomorphism, then the number of $H$-$\star$ orbits in $\Ra_c$ is the same as the number of $H^T$-$\star$ orbits in $\Ra^T_{c'}$.
\end{corollary}
\begin{proof}
    If both $\Ba_{c}$ and $\Ba_{c'}^T$ are empty, we are done.
    
    If $\Ba_{c'}^T$ is empty and $\Ba_c$ is nonempty, then by Lemma \ref{lem:empty-fibers} we see that $\Ra_c$ is empty, as every element of $\Ba_c$ belongs to $\Ea_c$. The case of $\Ba_{c}$ empty and $\Ba_{c'}^T$ nonempty is similar.
    
    So now we assume that both $\Ba_c$ and $\Ba_{c'}^T$ are nonempty.

    If $\Ea_c$ is nonempty, then the DSC says that  $\Ba_{c'}^T$ is empty and $\Ra_c$ is empty.
    
    Similarly, if $\Ea^T_{c'}$ is nonempty, then both $\Ba_{c}$ and $\Ra^T_{c'}$ are empty.
    
    If both $\Ea_c$ and $\Ea^T_{c'}$ are empty, then $\Ra_c = \Ba_c$ and $\Ra^T_{c'} = \Ba_{c'}^T$, which have matching $\star$ orbit sizes by Corollary \ref{cor:count-non-empty}.
\end{proof}

Thus for each $c \in \Ca_{\sigma}$, corresponding to a $c' \in \Ca_{\sigma}^T$, we can pick a bijection $\Phi_c$ between $H$-$\star$ orbits in $\Ra_c$ and $H^T$-$\star$ orbits in $\Ra^T_{c'}$, and assemble these into a bijection $\Phi$ between $H$-$\star$ orbits in $\Ra$ and $H^T$-$\star$ orbits in $\Ra^T$.

\subsection{Equivariance of the bijection of \texorpdfstring{$H$}{H} and \texorpdfstring{$H^T$}{HT} orbits} \label{sec:EV}
So far, we have described a basis $\Ba$ for the unprojected state space and chosen a subset $\Ra$ of that basis with the property that the orbits of the $\star$ action of $G$ on $\Ra$ correspond to a basis of the state space. So our goal is to find a bijection between the $G$-$\star$ orbits of $\Ra$ and the $G^\vee$-$\star$ orbits of $\Ra^T$. We found that we have a bijection between the $H$-$\star$ orbits and the $H^T$-$\star$ orbits, under certain conditions. Finally, we need to consider the action of $S$ on these orbits. 

In order to obtain an isomorphism, one more hypothesis is required, which we discuss here.
\begin{lemma} \label{lem:orbit-action}
	Define $\cO_H(b)\star\sigma := \cO_H(b \star \sigma)$. This is a well-defined action of $S$ on the $H$-$\star$ orbits in $\Ra$. There is a similar action of $S$ on the $H^T$-$\star$ orbits in $\Ra^T$. 
\end{lemma}
\begin{proof}
    If $b' = b \star \lambda$ is a different representative of the orbit, then $b' \star \sigma = b \star \lambda \sigma = b \star \sigma (\lambda.\sigma) = (b \star \sigma) (\lambda.\sigma)$, so $b'\star \sigma$ represents the same $H$ orbit as $b \star \sigma$.   
\end{proof}

\begin{condition}[Equivariant $\Phi$] \label{EV}
    We say that the pair $(W,G)$ satisfies the Equivariant $\Phi$ condition if
    \begin{enumerate}
        \item $(W,G)$ satisfies the DSC \ref{PC}.
        \item The bijection $\Phi$ between $H$ orbits and $H^T$ orbits from the previous section can be chosen so that it is $S$ equivariant for the $\star$ action of Lemma \ref{lem:orbit-action}.
    \end{enumerate}
\end{condition}
The examples in Sections \ref{sec:good-act-H-orbits} and \ref{sec:bad-act-H-orbits} may help illustrate why this condition is necessary. It is often the case that there is only one $G$-$\star$ orbit, in which case the there is a canonical choice for an isomorphism of state spaces. However, if there is more than one orbit, then the choice of isomorphism might not be canonical---at least if we are only considering the structure of graded vector spaces. It may be that considerations of higher Frobenius structures may force a canonical choice. 

\begin{lemma} \label{lem:hatPhi}
	Assume $(W,G)$ satisfies the Equivariant $\Phi$ condition (\ref{EV}).
	Then (an equivariant) $\Phi$ induces a bijection $\hat \Phi$ between $G$-$\star$ orbits in $\Ra$ and $G^T$-$\star$ orbits in $\Ra^T$.  
\end{lemma}
\begin{proof}
	Recall from \eqref{eq:orbit-union} that a $G$-$\star$ orbit is a union of $H$-$\star$ orbits. We set  and check
	\begin{align*}
		\hat \Phi( \cO_G(x) ) &:= \coprod_{\sigma \in S} \Phi(\cO_H(x)\star\sigma) \\
							&=\coprod_{\sigma \in S} \Phi(\cO_H(x))\star\sigma \\
							&= \coprod_{\sigma \in S} \cO_{H^T}(y)\star\sigma \\
							&= \cO_{G^T}(y)
	\end{align*}
for some $y \in \Ra^T$. It is easy to check that $\cO_{G^T}(y) \mapsto  \coprod_{\sigma \in S} \Phi^{-1}(\cO_{H^T}(y) \star \sigma)$ gives the inverse.
\end{proof}

Thus we arrive at our main theorem:

\begin{theorem} \label{thm:main}
	Assume $(W,G)$ satisfies the Equivariant $\Phi$  condition \ref{EV}. Then there is a bidegree-preserving isomorphism of state spaces $\mathcal A_{W,G} \cong \mathcal B_{W^T,G^T}$.
\end{theorem}
\begin{proof}
	This follows from Lemma \ref{lem:orb-inv}, Lemma \ref{lem:hatPhi}, and Corollary \ref{cor:bidegrees}.
\end{proof}

\section{Cyclic groups of prime order} \label{sec:prime-order}

The Equivariant $\Phi$ Condition \ref{EV} (including DSC \ref{PC}) is sufficient for establishing a mirror isomorphism, but it is quite technical and not easy to check. In this section, we will show that the Equivariant $\Phi$ Condition (and therefore also DSC) is automatically satisfied when $S$ is cyclic with order an odd prime. Then as a corollary to \label{thm:main} %this label is not showing up in pdf
we will obtain the following result. Interestingly, in this case, the PC of \cite{EG2} is also satisfied. 

\begin{theorem} \label{thm:main-cyc}
    If $S$ is cyclic with order an odd prime, then there is a bidegree-preserving isomorphism of state spaces $\mathcal A_{W,G} \cong \mathcal B_{W^T,G^T}$.
\end{theorem}

\subsection{DSC \ref{PC} for cyclic permutation groups}
First, we want to show that the DSC \ref{PC} is satisfied when $S$ is cyclic of odd prime order. Remember DSC says that if an element is ``bad,'' it is because of an element of $H$. These first two lemmas show that the action of certain symmetries are given by the action of a diagonal symmetry. 

\begin{lemma} \label{lem:omega}
 Suppose that $(\sigma^j \mu)^{-1}$, with $j \ge 0$, commutes with $\sigma \lambda$. Let 
 \[
 \omega = \left(\Big(\mu^{-1}.\sigma^{-j}\Big) \prod_{i=1}^j \Big(\lambda.\sigma^{-i}\Big) \right) \in H.
 \]
 Then
 \[
    z_{n, \sigma \lambda}. (\sigma^j \mu)^{-1} = z_{n, \sigma \lambda} .\omega^{-1}.
 \]
\end{lemma}
\begin{proof}
We apply Corollary~\ref{cor:act-z}.  Since $P_{\sigma, \sigma^j} = \text{id}$, we may take $\ell = j$. Then the left-hand side is 
\[
    \left(\mu^{-1}_{\sigma^{-j}(s(n))} \prod_{i=1}^j \lambda_{\sigma^{-i}(s(n))} \right) z_{n,\sigma \lambda}.
\]
The right-hand side is
\[
    \left(\mu^{-1}.\sigma^{-j} \prod_{i=1}^j \lambda.\sigma^{-i} \right)_{s(n)} z_{n, \omega \sigma \lambda \omega^{-1}}.
\]
The coefficients match, so it remains only to check that $\omega$ commutes with $\sigma \lambda$, so that $z_{n,\sigma\lambda}=z_{n,\omega\sigma\lambda\omega^{-1}}$. Indeed,
\begin{align*}
    \left( \left(\prod_{i=1}^j \lambda.\sigma^{-i}\right) (\mu^{-1}.\sigma^{-j})\right) \sigma \lambda &= \left(\prod_{i=1}^j \sigma^i \lambda \sigma^{-i}\right) \sigma^j (\mu^{-1} \sigma^{-j}) (\sigma \lambda) \\
    &= \left(\prod_{i=1}^j \sigma^i \lambda \sigma^{-i}\right) \sigma^j (\sigma \lambda) (\mu^{-1} \sigma^{-j})  \\
    &= \sigma \lambda \sigma^{-1} \left(\prod_{i=2}^{j-1} \sigma^i \lambda \sigma^{-i} \right) (\sigma^j \lambda \sigma^{-j}) \sigma^j (\sigma \lambda) (\sigma^{-j} \sigma^j) (\mu^{-1} \sigma^{-j}) \\
    &= (\sigma \lambda) \sigma^{-1} \left(\prod_{i=2}^{j-1} \sigma^i \lambda \sigma^{-i} \right) \sigma \sigma^{j-1} \lambda (\sigma \lambda) \sigma^{-j} (\mu^{-1}.\sigma^{-j}) \\
    &= (\sigma \lambda)  \left(\prod_{i=1}^{j-2} \sigma^i \lambda \sigma^{-i} \right) \sigma^{j-1} \lambda (\sigma^{1-j} \sigma^{j-1})(\sigma \lambda) \sigma^{-j} (\mu^{-1}.\sigma^{-j}) \\
    &= (\sigma \lambda)  \left(\prod_{i=1}^{j} \sigma^i \lambda \sigma^{-i} \right) (\mu^{-1}.\sigma^{-j}).
\end{align*}
\end{proof}

\begin{lemma} \label{lem:diag-sector-action}
 Suppose that $\pi\mu$ commutes with $\lambda$. Let $F_\lambda$ be the set of indexes $i$ so that $x_i$ is fixed by $\lambda$. If $m$ is a monomial in $\{x_i\}_{i \in F_\lambda}$ with $m.(\pi\mu)^{-1} = a\cdot m$ for some $a \in \CC$, then we have $m.(\pi\mu)^{-1} = m.\mu^{-1}$.
\end{lemma}
\begin{proof}
    We have $m.\mu^{-1} = a \cdot m.\pi$. Since $m.\mu^{-1}$ must be a scalar multiple of $m$, and $m.\pi$ cannot introduce any scalars, we have  $m.\pi = m$. So $m.\mu^{-1} = a\cdot m$ and $m.\pi^{-1} = m$. Then $m.(\pi\mu)^{-1} = m.\mu^{-1} \pi^{-1} = (a\cdot m).\pi^{-1} = a \cdot m = m.\mu^{-1}$. 
\end{proof}

\begin{corollary}
    If $S$ is cyclic of odd prime order, then the DSC \ref{PC} holds for any $H$.
\end{corollary}

\begin{proof}
We need to show that if $b\in \Ea$, then it is because of an element of $H$. Let $b = \lfloor m;\sigma\lambda \rceil \in \Ba$ and $\pi \mu \in G$ so that $b.\pi\mu = a \cdot b$ for some $a \in \CC$, $a \neq 1$. If this is the case, then $\tau\mu$ commutes with $\sigma\lambda$ and we know that $\pi\mu$ permutes the variables $\{z_{i,\sigma\lambda}\}_{i \in F_{\sigma \lambda}}$ according to $P_{\sigma,\pi}$ (and rescales them). The order of $P_{\sigma,\pi}$ must divide the order of $\pi$, so it is odd. It follows that $P_{\sigma,\pi}$ is an even permutation, so action on the volume form does not introduce a sign coming from the permutation of variables.

If $\sigma$ is not the identity, then $\pi = \sigma^j$ for some $j$. Lemma \ref{lem:omega} then shows that $\pi\mu$ acts on $m$ the same as an element of $H$, as desired.
    
If $\sigma$ is the identity, then $\pi\mu$ commutes with $\lambda$ and we may apply Lemma~\ref{lem:diag-sector-action} to see that $\pi\mu$ acts on $m$ the same as $\mu \in H$, as desired.
\end{proof}

\subsection{Reducing to the action on \texorpdfstring{$\Ba_c$}{Bc}} \label{sec:reducing}

Recall that if $b \in \Ba_c$, then $b \star \pi \in \Ba_{c\star \pi}$. Furthermore, Lemma~\ref{lem:kraw-comp} tells us that if $c$, $c'$ is a Krawitz pair in $\Ca_\sigma$ and $\Ca^T_\sigma$, then $c \star \pi$, $c' \star \pi$ is a Krawitz pair in $\Ca_{\bar \sigma}$ and $\Ca_{\bar \sigma}^T$. This compatibility  allows us to check Equivariant Condition \ref{EV} by only looking at the $S$ action on each $\Ba_c$. More precisely, assume the DSC \ref{PC} and assume $\Ba_c \cap \Ra$ and $\Ba^T_{c'}\cap \Ra^T$ are nonempty. Then, by Lemma \ref{cor:count-fibers-good},  there is a bijection between the set of $H$ orbits  in $\Ba_{c} \cap \Ra$ and the set of $H^T$ orbits in $\Ba^T_{c'}\cap \Ra^T$, call the bijection $\Phi_c$.

Let $\Stab_S(c) = \{\sigma \in S: c \star \sigma = c \}$. It follows from Lemma~\ref{lem:kraw-comp} that $\Stab_S(c) = \Stab_S(c')$. Call this group $S_c$. Then $S_c$ acts on the set of $H$ orbits  in $\Ba_c \cap \Ra$ and on the set of $H^T$ orbits in $\Ba^T_{c'} \cap \Ra^T$. 

\begin{lemma}[Equivariant fiber actions] \label{lem:efa}
	 Suppose $(W,G)$ satisfies the DSC \ref{PC}.
    Suppose that, for all Krawitz pairs $c$,$c'$ as above, the bijection $\Phi_c$ between the set of $H$ orbits  in $\Ba_c \cap \Ra$ and the set of $H^T$ orbits in $\Ba_{c'} \cap \Ra^T$ can be picked to be $S_c$ equivariant.
    
    Then we can satisfy the Equivariant Condition \ref{EV}, that is, we can pick a bijection $\Phi$ between $H$ orbits in $\Ra$ and $H^T$ orbits in $\Ra^T$ that is equvariant for the $S$-$\star$ action.
\end{lemma}
\begin{proof}
	For each $S$ orbit $\cO_S(c)$ in $\Ca$, fix a representative $t(c) \in \Ca$. 
	
	Now, for $b \in \Ba_c$, pick $\sigma$ so that $c \star \sigma = t(c)$. Then we define
	\[
	    \Phi(\cO_H(b)) := \Phi_{t(c)}(\cO_H(b) \star \sigma) \star \sigma^{-1}.
	\]
	This is independent of the choice of $\sigma$. Indeed, if there is also $\sigma'$ so that $c \star \sigma' = t(c)$, then $\sigma^{-1} \sigma'$ fixes $t(c)$. Thus, by hypothesis,   
	\[
	    \Phi_{t(c)}(\cO_H(b) \star \sigma) \star (\sigma^{-1} \sigma') = \Phi_{t(c)}(\cO_H(b) \star \sigma'),
	\]
	so we conclude that 
	\[
    \Phi_{t(c)}(\cO_H(b) \star \sigma) \star \sigma^{-1} = \Phi_{t(c)}(\cO_H(b) \star \sigma') \star (\sigma')^{-1},
	\]
	as desired.
	
	We also see that $\Phi$ is $S$ equivariant. Again, for $b \in \Ba_c$, pick $\sigma$ so that $c \star \sigma = t(c)$. For any $\tau \in S$,  we see that $(c \star \tau) \star (\tau^{-1} \sigma) = c \star \sigma  = t(c) =  t(c \star \tau)$. Hence
	\begin{align*}
	    \Phi(\cO_H(b) \star \tau) &= \Phi(\cO_H(b\star \tau)) \\
	    &= \Phi_{t(c \star \tau)}(\cO_H(b \star \tau) \star \tau^{-1} \sigma ) \star \sigma^{-1}  \tau \\
        &= \Phi_{t(c \star \tau)}(\cO_H(b)) \star \tau.
	\end{align*}
\end{proof}

\subsection{Trivial actions on \texorpdfstring{$\Ba_c$}{Bc}}
One special case when the hypothesis of Lemma \ref{lem:efa} is satisfied is when the action of $S_c$ is trivial.

Recall from \eqref{eq:num-orbits} that if $[H,\sigma] =\ker \beta_{\sigma} \cap H$ for all $\sigma \in S$, then each $\Ba_c$ has only one $H$-$\star$ orbit. This is the case, for example, when $H=G_W^{diag}$ (Lemma \ref{lem:psi-sur-ker}) or when $H$ is the trivial group.  

We now give conditions for the fiber actions to be trivial, even if $\Ba_c$ has more than one $H$ orbit.

\begin{lemma} Let $c \in \Ca_\sigma$, $b = \lfloor m, \sigma \lambda \rceil \in \Ba_c$, and $\tau \in S_c$. Then $\sigma$ and $\tau$ commute and $(\lambda.\tau) \lambda^{-1} \in \ker \beta_\sigma$. 

If we suppose further that $\lambda.\tau \lambda^{-1} \in [H,\sigma]$, then $\cO_H(b) \star \tau = \cO_H(b)$. \label{lem:triv-fib-act}
\end{lemma}
\begin{proof}
	Since $c \in \Ca_\sigma$ and $c \star \tau \in \Ca_{\tau^{-1}\sigma\tau}$, we see that $\sigma$ and $\tau$ commute. 
    We have
    \[
    b \star \tau =  \lfloor m \star \tau, (\tau^{-1}\sigma\tau)  (\lambda.\tau) \rceil. 
    \]
   	Hence $\beta_\sigma (\lambda) = \beta_\sigma (\lambda.\tau)$, which implies that $(\lambda.\tau)\lambda^{-1} \in \ker \beta_\sigma$. We also see that $m \star \tau = m$.
   	
    For any $\mu \in H$
    \[
    b \star \mu = \lfloor m ,\sigma (\mu^{-1}.\sigma) \lambda \mu \rceil.
    \]
    So, to see whether $b$ and $b \star \tau$ are in the same $H$ orbit, the question is whether we can find a $\mu \in H$ with
    \[
     \lambda.\tau = (\mu^{-1}.\sigma) \lambda \mu,
     \]
     which is equivalent to finding a $\mu$ such that
     \[ (\lambda.\tau) \lambda^{-1} = (\mu^{-1}.\sigma) \mu,
    \]
    which is saying that $(\lambda.\tau) \lambda^{-1}$ is an element of $[H,\sigma]$.
\end{proof}

\begin{corollary} \label{cor:trivial-fiber-action}
    Suppose that for all commuting pairs $\sigma, \tau \in S$, we have $[H,\tau] \cap \ker \beta_\sigma \subseteq [H,\sigma]$. Then, for each $c \in \Ca$, the action of $S_c$ on the $H$-$\star$ orbits in $\Ba_c$ is trivial.
\end{corollary}

\begin{proof}
    Let $b = \lfloor m ,\sigma \lambda \rceil$ and $\tau \in S_c$. Then by the lemma we have $\tau \in C_S(\sigma)$ and $(\lambda.\tau) \lambda^{-1} \in \ker \beta_\sigma$, so by hypothesis $\lambda.\tau \in [H,\sigma]$. Then apply the second part of the lemma.
\end{proof}

\begin{proposition}
The hypothesis of Corollary~\ref{cor:trivial-fiber-action} is satisfied when $S$ is a cyclic group of prime order. \label{prop:prime} 
\end{proposition}
\begin{proof}
Suppose $S =\left<\sigma\right>$ is cyclic and of prime order $p$. Let $j \in \mathbb Z$. We will show that $[H,\sigma] \subseteq [H, \sigma^j]$. A general element of $[H,\sigma]$ is $\lambda(\lambda^{-1}.\sigma)$ for some $\lambda \in H$. Let $k$ be a natural number so that $kj \equiv 1 \mod p$. We set 
\[
 \mu = \lambda \, (\lambda.\sigma^j) \, (\lambda.\sigma^{2j}) \cdots (\lambda.\sigma^{(k-1)j})
\]
and we see that
\begin{align*}
 \mu\, (\mu^{-1}.\sigma^j) &= \left(\lambda \, (\lambda.\sigma^j) \, (\lambda.\sigma^{2j}) \cdots (\lambda.\sigma^{(k-1)j})\right) \left( (\lambda^{-1}.\sigma^{j}) \, (\lambda^{-1}.\sigma^{2j}) \, (\lambda^{-1}.\sigma^{3j}) \cdots (\lambda^{-1}.\sigma^{kj})\right) \\
 &= \lambda \, (\lambda^{-1}.\sigma^{kj}) = \lambda\, (\lambda^{-1}.\sigma).
\end{align*}
This shows that $[H,\sigma] \subseteq [H,\sigma^j]$. It follows that $[H,\sigma^i] = [H,\sigma^j]$ for all $i$ and $j$. Hence the hypothesis of Corollary \ref{cor:trivial-fiber-action} holds. 
\end{proof}

The next two examples show that some first guesses as to how to generalize Proposition~\ref{prop:prime} don't work.
\begin{example} The hypothesis of Corollary~\ref{cor:trivial-fiber-action} does not automatically hold when $S$ is cyclic of composite order.
Let $W=x_1^4+x_2^4+x_3^4+x_4^4$, let $S = \left< (1234) \right>$ and let
\[
H = \left< (1/4,1/4,0,0), (0,1/4,1/4,0), (0,0,1/4,1/4)\right>.
\]
Then consider $\mu:=(1/4,1/4,0,0) - (0,1/4,1/4,0) = (1/4,0,-1/4,0) \in [H,\sigma] \cap \ker \beta_{\sigma^2}$. However, we see that $[H, \sigma^2]$ is generated by $(1/4,1/4,-1/4,-1/4)$ and $(0,1/2,0,-1/2)$, so it does not contain $\mu$. 

If we look a little more closely, we can consider $W_{\sigma^2}=y_1^4+y_2^4$, and take $c=\lfloor 1; (\tfrac{1}{2}, \tfrac{1}{2})\rceil$. Then we have $c'= \lfloor y_1y_2; (0,0)\rceil$. One can check that $\Ea^T_{c'}$ is not empty, but $\Ba_c$ is also not empty, hence the DSC \ref{PC} also does not hold. 
\end{example}

%\begin{example} The hypothesis of Corollary~\ref{cor:trivial-fiber-action} does not automatically hold when $S$ is generated by disjoint cycles of prime order.
%Let $S = \left<(12),(34)\right>$ and $H = \left<(-1/4,1/4,-1/4,1/4), (1/4,-1/4, -1/4, 1/4) \right>$. Then $[H,(12)] = \left<(1/2,1/2,0,0) \right>$ and $[H,(12)(34)] = \left< (1/2,1/2,1/2,1/2) \right>$, and $\ker \beta_{(12)(34)} = H$. Hence we have $(1/2,1/2,0,0) \in [H,(12)] \cap \ker \beta_{(12)(34)}$, but $(1/2,1/2,0,0) \notin[H,(12)(34)]$.
%
%\end{example}

\section{Examples} \label{sec:examples}
The examples in this section show why our hypotheses are necessary.
\subsection{\texorpdfstring{$S$}{S} acts with determinant \texorpdfstring{$-1$}{-1}}
Let $W=x^3+y^3$, $S = \left<(12)\right>$, and $H = G_W^{diag}$, so $H^T$ is trivial. On the transpose side, we have elements
\[
    b' =\lfloor dz_1\wedge dz_2; \text{id}\rceil.
\]
Then $b'.(12) = - b'$. Hence $b'$ cannot contribute to the projected state space.  For this element, the mirror map is the standard Krawitz map that matches this with
\[
    b = \lfloor 1; (\frac13, \frac13) \rceil,
\]
which is invariant under $S$ (and all of $G$).

\subsection{\texorpdfstring{$S$}{S} has nondiagonal elements that rescale basis elements} 
Let $W = \sum_{i=1}^6 x_i^3$, so $j_W = (\frac13,\frac13,\frac13,\frac13,\frac13,\frac13)$. Let $H = \left<j_W\right>$ and $S = \left<(123),(456)\right>$. Every element of $S$ has odd order, so we don't have problems with signs as in the previous section.

Consider the element
\[
    b_1:= \lfloor 1; (123)(456) j_W \rceil
\]
Recall that there is implicitly a $dz_1\wedge dz_4$. One can check using Corollary~\ref{cor:act-z} that $z_1.(123) = \zeta^{-1} z_1$, where $\zeta = \exp(2\pi i / 3)$. On the other hand, $z_4.(123) = z_4$. So we see that $b_1.(123) = \zeta^{-1} b_1$. But $b_1$ is invariant under $H$, so the DSC \ref{PC} fails. Let's see what goes wrong. More generally, one checks that 
\[
    b_i:= \lfloor 1; (123)(456) j_W^i \rceil.
\]
for $i=0,1,2$ satifies $b_i.(123) = \zeta^{-i} b_i$. 

Write $\sigma = (123)(456)$ and
\[
    c:= \lfloor 1; (0,0) \rceil \in \Ca_\sigma.
\]
Notice $b_i \in \Ba_c$ for each $i$. Equation \eqref{eq:num-orbits} shows that the number of $H$ orbits in  $\Ba_c$ is the index of $[H,\sigma]$ in $\ker \beta_\sigma$. In our example, $[H,\sigma]$ is trivial, while $\ker \beta_\sigma = H$ is of order 3. Each of the $b_i$ are invariant under $H$, so these are the three orbits in $\Ba_c$. However, only $b_0 \in \Ra$ and gives rise to an element of the projected state space.

On the other hand, the Krawitz mirror map matches $c$ with
\[
    c' = \lfloor 1; (1/3,1/3) \rceil \in \Ca^T_\sigma.
\]
 One can see that every element of $\Ba_{c'}^T$ is narrow, that is,  its group element has trivial fixed locus. It follows that every element of $\Ba_{c'}^T$ is in $\Ra$. By Lemma \ref{lem:count-fibers}, we know that $\Ba_{c'}^T$ has three $H^T$ orbits. We claim that $S$ acts trivially on the set of orbits in $\Ba_{c'}^T$. Hence, $\Ba_{c'}^T$ contains three $G^T$ invariants, in contrast to the one $G$ invariant in $\Ba_c$.

To see the claim about the $S$ action, one can check that $H^T$ consists of precisely those elements $\mu$ with $\sum_i \mu_i \equiv 0 \mod \ZZ$.  Then we define a function $I: H^T \rightarrow \QQ/\ZZ$ by  $I(\mu) := \sum_i i \mu_i \mod \ZZ$. One can check that $I$ is $\sigma$ invariant---indeed,
\[
 I(\mu.\sigma) = 3 \mu_1 + \mu_2 + 2\mu_3 + 6 \mu_4 + 4\mu_5 + 5\mu_6 = I(\mu) - \sum_i \mu_i + 3 \mu_1 + 3 \mu_4,
\]
and we know $\mu_1,\mu_4 \in \frac13 \ZZ$ and $\sum_{i=1}^6 \mu_i \equiv 0 \mod \ZZ$.
It follows then that $I(\mu)= 0$ for any $\mu \in [H^T, \sigma]$. We also apply $I$ to $\Ba_{c'}^T$ by setting $I(\lfloor m; \sigma \mu \rceil) = I(\mu)$. Recall that (the sectors of elements in) an $H^T$ orbit in $\Ba_{c'}^T$ corresponds to a coset of $[H^T,\sigma]$, so $I$ is constant on  $H^T$ orbits. The only possible values of $I$ are $0, 1/3, 2/3$. One can find elements for each of these values, and hence $I$ actually classifies $H^T$ orbits in $\Ba^T_{c'}$. Since $I$ is $\sigma$ invariant, we see that the $\sigma$ action on $H^T$ orbits is trivial. 

\subsection{The action on \texorpdfstring{$H$}{H} orbits is nontrivial, but still compatible} \label{sec:good-act-H-orbits}

Although the hypothesis of Theorem \ref{thm:main-cyc} implies that the $S$ action on $H$ orbits will be trivial, our more general conditions allow for examples where the action is nontrivial, and the isomorphism still holds. We explain in detail for one pair $(c,c')$.

Let 
\[
    W = x_1^6 + x_2^6 + x_3^6 + y_1^6 + y_2^6 + y_3^6 + w_1^6 + w_2^6 + w_3^6.
\]
Now, fix permutations $\sigma = (123)$ and $\tau = (xyw)$. That is, $x_1.\sigma = x_3$, $y_1.\sigma = y_3$, $x_i.\tau = w_i$, etc. Let $S= \left<\sigma, \tau\right>$ and 
\[
H =\left< \alpha_x=(1/3,1/3,1/3,0,0,0,0,0,0), \alpha_y=(0,0,0,1/3,1/3,1/3,0,0,0), \alpha_w=(0,0,0,0,0,0,1/3,1/3,1/3) \right>.
\]
We see $\tau$ and $\sigma$ commute, $H$ commutes with $\sigma$, and $\tau^{-1}\alpha_u\tau = \alpha_{\tau^{-1}(u)}$, for $u=x,y,w$.

Consider
\[
    c:=\lfloor z_x^2 z_y^2 z_w^2; \text{id} \rceil \in \Ca_\sigma.
\]
Here, we write $z_x$ for the coordinate of the fixed locus corresponding to $x_1$. 
Then we see that $\Ba_c$ is of size $27$, consisting of 
\[
    \{ \lfloor z_x^2 z_y^2 z_w^2; \sigma \lambda \rceil : \lambda \in H \}.
\]
$\Ba_c$ consists of $27$ singleton $H$-$\star$ orbits, $\sigma$ acts trivially on these orbits, and $\tau$ acts on these by conjugating the $\lambda$. Hence there are 10 $\tau$ orbits of $H$-$\star$ orbits in $\Ba_c$, represented by $k_x \alpha_x + k_y \alpha_y + k_w \alpha_w$ with $0 \le k_x \le k_y \le k_w  \le 2$.

The element $c$ is a Krawitz mirror to
\[
 c':= \lfloor 1 ; (1/2,1/2,1/2) \rceil.
\]
One can check that $H^T$ is precisely those elements $\lambda$ so that $\lambda_1 + \lambda_2 + \lambda_3 \in \frac12 \ZZ$, $\lambda_4 + \lambda_5 + \lambda_6 \in \frac12 \ZZ$, and $\lambda_7 + \lambda_8+ \lambda_9 \in \frac12 \ZZ$.  We know by Lemma~\ref{lem:count-fibers} that there are $27$ $H^T$-$\star$ orbits in $\Ba_{c'}$. Let's try to classify them. Since every element of $\Ba_{c'}^T$ is narrow, we need only consider the $H^T$ action by conjugation on $\sigma \beta_\sigma^{-1}((1/2,1/2,1/2))$. Each $H^T$ orbit is a $[H^T, \sigma]$ coset. We claim that these orbits are classified by the function
\[
I(\lambda) = (\lambda_1 + 2\lambda_2 + 3\lambda_3, \lambda_4 + 2\lambda_5 + 3\lambda_6, \lambda_8 + 2\lambda_8 + 3\lambda_9) \mod \frac12 \ZZ.
\]
Indeed, one can check that $I$ is invariant when one adds an element of $[H^T, \sigma]$. Furthermore, for integers $0 \le k_x, k_y, k_w  \le 2$, set
\[
    \lambda = (1/2, -k_x/6, k_x/6, 1/2, -k_y/6, k_y/6, 1/2, -k_w/6, k_w/6),
\]
and note that $\lambda \in \beta_\sigma^{-1}((1/2,1/2,1/2)) \cap H^T$, and $I(\lambda) = (k_x/6,k_y/6,k_w/6)$. This accounts for all 27 orbits, so $I$ classifies them as we claimed. Finally, observe that $I(\lambda.\tau)$ is the same as $I(\lambda)$ with the entries permuted by $\tau$. Hence the action corresponds exactly to the action on the $H$ oribts in $\Ba_{c}^T$.

\subsection{The action on \texorpdfstring{$H$-$\star$}{H star} orbits is nontrivial, but the bijection fails} \label{sec:bad-act-H-orbits} Unfortunately, it doesn't always work as nicely as in the previous example. 

%Idea: Let $S'$ be the subgroup of $S$ that preserves a fiber of $\phi$. Then $S'$ acts on the $H$ orbits within that fiber. We wish that the set of $H$-$\star$ orbits with this action is isomorphic (as an $S'$-set) to that on the other side. But this example show that this is not always the case.

Let 
\[
    W = x_1^9 + x_2^9 + x_3^9 + y_1^9 + y_2^9 + y_3^9 + w_1^9 + w_2^9 + w_3^9.
\]

Let $\lambda_x = (1/9,1/9,1/9,4/9,4/9,4/9,7/9,7/9,7/9)$, $\lambda_y = (7/9,7/9,7/9,1/9,1/9,1/9,4/9,4/9,4/9)$, and $\lambda_w= (4/9,4/9,4/9,7/9,7/9,7/9,1/9,1/9,1/9)$. Note $\lambda_y = \lambda_x^4$ and $\lambda_w = \lambda_x^7$. Set $H = \left< \lambda_x\right>$. 

Now, fix permutations $\sigma = (123)$ and $\tau = (xyw)$. That is, $x_1.\sigma = x_3$, $y_1.\sigma = y_3$, $x_i.\tau = w_i$, etc. Let $S = \left<\tau, \sigma\right>$. We have that $\sigma$ commutes with everything, and $\tau^{-1} \lambda_u \tau = \lambda_{\tau^{-1}(u)}$, for $u=x,y,w$.

Consider the elements
\[
 b_u:=\lfloor 1; \sigma \lambda_u \rceil.
\]
Since there is no fixed locus here, of course there is no $g \in G$ with $b.g = \alpha g$, $\alpha \neq 1$. We see that $b_u \star \tau = b_{\tau^{-1}(u)}$ and $b$ is fixed by $H$ and $\sigma$.

Notice that the $b_u$ are all in $\Ba_c$, with
\[
c:=\lfloor 1; (1/3,1/3,1/3) \rceil.
\]
We have $[H,\sigma]$ is trivial but $\ker \beta_\sigma$ is not trivial---it is generated by $\lambda_x^3 = (1/3,1/3,1/3,1/3,1/3,1/3,1/3,1/3,1/3)$. Hence, by \eqref{eq:num-orbits}, $\Ba_c$ has 3 $H$-$\star$ orbits, which are the three singletons $\{b_x\},\{b_y\},\{b_w\}$. They are fixed by $\sigma$ and permuted by $\tau$.

The $c$ is  Krawitz dual to
\[
   c':=\lfloor z_x^2 z_y^2 z_w^2 ;\text{id} \rceil.
\]
In $\Ba_{c'}^T$, by Lemma~\ref{lem:count-fibers}, there are 3 $H^T$-$\star$ orbits.  What is the $S$ action on these? We have for example:
\[
 b':=\lfloor z_x^2 z_y^2 z_w^2 ; \sigma  \rceil \in \Ba_{c'}^T.
\]
Here, we write $z_u = c(u_1 + u_2 + u_3)$, so $z_u.\sigma = z_u$, and $z_u.\tau = z_{\tau^{-1}(u)}$. So we see that $b'$ is fixed by $S$. (Notice the sign after permuting the volume form is 1.) Hence the $H$-$\star$ orbit of $b'$ is fixed by $S$. An action of $\tau$ on a set of size 3 that has a fixed point is of course trivial. 

Hence, no matter how we match up the $H$-$\star$ orbits in $\Ba_c$ with the $H^T$-$\star$ orbits in $\Ba_{c'}$, we will not get a bijection of $G$-$\star$ orbits with $G^T$-$\star$ orbits.

%Why is this example a problem?
%We are trying to match up $G$ orbits on the B side with $G^T$ orbits on the A side. Here, on the first side there is a nontrivial action of $\sigma$ on the three $H$ orbits in the fiber. So in this fiber there is actually only 1 $G$ orbit. But on the other side, the action of $G^T$ is trivial. So there are 3 $G^T$ orbits.

%In this example, the whole $G$ or $G^T$ orbits are contained in a single $\Ba_c$ or $\Ba_{c'}^T$ (that may not be true in general). In that case, our proof strategy wants us to match up $G$-$\star$ orbits with $G^T$of matching up fibers, in this case, we would have expected that we could match up orbits between the corresponding fibers.

%Furthermore, I believe (thought I should check this more carefully) that the other fibers (besides ones that do this) match up nicely. So in fact, the desired theorem is not true (and not just that our proof approach is wrong).

\bibliographystyle{abbrv}
\bibliography{references}

\end{document}